\documentclass[letterpaper,12pt,unknownkeysallowed, dvipsnames]{article} 
    \usepackage{algorithm}
    \usepackage[noend]{algpseudocode}
    \usepackage{booktabs}
    \usepackage{filemod}
	\usepackage[utf8]{inputenc}
	\usepackage[english]{babel}
	\usepackage{amsfonts}
	\usepackage{amssymb}
	\usepackage{enumitem}
	\usepackage{bbm}
	\usepackage{fullpage}
	\usepackage{amsmath}
	\usepackage{mathabx}
	\usepackage{mathtools}
    \usepackage{scalerel}
    \usepackage{xcolor}
    \usepackage{array}
    \usepackage{tablefootnote}
    \usepackage{threeparttable}
    \newcolumntype{L}[1]{>{\raggedright\let\newline\\\arraybackslash\hspace{0pt}}m{#1}} 
    \usepackage{multirow}

	\usepackage{natbib}
\renewcommand\bibpreamble{\vspace{.8\baselineskip}}

\renewcommand{\appendix}{\small\parindent 0cm\parskip 5pt\setcounter{equation}{0}
\setcounter{section}{0}
\renewcommand{\thesection}{A.\arabic{section}}
\renewcommand{\theequation}{A.\arabic{equation}}
\setcounter{lemma}{0}\renewcommand{\thelemma}{A.\arabic{lemma}}
\setcounter{theorem}{0}\renewcommand{\thetheorem}{A.\arabic{theorem}}
}

\makeatletter
\def\blfootnote{\xdef\@thefnmark{}\@footnotetext}
\makeatother

\makeatletter
\def\@sect#1#2#3#4#5#6[#7]#8{\ifnum #2>\c@secnumdepth
     \let\@svsec\@empty\else
     \refstepcounter{#1}\edef\@svsec{\csname the#1\endcsname. \hskip 0.4em}\fi
     \@tempskipa #5\relax
      \ifdim \@tempskipa>\z@
        \begingroup #6\relax
          \@hangfrom{\hskip #3\relax\@svsec}{\interlinepenalty \@M #8\par}
        \endgroup
       \csname #1mark\endcsname{#7}\addcontentsline
         {toc}{#1}{\ifnum #2>\c@secnumdepth \else
                      \protect\numberline{\csname the#1\endcsname}\fi
                    #7}\else
        \def\@svsechd{#6\hskip #3\relax  
                   \@svsec #8\csname #1mark\endcsname
                      {#7}\addcontentsline
                           {toc}{#1}{\ifnum #2>\c@secnumdepth \else
                             \protect\numberline{\csname the#1\endcsname}\fi
                       #7}}\fi
     \@xsect{#5}}
\makeatother

\makeatletter
\renewcommand{\section}{\@startsection{section}{1}{0mm}{-\baselineskip}{0.25\baselineskip}{\centering\normalfont\normalsize\bf}}
\renewcommand{\subsection}{\@startsection{subsection}{2}{0mm}{-\baselineskip}{0.25\baselineskip}{\raggedright\normalfont\normalsize\bf}}
\renewcommand{\subsubsection}{\@startsection{subsubsection}{3}{0mm}{-\baselineskip}{0.25\baselineskip}{\raggedright\normalfont\small}}
\def\@begintheorem#1#2{\trivlist \item[\hskip \labelsep{\bf #1\ #2}]\it}
\makeatother



\renewcommand{\thesection}{\arabic{section}}
\renewcommand{\thesubsection}{\arabic{section}.\arabic{subsection}}

 \newcommand*\mwidebar[1]{%
   \hbox{%
     \vbox{%
       \hrule height 0.5pt %
       \kern0.30ex
       \hbox{%
         \kern-0.02em
         \ensuremath{#1}%
         \kern-0.02em
       }%
     }%
   }%
} 

\newcommand*\lowwidebar[1]{%
   \hbox{%
     \vbox{%
       \hrule height 0.5pt %
       \kern0.20ex
       \hbox{%
         \kern-0.02em
         \ensuremath{#1}%
         \kern-0.02em
       }%
     }%
   }%
} 

\allowdisplaybreaks

\begin{document}
\vspace*{0.2cm}

\thispagestyle{empty}
\vskip 20pt
\centerline{\large\bf When Should You Adjust Standard Errors for Clustering?}
  \begin{center}
  
    \vskip 10pt
    \blfootnote{\hspace*{-0.25in}\!The questions addressed in this article partly originated in discussions with Gary Chamberlain. We are grateful for questions raised by Chris Blattman and seminar audiences, and for insightful comments by Colin Cameron,  Vicente Guerra, four reviewers, Larry Katz, and Jesse Shapiro. Jaume Vives-i-Bastida provided expert research assistance. This work was supported by the Office of Naval Research under grants N00014-17-1-2131 and N00014-19-1-2468.}
    
    {\large
     \lineskip .5em
   \begin{tabular}[t]{ccc}
       Alberto Abadie&&Susan Athey\\[.2ex]
       MIT&&Stanford\\[1ex]
       Guido W. Imbens&&Jeffrey M. Wooldridge\\[.2ex]
       Stanford&&MSU\\
    \end{tabular}
      \par}
      \vskip 1em
      {\large \today } \par
       \vskip 1em
  \end{center}\par

\bigskip
\begin{center}\normalsize\bf\text{Abstract}
 \end{center}\begin{quote}\normalsize
\noindent  
Clustered standard errors, with clusters defined by factors such as geography, are widespread in empirical research in economics and many other disciplines.
Formally, clustered standard errors adjust for the correlations induced by sampling the outcome variable from a data-generating process with unobserved cluster-level components. However, the standard econometric framework for clustering leaves important questions unanswered: {\em (i)} Why do we adjust standard errors for clustering in some ways but not others, e.g., by state but not by gender, and in observational studies, but not in completely randomized experiments? {\em (ii)} Why is conventional clustering an ``all-or-nothing'' adjustment, while within-cluster correlations can be strong or extremely weak? {\em (iii)} In what settings does the choice of whether and how to cluster make a difference? We address these and other questions using a novel framework for  clustered inference on average treatment effects. In addition to the common sampling component, the new framework incorporates a design component that accounts for the variability induced on the estimator by the treatment assignment mechanism.
We show that, when the number of clusters in the sample is a non-negligible fraction of the number of clusters in the population, conventional cluster standard errors can be severely inflated, and propose new variance estimators that correct for this bias.
\end{quote}

\newpage
\addtolength{\baselineskip}{0.5\baselineskip}
 
\section{Introduction} \label{sec:introduction}
\medskip

Imagine you estimated the effect of attending college on labor earnings using linear regression on a cross-section of U.S. workers. How should you calculate the standard error? Empirical studies in economics often report heteroskedasticity-robust standard errors (henceforth ``robust'') associated with the work by \cite{eicker1963}, \cite{huber1967behavior}, and \cite{white1980heteroskedasticity}. A common alternative is to report cluster-robust standard errors (henceforth ``cluster'') associated with the work by \cite{liang1986longitudinal} and \cite{arellano1987practitioners}, with clustering often applied within geographic units such as states or counties. \citet{moulton1986random, moulton1987diagnostics} and \citet*{Bertrand2004did} have shown that clustering adjustments can make a substantial difference, and since the 1980s cluster standard errors have become commonplace in empirical economics.   

Later in this section, we estimate a log-linear regression of earnings on an indicator for some college using data from the 2000 U.S.\ Census. We find that standard errors clustered at the state level are more than 20 times larger than robust standard errors. Which ones should a researcher report? The conventional framework for clustering \citep*[see][for recent reviews]{cameron2015practitioner, mackinnon2021cluster} 
suggests that if the clustering adjustment matters, in the sense that the cluster standard errors are substantially larger than the robust standard errors, one should use the cluster standard errors. In this article, we develop a new framework for cluster adjustments to standard errors that nests the conventional framework as a limiting case. The new framework suggests novel standard error formulas that can substantially improve over robust and cluster standard errors in settings like the earnings regression described above. 

Our proposed clustering framework differs from the standard one in that it includes a design component that accounts for between-clusters variation in treatment assignments. We argue that the new design component is important because between-cluster variation in treatment assignments often motivates the use of clustered standard errors in empirical studies \citep[see, e.g.,][]{gentzkow2008preschool,cohen2010free}. In addition, our framework shifts the focus of interest from features of infinite super-populations/data-generating processes to average treatment effects defined for the finite (but potentially large) population at hand. As a result of this shift, it is the sampling process and the treatment assignment mechanism that solely determine the correct level of clustering; the presence of cluster-level unobserved components of the outcome variable becomes irrelevant for the choice of clustering level. Moreover, by focusing on finite populations (which could be entirely or substantially sampled in the data) we obtain standard errors smaller than those aiming to measure uncertainty with respect to features of infinite super-populations. We derive the large sample variances for the least squares and fixed effect estimators under our proposed framework and show that they differ in general from both the robust and the cluster variances. We also propose two estimators for the large sample variances, one analytic and one based on a re-sampling (bootstrap) approach. For the U.S. earnings application, our proposals produce standard errors that are substantially larger than the robust standard errors, but also substantially smaller than the conventional version of cluster standard errors.

We use our framework to highlight three common misconceptions surrounding clustering adjustments. The first misconception is that the need for clustering hinges on the presence of a non-zero correlation between residuals for units belonging to the same cluster. We show that  the presence of such correlation does not imply the need to use cluster adjustments, and that the absence of such correlation does not imply that clustering is not required. The second misconception is that there is no harm in using clustering adjustments when they are not required, with the implication that if clustering the standard errors makes a difference, one should do so. To see that both of these claims are in fact incorrect, consider the following simple example. Suppose that, based on a random sample from the population of interest, we use the sample average of a variable to estimate its population mean. Suppose also that the population can be partitioned into clusters such as geographical units. If outcomes are positively correlated within clusters, the cluster variance will be larger than the robust variance. However, standard sampling theory directly implies that if the units are sampled randomly from the population there is no need to cluster. The harm in clustering in this case is that confidence intervals will be unnecessarily conservative, possibly by a wide margin. A third misconception is that researchers have only two choices: either fully adjust for clustering and use the cluster standard errors, or not adjust the standard errors at all and use the robust standard errors. 
We show that a combination of the robust and the cluster variance estimators can substantially improve accuracy over its two components.

The new clustering framework in this article has the advantage of providing actionable guidance on a question of substantial consequence for empirical practice in econometrics: When should standard errors be clustered, and at what level? In the conventional model-based econometric framework, the researcher takes a stand on the error component structure of a model for the outcome variable. For example, suppose that, following \citet{moulton1986random, moulton1987diagnostics}, the researcher posits a random effects model, with random effects at the state level. In this setting, a repeated sampling thought experiment entails that, for each sample, different values of the state random effects are drawn from their distributions. This model-based approach implies that if we are estimating a population mean using a sample average one needs to cluster the standard errors at the state level {\em even if the sample is a random sample of individuals} and not a clustered sample. A drawback of the model-based econometric framework for clustering is that empirical researchers 
need to take a stand on the structure of the error components of their models.

A second, closely related, framework for clustering that is often invoked in the econometrics literature is motivated by a sampling mechanism that in a first stage selects clusters at random from an infinite population, followed by a second stage of random sampling of units from the sampled clusters (or keeping all units in a cluster). Although this framework is appropriate for some applications in the analyses of surveys, where it originated \citep{kish1995survey, thompson2012sampling}, we argue that it is not appropriate for many of the data sets economists and other social scientists analyze. In many applications in economics, researchers do observe units from {\it all} the clusters they are interested in, e.g., all the states in the U.S., and a framework based on randomly sampling a small fraction of a large population of clusters does not apply.  

Neither of the two conventional frameworks for clustered inference described above fully incorporates the design aspect of clustering. And it is the lack of a design component that makes them inappropriate for inference on treatment effects. To gain insight on the importance of the assignment mechanism for the standard errors of treatment effects estimators, consider a setting with individuals sampled at random from a population, but where treatment is assigned at the cluster level, with the same treatment value for all the individuals in the same cluster. Assume that the quantity of interest is the population average treatment effect. Clustered assignment to treatment is equivalent to clustered sampling of potential outcomes. Because the parameter of interest depends on averages of potential outcomes, which are sampled in a clustered manner, clustering of the standard errors is required in this setting, even when the individual observations are sampled at random. Our framework for clustered inference in this setting is close in spirit to the sampling framework described in the previous paragraph, but it incorporates a design component.  

By shifting the attention from parameters of a data generating process for the outcomes to the average treatment effect for the population at hand, a researcher applying the proposals in this article does not need to take a stand on the error component structure of a model for the outcome variable to calculate standard errors. Instead, all the relevant variability of the estimator with respect to the average treatment effect is generated by the sampling mechanism, which extracts the sample from the population, and the assignment mechanism, which determines which units are exposed to the treatment. We see this as an intrinsic advantage of the framework proposed in this article in settings where it is difficult  to justify a particular error component structure.


In this article we make three contributions. The first one is a novel framework for clustering, building on the one developed by \citet{abadie2020sampling} for the analysis of regression estimators from a design perspective. We allow for clustering both in the sampling process and in the assignment process. As a result, the framework nests both the traditional case of clustered sampling and the case of clustered treatment assignment in experiments as special cases. It also allows for intermediate cases. In particular, treatment assignment may depend on cluster but not perfectly so, and there remains variation in treatments within-clusters. This framework clarifies the separate roles of clustering in the sampling process and clustering in the assignment process. It also clarifies what we can learn from the data about the need to adjust standard errors for clustering. In our framework, the data are \emph{not} informative about the need to adjust for clustering in the sampling process, but they \emph{are} informative about the need to adjust for clustering in the assignment process.

In our second  contribution, we derive central limit theorems and large sample variances for the least squares and the fixed effect estimators of average treatment effects that take into account variation both from sampling and assignment. Comparing these variances to limit versions of the robust and cluster variances shows that the robust standard errors are generally too small, and the cluster standard errors are unnecessarily conservative. These comparisons also highlight how heterogeneity in treatment effects affects inference in the estimation of average treatment effects. Often researchers specify models that implicitly assume constant treatment effects without appreciating the implications for inference. We show, however, that heterogeneity in treatment effects introduces additional variance components that affect the need for clustering adjustments. 

In our third contribution, we propose new variance formulas and bootstrap procedures for treatment effects estimators in the presence of clustering.  We use the term Causal Cluster Variance (CCV) for the analytic variance formulas. For the case of a least squares estimator of average treatment effects, the intuition for the CCV variance formula is as follows. The error of the least squares estimator is approximately equal to a sum, over all units, of an expression involving products of regression errors and regressors values. The robust variance is approximately equal to a sum, over all units, of the squares of these products. In contrast, the conventional cluster variance estimator is approximately equal to a sum, over all clusters, of squares of within-cluster sums of the same products. Although the sum over all clusters of the expectation of the within-cluster sums of these products is zero, the expectation for each cluster separately is not. For each cluster in the sample, it is possible to estimate the expectation of the sum of the products between regression errors and regressors values. The CCV formula uses these estimates to correct the bias of the conventional cluster variance. The CCV correction does not help much if only a small fraction of clusters are sampled. However, when a large fraction of the clusters are represented in the sample, the CCV correction can lead to substantial improvements. 
This adjustment relies on estimates of cluster-level treatment effects, and thus requires within-cluster variation in treatment assignment. 
In addition, we propose a bootstrap version of the variance estimator., which we compare to two benchmarks. In contrast to conventional bootstrap procedures, which are based on resampling individual units or entire clusters of units, our proposed Two-Stage-Cluster-Bootstrap (TSCB) conducts resampling in two stages. In the first stage, the fraction treated for each cluster is drawn from the empirical distribution of cluster-specific treatment fractions. In the second stage, the researcher samples the treated and control units from each cluster, with their number of units determined in the first stage. The CCV and TSCB variance estimators are designed for applications with large number of observations and substantial variation in treatment assignment within clusters. 

To illustrate the empirical relevance of our results, we  analyze a sample from the  2000 U.S. Decennial Census, which includes 2,632,838 individuals. We define 52 clusters according to residency in the 50 states, Puerto Rico, and the District of Columbia.
We consider two log-linear regressions of individual earnings on a treatment variable that encodes information on college attendance. In the first specification, the treatment variable is measured as an average, at the state level. In a second specification, we measure college attendance at the individual level. 

\begin{table}[t]
    \centering
    \caption{College effects in the Census sample}
    \vspace*{0.4cm}
    \begin{tabular}{lcccc}
    \multicolumn{5}{l}{{\it Dependent variable:} Log labor earnings}\\
    \hline
   
   \multicolumn{5}{l}{\it Panel A}\\
    
    \multicolumn{1}{l}{\hspace{0.2cm}\it Treatment:}
    &\multicolumn{4}{l}{\hspace{-4.8cm}State indicator for share of some}\\
    
    &\multicolumn{4}{l}{\hspace{-4.8cm}college greater than 0.55}\\[1.5ex]
         &&OLS&&\\[.5ex] 
         coefficient&&0.1022\\
         standard error:\\
         \ \,robust&&(0.0012)\\
         \ \,cluster&&(0.0312)\\
    \hline
    \multicolumn{5}{l}{{\it Panel B}}\\
\multicolumn{1}{l}{\hspace{0.2cm}\it Treatment:}
    &\multicolumn{4}{l}{\hspace{-4.8cm}Individual indicator for some college}\\[1.5ex]
         &&OLS&&FE\\ 
         coefficient&&0.4656&&0.4570\\
         standard error:\\
         \ \,robust&&(0.0012)&&(0.0012)\\
         \ \,cluster&&(0.0269)&&(0.0276)\\
         \ \,causal cluster variance (CCV)&&
         (0.0035)
         &&(0.0014)\\
         \ \,two-stage cluster bootstrap (TSCB)&&(0.0036)&&(0.0014)\\[.5ex]
    \hline
    \end{tabular}
    \label{intro_table}
\end{table}

In Panel A of Table \ref{intro_table}, we report results for a regression where the only explanatory variable is a binary treatment that takes value one if the fraction of individuals with at least some college residing in the state is 0.55 or higher, and value zero otherwise (we chose the 0.55 value to ensure sufficient variation in the treatment over the 52 clusters). Notice that  the treatment is constant within states. We report the ordinary least squares (OLS) estimate, as well as robust and cluster standard errors. Since the late 1980s, it has been common practice to report cluster standard errors in settings where the regressors are constant within a cluster. Clustering at the state level makes a substantial difference relative to using robust standard errors, with the cluster standard errors approximately twenty-six times larger than the robust standard errors. 

In Panel B of Table \ref{intro_table}, the sole regressor is an individual-level indicator for at least some college. In addition to OLS, we report the fixed effects (FE) estimate (with fixed effects for the 50 states, plus Washington DC and Puerto Rico) and robust, cluster,  CCV, and TSCB standard errors in parentheses. Like for the regression of the first panel, clustering at the state level makes a substantial difference in the standard errors, with the cluster standard errors approximately twenty-three times larger than the robust standard errors, both for the OLS and the FE regressions. In Panel B, our proposed CCV and TSCB  standard errors for the OLS estimate are  0.0035 and 0.0036 respectively, in between the robust standard errors (0.0012) and the cluster standard errors (0.0269), and substantially different from both. The same holds for the FE estimator. The cluster standard error is 0.0276,  quite different from the robust standard errors, 0.0012. The CCV and TSCB standard errors are 0.0014, in between robust and cluster but much closer to robust.

\section{A Framework for Clustering}
\label{section:set_up}

In this section, we describe in detail the framework for our analysis. There are multiple components to our set-up that are not explicitly modeled in the usual analysis of the variance of econometric estimators. In general, quantifying the uncertainty of parameter estimates requires describing the population and articulating the assumptions that describe how the sample was generated from that population (that is, building a model for the data generating process). In our framework, there are three distinct sources of sampling variation that lead to variation in the estimates. First, there is variation across samples in which units are observed in each cluster. Second, there is potentially variation in which clusters are observed (which leads to different units being observed). Third, there is variation in the treatment assignment across units. Whereas the standard framework for clustering focuses solely on the first two (sampling)
sources of uncertainty, our proposed framework allows for all three. How much these three components matter for the variance of the least squares and fixed effects estimators of the average treatment effect depends on {\em (i)} the sampling process, {\em (ii)} the assignment process, and  {\em (iii)} the heterogeneity in the treatment effects across clusters. To facilitate the calculation of asymptotic approximations in a range of relevant settings for empirical practice, it is convenient to formally consider a sequence of populations where we can separately control the fraction of units in the population that are sampled and the fraction of clusters in the population that is sampled, as well as the assignment mechanism.

\subsection{A Sequence of Populations}

We have a sequence of populations indexed by $k$. The $k$-th population has $n_k$ units, indexed by $i=1,\ldots, n_k$. The population is partitioned into $m_k$  clusters. Let $m_{k,i}\in \{1,\ldots, m_k\}$ denote the cluster that unit $i$ of population $k$ belongs to. The number of units in cluster $m$ of population $k$ is $n_{k,m}\geq 1$. For each unit, $i$, there are two potential outcomes, $y_{k,i}(1)$ and $y_{k,i}(0)$, corresponding to treatment and no treatment.
Thus the population is characterized by the set of triples $(m_{k,i},y_{k,i}(0),y_{k,i}(1))$, for units $1,\ldots,n_k$ and clusters $1,\ldots, m_k$.
The object of interest is the population average treatment effect
\begin{equation*}
\label{estimand}
\tau_k = \frac{1}{n_k}\sum_{i=1}^{n_k}\bigl(y_{k,i}(1)-y_{k,i}(0)\bigr).
\end{equation*}
The population average treatment effect by cluster is
\[
\tau_{k,m} = \frac{1}{n_{k,m}}\sum_{i=1}^{n_k} 1\{m_{k,i}=m\}(y_{k,i}(1)-y_{k,i}(0)).
\]
Therefore,
\[
\tau_k = \sum_{m=1}^{m_k}\frac{n_{k,m}}{n_k}\tau_{k,m}.
\]
We assume that potential outcomes, $y_{k,i}(1)$ and $y_{k,i}(0)$, are bounded in absolute value, uniformly for all $(k,i)$.

For each unit in the population, we define the stochastic treatment indicator, $W_{k,i}\in\{0,1\}$. The realized outcome for unit $i$ in population $k$ is $Y_{k,i}=y_{k,i}(W_{k,i})$. For a random sample of  the population, we observe the triple $(Y_{k,i},W_{k,i},m_{k,i})$. Inclusion in the sample is represented by the random variable $R_{k,i}$, which takes value one if unit $i$ belongs to the sample, and value zero if not. We next describe the two components of the stochastic nature of the sample: the sampling process that determines the values of $R_{k,i}$, and the assignment process that determines the values of $W_{k,i}$. 
 
\subsection{The Sampling Process}
\label{section:sample}

The sampling process that determines the values of $R_{k,i}$ is independent of the potential outcomes and the assignments. It consists of two stages. First, clusters are sampled with cluster sampling probability $q_k\in (0,1]$. Second, units are sampled from the subpopulation consisting of all the sampled clusters, with unit sampling probability equal to $p_k\in (0,1]$. Both $q_k$ and $p_k$ may be equal to one, or close to zero. If $q_k=1$, we sample all clusters. If $p_k=1$, we sample all units from the sampled clusters. If $q_k=p_k=1$, all units in the population are sampled. The standard framework for analyzing clustering focuses on the special case where $q_k\rightarrow 0$, so only a small fraction of the clusters in the population are sampled. The case $q_k=1$ and $p_k\rightarrow 0$ corresponds to taking a relatively small random sample of units from the population. While this is an important special case, there are also many applications where the sampled clusters comprise a large fraction of the overall set of clusters. We refer to the case of $q_k=1$ as {\em random sampling} and to the case of $q_k<1$ as {\em clustered sampling}. 

\subsection{The Assignment Process}
\label{section:assignment_process}

The assignment process that determines the values of $W_{k,i}$ also consists of two stages. In the first stage of the assignment process, for cluster $m$ in population $k$, an assignment probability $A_{k,m}\in[0,1]$ is drawn randomly from a distribution with mean $\mu_k$, bounded away from zero and one uniformly in $k$, and variance $\sigma^2_k$, independently for each cluster. The variance $\sigma^2_k$ is key. If $\sigma^2_k$ is zero, then $A_{k,m}$ is the same for all $m$, and $W_{k,i}$ is randomly assigned across clusters. We refer to this case as {\em random assignment}. For positive values of $\sigma_k^2$ assignment probabilities depend on cluster. Because $A_{k,m}^2\leq A_{k,m}$, it follows that $\sigma_k^2$ is bounded above by $\mu_k(1-\mu_k)$ and that the bound is attained when $A_{k,m}$ can only take values zero or one, so all units within a cluster have the same values for the treatment. 
We use the term {\em clustered assignment} to refer to the case  $\sigma_k^2=\mu_k(1-\mu_k)$, when there is no within-cluster variation in $W_{k,i}$. We use the term {\em partially clustered assignment} to refer to the case  $0<\sigma_k^2<\mu_k(1-\mu_k)$, where assignment depends on cluster but not all units in the same cluster necessarily have the same value of $W_{k,i}$. In the second stage of the assignment process, each unit in cluster $m$ is assigned to the treatment independently, with cluster-specific probability $A_{k,m}$.

\section{The Least Squares Estimator and its Variance}
\label{section:least_squares_estimator}

Let 
\[
N_{k,1}=\sum_{i=1}^{n_k} R_{k,i}W_{k,i}\quad\mbox{ and }\quad N_{k,0}=\sum_{i=1}^{n_k} R_{k,i}(1-W_{k,i})
\]
be the number of treated and untreated units in the sample, respectively; these are random variables. The total sample size is $N_k=N_{k,1}+N_{k,0}$.

We first analyze the OLS estimator of a regression of the outcome $Y_{k,i}$ on an intercept and the treatment indicator $W_{k,i}$. The OLS estimator (modified so it is well-defined even when $N_{k,1}=0$ or $N_{k,0}=0$) is equal to the difference in means:
\begin{equation}
\label{ls_estimator}
\widehat\tau_k = \frac{1}{N_{k,1}\vee 1}\sum_{i=1}^{n_k} R_{k,i}W_{k,i} Y_{k,i}
-\frac{1}{N_{k,0}\vee 1}\sum_{i=1}^{n_k}R_{k,i}(1-W_{k,i}) Y_{k,i},
\end{equation}
where $N_{k,1}\vee 1$ and $N_{k,0}\vee 1$ are the maxima of $N_{k,1}$ and 1 and of $N_{k,0}$ and 1, respectively.

We make the following assumptions about the sampling process and the cluster sizes: {\em (i)} $m_kq_k\rightarrow \infty$, {\em (ii)} $\liminf_{k\rightarrow\infty}p_k \min_m n_{k,m}>0$, and {\em (iii)} $\limsup_{k\rightarrow\infty} \max_m n_{k,m}/\min_m n_{k,m}<\infty$. The first assumption implies that the expected number of sampled clusters goes to infinity as $k$ increases. The second assumption implies that the average number of observations sampled per cluster, conditional on the cluster being sampled, does not go to zero. The third assumption restricts the imbalance between the number of units across clusters. Notice that assumptions {\em (i)} and {\em (ii)} imply $n_kp_kq_k\rightarrow\infty$, so the sample size becomes larger in expectation as $k$ increases.

\subsection{Large $k$ Distribution of the Least Squares Estimator}

Our first main result derives the large $k$ distribution of $\widehat\tau_k$.
Let $\alpha_k = (1/n_k)\sum_{i=1}^{n_k} y_{k,i}(0)$, $u_{k,i}(1)=y_{k,i}(1)-(\alpha_k+\tau_k)$, and $u_{k,i}(0)=y_{k,i}(0)-\alpha_k$. Under additional regularity conditions in the Appendix,
\begin{equation*}
\sqrt{N_k}(\widehat\tau_k-\tau_{k})/v_k^{1/2}\stackrel{d}{\longrightarrow} N(0,1),
\end{equation*}
where
\begin{align}\label{vk}
v_k&=\frac{1}{n_{k}}\sum_{i=1}^{n_k} \bigg(\frac{u^2_{k,i}(1)}{\mu_k}+\frac{u^2_{k,i}(0)}{1-\mu_k}\bigg)\nonumber\\
&-p_k\frac{1}{n_{k}}\sum_{i=1}^{n_k}\big(u_{k,i}(1)-u_{k,i}(0)\big)^2-p_k\sigma_k^2
\frac{1}{n_{k}}\sum_{i=1}^{n_k}\bigg(\frac{u_{k,i}(1)}{\mu_k}+\frac{u_{k,i}(0)}{1-\mu_k}\bigg)^2\nonumber\\
&+p_k(1-q_k)\frac{1}{n_{k}}\sum_{m=1}^{m_k}\Bigg(\sum_{i=1}^{n_k} 1\{m_{k,i}=m\}\big(u_{k,i}(1)-u_{k,i}(0)\big)\Bigg)^2
\nonumber\\
&+p_k\sigma^2_k\frac{1}{n_{k}}\sum_{m=1}^{m_k}\Bigg(\sum_{i=1}^{n_k}1\{m_{k,i}=m\}\bigg(\frac{u_{k,i}(1)}{\mu_k}+\frac{u_{k,i}(0)}{1-\mu_k}\bigg)\Bigg)^2.
\end{align}
The expression for the variance $v_k$ has multiple terms that make its interpretation challenging. We first interpret $v_k$ in some special cases to highlight the implications of clustered sampling and clustered assignment. In Section \ref{section:discussion}, we compare $v_k$ to the large-$k$ form of the robust and cluster variance estimators.

For the case of random sampling ($q_k=1$) and random assignment ($\sigma_k^2=0$), the variance simplifies to
\[
\frac{1}{n_{k}}\sum_{i=1}^{n_k} \bigg(\frac{u^2_{k,i}(1)}{\mu_k}+\frac{u^2_{k,i}(0)}{1-\mu_k}\bigg)
-p_k\frac{1}{n_{k}}\sum_{i=1}^{n_k}\big(u_{k,i}(1)-u_{k,i}(0)\big)^2.
\]
As we show in Section \ref{section:robust_variances} below, the first term in this variance is estimated by the robust  variance estimator. The second term is a finite sample correction that is familiar from the literature on randomized experiments 
\citep[e.g.,][]{neyman1923, imbens2015causal, abadie2020sampling}. This finite sample correction vanishes if there is either no heterogeneity in the treatment effects (so $u_{k,i}(1)-u_{k,i}(0)=y_{k,i}(1)-y_{k,i}(0)-\tau_k=0$), or if the sample is a small fraction of the population ($p_k\approx 0$).

Adding clustered sampling, $q_k<1$, increases the variance by
\[ 
p_k(1-q_k)\frac{1}{n_{k}}\sum_{m=1}^{m_k}\Bigg(\sum_{i=1}^{n_k} 1\{m_{k,i}=m\}\big(u_{k,i}(1)-u_{k,i}(0)\big)\Bigg)^2,
\]
which is the same as
\[
p_k(1-q_k)\frac{1}{n_{k}}\sum_{m=1}^{m_k} n_{k,m}^2 (\tau_{k,m}-\tau_k)^2.
\]
This term vanishes if there is no heterogeneity in the average treatment effect across clusters. Although the sample is informative about heterogeneity in cluster average treatment effects, it is not informative about the value of $q_k$. Information about the need to adjust for clustered sampling ($q_k<1$) must come from outside the sample.

Clustered assignment, $\sigma^2_k>0$, adds two terms to the variance,
\[
-p_k\sigma_k^2
\frac{1}{n_{k}}\sum_{i=1}^{n_k}\bigg(\frac{u_{k,i}(1)}{\mu_k}+\frac{u_{k,i}(0)}{1-\mu_k}\bigg)^2\!
+p_k\sigma^2_k\frac{1}{n_{k}}\sum_{m=1}^{m_k}\Bigg(\sum_{i=1}^{n_k}1\{m_{k,i}=m\}\bigg(\frac{u_{k,i}(1)}{\mu_k}+\frac{u_{k,i}(0)}{1-\mu_k}\bigg)\Bigg)^2.
\]
As we explain in more detail in section \ref{section:discussion}, the sign of this expression depends on the amount of variation in potential outcomes that can be explained by the clusters. Note that in contrast to the lack of sample information about the need to adjust for clustered sampling, the sample is potentially informative about the need to account for clustered assignment.

The five terms making up the asymptotic variance $v_k$ can be of different order.
The first term is an average of bounded terms, and so under our assumptions will be of order $\mathcal{O}(1)$. The second and third terms will be at most of the same order as the first one. If $p_k\approx 0$ so we can think of the sample as small relative to the population of sampled clusters, the first term dominates the second and third terms. If cluster sizes are bounded as $k$ increases, the fourth and fifth terms in are also order $\mathcal{O}(1)$. If, on the other hand, cluster sizes increase with $k$, these terms can be of higher order and dominate the variance. Whether they do so or not depends on the {\em (i)}  magnitude of $p_k$, {\em (ii)}  presence of clustering in sampling, {\em (iii)}  presence of clustering in assignment, and {\em (iv)} heterogeneity in potential outcomes. 

\subsection{The Robust and Cluster Robust Variance Estimators}
\label{section:robust_variances}

Let $\widehat U_{k,i}=Y_{k,i}-\widehat\alpha_k-\widehat\tau_kW_{k,i}$ be the residuals from the regression of $Y_{k,i}$ or a constant and $W_{k,i}$. Here, $\widehat\alpha_k$ is the intercept of the regression and $\widehat\tau_k$ is the coefficient on $W_{k,i}$ (equal to the expression in (\ref{ls_estimator}) with probability approaching one).

There are two common estimators of the variance of $\sqrt{N_k}(\widehat\tau_k-\tau_{k})$.
First, the conventional robust variance estimator (\citet{eicker1963, huber1967behavior, white1980heteroskedasticity}):
\begin{equation}\label{var_ehw}
\widehat V^{\rm robust}_k
=\frac{1}{\mwidebar W_k^2(1-\mwidebar W_k)^2}\left\{\frac{1}{N_k}\sum_{i=1}^{n_k}R_{k,i} \widehat{U}_{k,i}^2(W_{k,i}-\mwidebar W_k)^2\right\},
\end{equation}
where 
\[
\mwidebar W_k = \frac{1}{N_k\vee 1}\sum_{i=1}^{n_k} R_{k,i}W_{k,i}.
\]
Let
\[
v^{\rm robust}_k=\frac{1}{n_{k}}\sum_{i=1}^{n_k} \bigg(\frac{u^2_{k,i}(1)}{\mu_k}+\frac{u^2_{k,i}(0)}{1-\mu_k}\bigg). 
\]
Under regularity conditions (see appendix), $\widehat V_k^{\rm{robust}}$ and $v_k^{\rm robust}$ are close in the following sense,
\[
\frac{\widehat V_k^{\rm{robust}}}{v_k}=\frac{v_k^{\rm robust}}{v_k}+\scaleto{\mathcal{O}}{5pt}_p(1),
\]
motivating our focus on the comparison of $v_k^{\rm robust}$ and $v_k$.
In general the difference $v^{\rm robust}_k-v_k$ can be positive or negative, so the robust variance estimator can be invalid in large samples.

The second common variance estimator is the cluster variance \citep{liang1986longitudinal, arellano1987practitioners}, 
\begin{equation}\label{clustervar}
\widehat V^{\rm cluster}_k
=\frac{1}{\mwidebar W_k^2(1-\mwidebar W_k)^2}\left\{\frac{1}{N_k}\sum_{m=1}^{m_k}\left(\sum_{i=1}^{n_k}
1\{m_{k,i}=m\} R_{k,i}\widehat U_{k,i}(W_{k,i}-\mwidebar W_k)\right)^2\right\}.
\end{equation}
Define
\begin{align*}
v^{\rm cluster}_k&=\frac{1}{n_{k}}\sum_{i=1}^{n_k} \bigg(\frac{u^2_{k,i}(1)}{\mu_k}+\frac{u^2_{k,i}(0)}{1-\mu_k}\bigg)\\
&-p_k\frac{1}{n_{k}}\sum_{i=1}^{n_k}\big(u_{k,i}(1)-u_{k,i}(0)\big)^2-p_k\sigma_k^2
\frac{1}{n_{k}}\sum_{i=1}^{n_k}\bigg(\frac{u_{k,i}(1)}{\mu_k}+\frac{u_{k,i}(0)}{1-\mu_k}\bigg)^2\\
&+p_k\frac{1}{n_{k}}\sum_{m=1}^{m_k}\Bigg(\sum_{i=1}^{n_k} 1\{m_{k,i}=m\}\big(u_{k,i}(1)-u_{k,i}(0)\big)\Bigg)^2\\
&+p_k\sigma^2_k\frac{1}{n_{k}}\sum_{m=1}^{m_k}\Bigg(\sum_{i=1}^{n_k}1\{m_{k,i}=m\}\bigg(\frac{u_{k,i}(1)}{\mu_k}+\frac{u_{k,i}(0)}{1-\mu_k}\bigg)\Bigg)^2.\label{v_vijf}
\end{align*}
Then,  $\widehat V_k^{\rm{cluster}}$ is close to $v_k^{\rm cluster}$ in the sense that
\[
\frac{\widehat V_k^{\rm{cluster}}}{v_k}=\frac{v_k^{\rm cluster}}{v_k}+\scaleto{\mathcal{O}}{5pt}_p(1).
\]
The difference $v^{\rm cluster}_k-v_k$ is always nonnegative. Therefore, for large $k$, the cluster variance estimator can be conservative but cannot underestimate the variance of $\widehat\tau_k$.

\subsection{Discussion}
\label{section:discussion}

From the formulas for $v_k$, $v_k^{\rm{robust}}$, and $v_k^{\rm{cluster}}$ it follows that if $p_k$ is small enough, then $v_k^{\rm{robust}}$ and $v_k^{\rm{cluster}}$ are approximately equal to $v_k$. In this case, clustered sampling and clustered assignment do not matter much because the probability that two sample units belong to the same cluster is small. 

The difference $v_k^{\rm{robust}}-v_k$ depends on two terms. The first term,
\begin{equation}
p_k\frac{1}{n_k}\Bigg[\sum_{i=1}^{n_k}\big(u_{k,i}(1)-u_{k,i}(0)\big)^2-(1-q_k)\sum_{m=1}^{m_k}n_{k,m}^2 (\tau_{k,m}-\tau_k)^2\Bigg],
\label{equation:r1}
\end{equation}
is equal to zero when treatment effects are constant (in which case, $u_{k,i}(1)-u_{k,i}(0)=0$ for $i=1, \ldots, n_k$ and $\tau_{k,m}-\tau_k=0$ for all $m=1,\ldots, m_k$). If all clusters are sampled, so $q_k=1$, and treatment effects are heterogeneous, \eqref{equation:r1} is positive. When only a fraction of the clusters are sampled, $q_k<1$, the sign of \eqref{equation:r1} depends on the extent to which heterogeneity in treatment effects can be explained by the clusters. If there is no variation in average treatment effects across clusters, the expression in \eqref{equation:r1} is non-negative. However, when clusters explain much of the variation in treatment effects, the expression in \eqref{equation:r1} can be negative and very large in magnitude because of the factor $n_{k,m}^2$.
The second term of $v_k^{\rm{robust}}-v_k$ is equal to
\begin{align}
p_k\sigma^2_k\sum_{m=1}^{m_k}\frac{n_{k,m}}{n_k}\Bigg[\frac{1}{n_{k,m}}\sum_{i=1}^{n_k}1\{m_{k,i}&=m\}\bigg(\frac{u_{k,i}(1)}{\mu_k}+\frac{u_{k,i}(0)}{1-\mu_k}\bigg)^2\nonumber\\&-n_{k,m}\Bigg(\frac{1}{n_{k,m}}\sum_{i=1}^{n_k}1\{m_{k,i}=m\}\bigg(\frac{u_{k,i}(1)}{\mu_k}+\frac{u_{k,i}(0)}{1-\mu_k}\bigg)\Bigg)^2\Bigg].
\label{equation:r2}
\end{align}
This term is equal to zero if there is no clustered assignment, that is, $\sigma^2_k=0$.
If $\sigma^2_k>0$, the sign of \eqref{equation:r2} depends on how much of the heterogeneity in potential outcomes is explained
by the clusters. 
The expression in \eqref{equation:r2} is close to zero when there is little heterogeneity in potential outcomes, so $u_{k,i}(1)$ and $u_{k,i}(0)$ are close to zero. If there is heterogeneity in potential outcomes but average potential outcomes are nearly constant across clusters, \eqref{equation:r2} is positive. When the clusters explain enough heterogeneity in potential outcomes \eqref{equation:r2} can be negative and potentially very large in magnitude because of the factor $n_{k,m}$ multiplying the second term of the sum in \eqref{equation:r2}. That is, the robust variance formula can severely underestimate the variance of $\widehat\tau_k$.

Clustered standard errors are conservative in general, that is, $v_k^{\rm{cluster}}\geq v_k$. In particular, the difference $v_k^{\rm{cluster}}-v_k$ is
\[v_k^{\rm{cluster}}-v_k=
p_kq_k\frac{1}{n_k}\sum_{m=1}^{m_k}\Bigg(\sum_{i=1}^{n_k} 1\{m_{k,i}=m\}\big(u_{k,i}(1)-u_{k,i}(0)\big)\Bigg)^2, 
\]
which can be rewritten as
\begin{equation}\label{var-adj-cluster}
v_k^{\rm{cluster}}-v_k=\left(\frac{p_k n_k}{m_k}\right) q_k\left\{\frac{1}{m_k}\sum_{m=1}^{m_k} \left(\frac{n_{k,m}m_k}{n_k}\right)^2 (\tau_{k,m}-\tau_k)^2\right\}.
\end{equation}
When the expected fraction of clusters in the sample, $q_k$, is small, or when the average treatment effect is nearly constant between clusters, then $v_k^{\rm{cluster}}\approx v_k$. Aside from these special cases, the $p_kn_{k}/m_k$ factor in the formula above indicates that cluster standard errors can be extremely conservative in general. 

\section{Two New Variance Estimators}
\label{section:new_variance}

Estimation of the variance of $\widehat\tau_k$ is challenging because the different terms in $v_k$ can be of different orders of magnitude. In this section, we propose two estimators of the variance of $\widehat\tau_k$ that allow us to correct the bias of the cluster variance estimator, one analytic, and one resampling-based.
As the expression for the bias of the cluster variance in (\ref{var-adj-cluster}) shows, the cluster variance is heavily biased if the fraction of the sampled clusters is large and there is substantial variation in the cluster-specific treatment effects. Although the proposed analytic variance estimator is defined irrespective of the value of $\sigma_k^2$, in order to for the correction to be effective we need to be able to estimate the cluster-specific treatment effects, and thus we need $\sigma_k^2$ to be less than its maximum value of $\mu_k(1-\mu_k)$ to ensure that there is variation in the treatment assignment within clusters. One of the proposed variance estimators is based on a correction to $\widehat V^{\rm cluster}_k$, and the other is based on resampling methods. An alternative would be to directly estimate the bias term in (\ref{var-adj-cluster}) and subtract that from the cluster variance. A challenge with this approach is that the estimation error for the adjustment term is large (often leading  to negative variances estimates) because the order of magnitude of the correction is itself large and this approach did not work well in our simulations.
We do not report formal results for the variance estimators in the current paper. We demonstrate their performance in the simulations in Section \ref{section:simulations}. There may well be further refinements possible.

If $q_k$ is close to zero, the proposed variance estimators are close to $\widehat V^{\rm cluster}_k$, which has little bias in that case. If $\sigma_k^2=\mu_k(1-\mu_k)$ (that is, when $W_{k,i}$ is constant within clusters), the proposed resampling variance estimator is not defined. 
To be effective both variance estimators rely on estimating the variation in  treatment effects across clusters, and therefore require a substantial number of both  treated and control observations per cluster. The proposed variance estimators lead to substantial improvements over $\widehat V^{\rm cluster}_k$ in cases where $\widehat V^{\rm cluster}_k$ has a large upward bias. The downside of the proposed variance estimators is that they can be   conservative when there is no need to cluster because there is no heterogeneity in treatment effects, or when there are too few treated and control observations per cluster to estimate the heterogeneity in the treatment effects precisely. 

We first consider in Section \ref{section_all_clusters} the case with $q_k=1$ so we have random sampling. Next we consider in Section \ref{section_some_clusters} the case with clustered sampling $q_k<1$. In Section \ref{section:bootstrap} we propose a bootstrap procedure for estimating the variance. The proposed variance estimators perform very well in the simulation study of Section \ref{section:simulations}. The derivation of their formal properties is left for future work.

\subsection{The Case with All Clusters Observed}\label{section_all_clusters}

First we focus on the case with $q_k=1$ (all clusters observed), but allowing for general $p_k$. Let $U_{k,i}=W_{k,i} u_{k,i}(1)+(1-W_{k,i})u_{k,i}(0)$. The first step is to approximate the normalized error of the least squares estimator $\widehat\tau_k$ by a normalized sample average over clusters,
\begin{equation} 
\label{new_een}
\sqrt{N_k}(\widehat\tau_k-\tau_k)/v_k^{1/2}
= \frac{1}{\sqrt{n_kp_kv_k}\mu_k(1-\mu_k)}
\sum_{m=1}^{m_k}C_{k,m}+ 
\scaleto{\mathcal{O}}{5pt}_p(1),
\end{equation}
where the terms
\[ 
C_{k,m}=\sum_{i=1}^{n_k}1\{m_{k,i}=m\}R_{k,i} (W_{k,i} -\mu_k) U_{k,i}
\]
are independent across clusters. In the appendix, we show
\begin{align}
\widehat V_k^{\rm{cluster}}/v_k 
=\frac{1}{n_kp_kv_k}\left(\frac{1}{\mu_k(1-\mu_k)}\right)^2
\sum_{m=1}^{m_k} C_{k,m}^2+\scaleto{\mathcal{O}}{5pt}_p(1).
\label{equation:vrobustsum}
\end{align}
The expectation of $C_{m,k}$ is 
\[
E[C_{k,m}]=n_{k,m}p_k\mu_k(1-\mu_k)(\tau_{k,m}-\tau_{k}),
\]
with sum over clusters
\begin{equation}
\sum_{m=1}^{m_k} E[C_{k,m}]
= p_k\mu_k(1-\mu_k)\sum_{m=1}^{m_k} n_{k,m}(\tau_{k,m}-\tau_{k})=0.
\label{equation:Cmksum}
\end{equation}
That is, although the sum of the expectations of $C_{k,m}$ over clusters is equal to zero, these expectations are not equal to zero in general for each cluster separately. Because $\mbox{var}(C_{k,m})\leq E[C_{k,m}^2]$, the first term on the right-hand side of (\ref{equation:vrobustsum}) is conservative on expectation relative to the variance of $\sqrt{N_k}(\widehat\tau_k-\tau_k)/v_k^{1/2}$, which explains the conservativeness of $\widehat V_k^{\rm{cluster}}$.

Because of \eqref{equation:Cmksum}, we can replace the terms $C_{k,m}$ in (\ref{new_een}) by $C_{k,m}-E[C_{k,m}]=C_{k,m,1}+C_{k,m,2}$, where
\begin{align*}
C_{k,m,1}&=\sum_{i=1}^{n_k}1\{m_{k,i}=m\} (R_{k,i}-p_k)(\tau_{k,m}-\tau_{k}) \mu_k(1-\mu_k),
\shortintertext{and}
C_{k,m,2}&=\sum_{i=1}^{n_k}1\{m_{k,i}=m\} R_{k,i}\Bigl((W_{k,i}-\mu_k)U_{k,i}-(\tau_{k,m}-\tau_{k}) \mu_k(1-\mu_k)
\Bigr).
\end{align*}
Therefore,
\begin{equation}
\label{new_eent}
\sqrt{N_k}(\widehat\tau_k-\tau_k)/v_k^{1/2}=
 \frac{1}{\sqrt{n_kp_kv_k}\mu_k(1-\mu_k)}\left(
\sum_{m=1}^{m_k}C_{k,m,1}+
\sum_{m=1}^{m_k}C_{k,m,2}\right)
+\scaleto{\mathcal{O}}{5pt}_p(1).
\end{equation}
It can be shown that $C_{k,m,1}$ and $C_{k,m,2}$ have means equal to zero and are uncorrelated. In addition,  $C_{k,m,1}$ and $C_{k,m,2}$ are uncorrelated across clusters. The variance of $\sum_{m=1}^{m_k}C_{k,m,1}/(\sqrt{n_kp_k}\mu_k(1-\mu_k))$ is
\[ 
(1-p_k)\sum_{m=1}^{m_k} \frac{n_{k,m}}{n_k} (\tau_{k,m}-\tau_{k})^2.
\]
Let $\widehat\tau_{k,m}$ be difference between the sample average of the outcome for treated and nontreated units in cluster $m$. A direct estimator the variance of $\sum_{m=1}^{m_k}C_{k,m,2}$ is
\begin{equation}  
\sum_{m=1}^{m_k}\left(\sum_{i=1}^{n_k}1\{m_{k,i}=m\}R_{k,i}\Bigl((W_{k,i}-\mwidebar W_k)\widehat U_{k,i}-(\widehat\tau_{k,m}-\widehat\tau_k) \mwidebar W_k(1-\mwidebar W_k)\Bigr)
\right)^2,
\label{equation:second_term_CCV}
\end{equation}
In practice, the estimator in (\ref{equation:second_term_CCV}) is biased from the correlations between the estimation errors of its components. We apply sampling splitting to address this bias. 
We first split the sample randomly into two subsamples. Let ${Z}_{k,i}\in\{0,1\}$ be the indicator that unit $i$ belongs to the second subsample, and let $\overline{Z}_k$ be the mean of $Z_{k,i}$. Using the subsample with ${Z}_{k,i}=0$, we obtain estimates  $\widehat\tau_{k,m}^{\,*}$, $\widehat\alpha_{k}^{\,*}$, and $\widehat\tau_k^{\,*}$ of $\tau_{k,m}$, $\alpha_k$, and $\tau_k$, respectively. 
Next, for observations with $Z_{k,i}=1$, we calculate the residuals $\widehat U_{k,i}^{\,*}=Y_{k,i}-\widehat\alpha_{k}^{\,*}-\widehat\tau_k^{\,*} W_{k,i}$.
Finally, we estimate the normalized variance for the case with $q_k=1$ as
\begin{align}
    \widehat V^{\rm CCV}_{k} (1)&=
    \frac{1}{N_k \mwidebar W_k^2(1-\mwidebar W_k)^2}\sum_{m=1}^{m_k}\Bigg[ \frac{1}{\mwidebar{Z}_k^2}\Bigg(\sum_{i=1}^{n_k} 1\{m_{k,i}=m\}R_{k,i} Z_{k,i} \Big((W_{k,i}-\mwidebar W_k)\widehat{U}_{k,i}^{\,*}\nonumber\\
    &-(\widehat\tau_{k,m}^{\,*}-\widehat\tau_k^{\,*}) \mwidebar W_k(1-\mwidebar W_k)\Big)\Bigg)^2\nonumber\\
    &-\frac{1-\mwidebar{Z}_k}{\mwidebar{Z}_k^2}\sum_{i=1}^{n_k} 1\{m_{k,i}=m\}R_{k,i}Z_{k,i}\Big((W_{k,i}-\mwidebar W_k)\widehat{U}_{k,i}^{\,*}-(\widehat\tau_{k,m}^{\,*}-\widehat\tau_k^{\,*}) \mwidebar W_k(1-\mwidebar W_k)\Big)^2\Bigg]\nonumber\\
    &+ (1-p_k)\sum_{m=1}^{m_k} \frac{\mwidebar N_{k,m}}{N_k} (\widehat\tau_{k,m}-\widehat\tau_k)^2,
    \label{equation:ccv}
\end{align}
where $\mwidebar N_{k,m}$ is the size of the sample in cluster $m$. 
For clusters with no variation in the treatment variable, we replace $\widehat\tau_{k,m}$ in \eqref{equation:ccv} with $\widehat\tau_k$. For clusters with no variation in the treatment variable for a particular subsample, we replace $\widehat\tau_{k,m}^*$ in \eqref{equation:ccv} with $\widehat\tau_k^*$. 
We derive the form of the CCV estimator in the appendix.
To improve the precision of $\widehat V^{\rm CCV}_k(1)$, we re-estimate it multiple times with new sample splits (new values for $Z_{k,i}$) and then average the corresponding variance estimators.
In our simulations of section \ref{section:simulations}, we re-estimate the variance estimator four times, and use  sample splits with in expectation an equal number of units in each subsample, so $E[\overline{Z}_k]=1/2$.

\subsection{The Case When Not All Clusters Are Sampled}\label{section_some_clusters}

To motivate the modification of the variance estimator $\widehat V_k^{\rm CCV}(1)$ for the $q_k<1$ case, notice that
\[ v_k(q_k)-v^{\rm cluster}_k=q_k \times 
(v_k(1)-v^{\rm cluster}_k),
\]
where $v_k(q_k)$ denotes the value of the true variance $v_k$ evaluated at $q_k$.
That is, the variance for the general $q_k$ case is a convex combination of the true variance at $q_k=1$ and the cluster variance,
\[ 
v_k(q_k)=q_k\times
v_k(1)+(1-q_k)\times v^{\rm cluster}_k.
\]
Let $\widehat q_k$ be the ratio between the number of sampled clusters and the total number of clusters in the population.  
The proposed variance estimator, $\widehat V^{\rm CCV}_k$, is a convex combination of $\widehat V^{\rm CCV}_k(1)$ and $\widehat V^{\rm cluster}_k$ with weights $\widehat q_k$ and $1-\widehat q_k$,
\begin{equation}\label{ccv_general}
\widehat V^{\rm CCV}_k= \widehat q_k\times  \widehat V^{\rm CCV}_k(1)+(1-\widehat q_k)\times \widehat V^{\rm cluster}_k.
\end{equation}
Computation of $\widehat q_k$ requires knowledge of $m_k$, the total number of clusters in the population.

\subsection{A Bootstrap Variance Estimator}
\label{section:bootstrap}

In the previous sections, we have discussed an analytic variance estimator. Here we suggest a resampling-based variance estimator, initially for the case with $q_k=1$. Like the causal bootstrap in \citet{imbens2021causal}, the proposed bootstrap procedure takes into account the causal nature of the estimand and creates bootstrap samples where units (in this case clusters) have different assignments and assignment probabilities than they have in the original sample. It differs from earlier bootstrap variance estimators for clustered settings \citep[{e.g.,}][]{cameron2015practitioner, menzel2021bootstrap} in that it allows for the possibility that a large fraction of clusters are observed.

The specific resampling procedure, which we call the two-stage-cluster-bootstrap (TSCB), consists of two stages. For each of the clusters, let $\mwidebar N_{k,m}$ be the cluster-level sample size and $\mwidebar W_{k,m}=N_{k,m,1}/(\mwidebar N_{k,m}\vee 1)$ the cluster-level fraction of treated units. In the first stage of the bootstrap procedure, for each cluster we draw ${\overline{W}}^{\, b}_{k,m}$ with replacement from the empirical distribution of the cluster-level fractions of treated units, that is with probability $1/m_k$ from the set $\{{\overline{W}}_{k,1},\ldots,{\overline{W}}_{k,m_k}\}$. In the second stage, we draw $\mwidebar N_{k,m} \overline W^{\, b}_{k,m}$ units with replacement from the set of treated units in cluster $m$ and 
$\mwidebar N_{k,m}(1- \overline W^{\, b}_{k,m})$ units with replacement from the set of untreated units in cluster $m$. 
In order for this to be well-defined we do need all the ${\overline{W}}_{k,1}$ to be strictly between zero and one. 
We do this for all clusters to create the bootstrap sample, and calculate the bootstrap standard errors as the standard deviation of the treatment effect estimates across bootstrap iterations. 

Next, consider the case with $q_k<1$. In this case, we need to take into account the fact that we see a fraction of the clusters in the population. We follow the approach proposed in \citet{chao1985bootstrap}. Suppose $q=1/2,$ so we observe half the clusters in the population. The bootstrap procedure first creates a pseudo population consisting  of the original population of clusters, plus one additional replica of each cluster. Then, to get a bootstrap sample, we sample randomly, without replacement, from the clusters in this pseudo population. Given the clusters in the bootstrap sample, we proceed as before, and ultimately calculate the bootstrap variance as the variance of the estimator over the bootstrap samples. \citet{chao1985bootstrap} provide details and extensions to the case for the case where $1/q_k$ is not an integer.

The algorithm for the TSCB is summarized here.

\begin{algorithm}
\caption{Two Stage Cluster Bootstrap}\label{algo:eif}
\begin{algorithmic}
\\ {\bf Input:  }
\\ \hskip0.6cm {Sample $(Y_{k,i},W_{k,i},m_{k,i})$}
\\ \hskip0.6cm {Fraction sampled clusters $q_k$}
\\ \hskip0.6cm {Number of bootstrap replications $B$}
\\ {\bf Stage 1:}
\\  \leftskip0.6cm{1a: Create pseudo population by replicating each cluster $1/q_k$ times}
\\  {1b: For each cluster in the pseudo population, calculate the assignment probability $\overline{W}_{k,m}$}
\\  {1c: Create a bootstrap sample of clusters by randomly drawing clusters from the pseudo population from Stage 1a 
}
\\  {1d: For each sampled cluster, draw an assignment probability $A_{k,m}$ from the empirical distribution of the $\overline{W}_{k,m}$ from Stage 1b}
\\ \leftskip0cm {\bf Stage 2:}
\\ \leftskip0.6cm {2a: Randomly draw from the set of treated units in cluster $m$, $\lfloor N_{k,m} A_{k,m}\rfloor$ units}
\\  {2b: Randomly draw from the set of control units in cluster $m$, $\lfloor N_{k,m} (1-A_{k,m})\rfloor$ units}
\\ \leftskip0cm{\bf Calculations:}
\\ \leftskip0.6cm {For the units in the bootstrap sample  constructed in Stage 2, collect the values for $(Y_{k,i},W_{k,i},m_{k,i})$ and calculate the least squares or fixed effect estimator}
\\ \leftskip0.6cm {Calculate the standard deviation of the least squares or fixed effect estimator over the $B$ bootstrap samples}
\end{algorithmic}
\end{algorithm}

\section{The Fixed Effect Estimator}
\label{section:fixed_effect}

In this section, we report results for the fixed effect estimator often used in empirical research in economics. \citet{arellano1987practitioners}, \cite*{Bertrand2004did}, \cite{cameron2015practitioner} and \citet*{mackinnon2021cluster}
have pointed out that cluster adjustments may still be necessary in fixed effects regressions. However, a view of clustering based on models with cluster-specific variance components creates ambiguity in the role of clustered standard errors for estimators with cluster fixed effects, which are specifically aimed to absorb cluster-level variation.

We first characterize the fixed effect estimator and derive its large $k$ distribution. Then, we discuss the properties of the two conventional variance estimators, the robust and cluster robust variance estimators. As in the least squares case, we find that the robust standard errors may be too small and the cluster standard errors may be unnecessarily large, especially in cases when the number of observations per cluster is large. We propose CCV and TSCB variance estimators. The CCV estimator for fixed effects has a different form than the one for least squares in section \ref{section:new_variance}.

The fixed effect estimator is based on a regression of the outcome on the treatment indicator and indicators for each of the clusters in the sample. It can be written as the least squares estimate for a regression of the outcome on the treatment, with both variables measured in deviation from cluster means,
\begin{align} \widehat\tau_k^{\rm{\,fixed}}&= \frac{\displaystyle\sum_{m=1}^{m_k}\sum_{i=1}^{n_k} 1\{m_{k,i}=m\}R_{k,i}Y_{k,i}(W_{k,i}-  \mwidebar W_{k,m})}{\displaystyle\sum_{m=1}^{m_k} \sum_{i=1}^{n_k} 1\{m_{k,i}=m\}R_{k,i}W_{k,i}(W_{k,i}-  \mwidebar W_{k,m})}. 
\label{equation:fixed_text}
\end{align} 

Like in section \ref{section:least_squares_estimator}, we assume that that potential outcomes are bounded, $m_kq_k\rightarrow \infty$, and $\limsup_{k\rightarrow\infty} \max_m n_{k,m}/\min_m n_{k,m}<\infty$. In addition, we assume {\em (i)} $(m_kq_k)/((p_k n_k)/m_k)\allowbreak\rightarrow 0$, and {\em (ii)} the supports of the cluster probabilities, $A_{k,m}$, are bounded away from zero and one (uniformly in $k$ and $m$). Assumption {\em (i)} restricts the focus of our analysis in this section to settings where the expected number of sampled clusters is small relative to the expected number of sampled observations per sampled cluster. Together with the previous assumptions, assumption {\em (i)} implies $(p_k n_k)/m_k\rightarrow\infty$, $n_kp_kq_k\rightarrow\infty$, and $p_k\min_m n_{k,m}\rightarrow\infty$. This last result, along with assumption {\em (ii)}, ensures that $\widehat \tau_k^{\rm{\,fixed}}$ in \eqref{equation:fixed_text} is well-defined with probability approaching one. 

Let $\alpha_{k,m}=(1/n_{k,m})\sum_{i=1}^{n_k}1\{m_{k,i}=m\}y_{k,i}(0)$. For an observation, $i$, with $m_{k,i}=m$, we define the within-cluster residuals
$e_{k,i}(0)=y_{k,i}(0)- \alpha_{k,m}$ and $e_{k,i}(1)=y_{k,i}(1)-\tau_{k,m}-\alpha_{k,m}$. 
Let
\begin{equation}\label{vtildek}
\tilde v_k = f_k/(\mu_k(1-\mu_k)-\sigma^2_k)^2
\end{equation}
where
\begin{align*}
f_k&=E[A_{k,m}(1-A_{k,m})^2]\frac{1}{n_k}\sum_{i=1}^{n_k} e_{k,i}^2(1)
+E[A^2_{k,m}(1-A_{k,m})]\frac{1}{n_k}\sum_{i=1}^{n_k} e_{k,i}^2(0)\nonumber\\
&- p_k E[A^2_{k,m}(1-A_{k,m})^2]\frac{1}{n_k}\sum_{i=1}^{n_k}(e_{k,i}(1)-e_{k,i}(0))^2\nonumber\\
&+\Big(E[A_{k,m}(1-A_{k,m})]-(5+p_k)E[A^2_{k,m}(1-A_{k,m})^2]\nonumber\\
&\qquad\qquad + 2q_k(E[A_{k,m}(1-A_{k,m})])^2 \Big)\sum_{m=1}^{m_k}\frac{n_{k,m}}{n_k}(\tau_{k,m}-\tau_k)^2\nonumber\\
&+\Big(p_kE[A^2_{k,m}(1-A_{k,m})^2]-p_kq_k(E[A_{k,m}(1-A_{k,m})])^2\Big)\sum_{m=1}^{m_k}\frac{n^2_{k,m}}{n_k}(\tau_{k,m}-\tau_k)^2.
\end{align*}
Under additional regularity conditions, which are described in the Appendix, we obtain the large $k$ distribution of the fixed effects estimator, 
\begin{equation}
    \label{vfixed}
    \sqrt{N_k}(\widehat\tau_k^{\rm{\,fixed}} -\tau_k)/\tilde v_k^{1/2}\stackrel{d}{\longrightarrow} N(0,1).
\end{equation}

Let $\widetilde U_{k,i}=\widetilde Y_{k,i} - \widehat\tau_k^{\,{\rm fixed}} \widetilde W_{k,i}$, where $\widetilde Y_{k,i}=Y_{k,i}-\mwidebar Y_{k,m_{k,i}}$, $\widetilde W_{k,i}=(W_{k,i}-\mwidebar W_{k,m_{k,i}})$. The robust estimator of the variance of 
$\sqrt{N_k}(\widehat\tau^{\,{\rm fixed}}_k-\tau_k)$ is 
\begin{align}
\label{vrobustfixed}
\widetilde V_k^{\rm{robust}} =\left.\frac{1}{N_k}\sum_{i=1}^{n_k}R_{k,i}\widetilde W_{k,i}^2\widetilde U_{k,i}^2\right/\left(\frac{1}{N_k}\sum_{i=1}^{n_k} R_{k,i}\widetilde W_{k,i}^2\right)^2.
\end{align}
Now let,
\[
\tilde v_k^{\rm robust} = f_k^{\rm robust}/(\mu_k(1-\mu_k)-\sigma^2_k)^2.
\]
with 
\begin{align*}
f_k^{\rm robust}&=E[A_{k,m}(1-A_{k,m})^2]\frac{1}{n_k}\sum_{i=1}^{n_k} e^2_{k,i}(1) + E[A^2_{k,m}(1-A_{k,m})]\frac{1}{n_k}\sum_{i=1}^{n_k} e^2_{k,i}(0)\\
&+E[A_{k,m}(1-A_{k,m})(1-3A_{k,m}(1-A_{k,m}))]\sum_{m=1}^{m_k}\frac{n_{k,m}}{n_k}
(\tau_{k,m}-\tau_k)^2.
\end{align*}
Notice that all terms of $f_k^{\rm robust}$ are bounded. In the appendix, we show that
\[
\widetilde V_k^{\rm{robust}}=\tilde v_k^{\rm robust}+\scaleto{\mathcal{O}}{5pt}_p(1).
\]
The cluster variance estimator for fixed effects is
\begin{align}\label{clusterfixedvar}
\widetilde V_k^{\rm{cluster}} =\left.\frac{1}{N_k}\sum_{m=1}^{m_k} \left(\sum_{i=1}^{n_k}1\{m_{k,i}=m\}R_{k,i}\widetilde W_{k,i}\widetilde U_{k,i}\right)^2\right/\left(\frac{1}{N_k}\sum_{i=1}^{n_k} R_{k,i}\widetilde W_{k,i}^2\right)^2.
\end{align}
Let,
\[
\tilde v_k^{\rm cluster} = f_k^{\rm cluster}/(\mu_k(1-\mu_k)-\sigma^2_k)^2.
\]
with 
\begin{align*}
f_k^{\rm cluster}&=E[A_{k,m}(1-A_{k,m})^2]\frac{1}{n_k}\sum_{i=1}^{n_k} e_{k,i}^2(1)
+E[A^2_{k,m}(1-A_{k,m})]\frac{1}{n_k}\sum_{i=1}^{n_k} e_{k,i}^2(0)\\
&- p_kE[A^2_{k,m}(1-A_{k,m})^2]\frac{1}{n_k}\sum_{i=1}^{n_k}(e_{k,i}(1)-e_{k,i}(0))^2\\
&+ (E[A_{k,m}(1-A_{k,m})]-(5+p_k)E[A^2_{k,m}(1-A_{k,m})^2])\sum_{m=1}^{m_k}\frac{n_{k,m}}{n_k}(\tau_{k,m}-\tau_k)^2\\
&+p_kE[A^2_{k,m}(1-A_{k,m})^2]\sum_{m=1}^{m_k}\frac{n^2_{k,m}}{n_k}(\tau_{k,m}-\tau_k)^2.
\end{align*}
We obtain in the appendix,
\[
\frac{\widetilde V_k^{\rm{cluster}}}{\tilde v_k}=\frac{\tilde v_k^{\rm cluster}}{\tilde v_k}+\scaleto{\mathcal{O}}{5pt}_p(1).
\]

Similar to the least squares case, the robust variance can underestimate the true variance, and the cluster variance is generally too large. Our proposed variance estimator is a convex combination of $\widetilde V^{\rm cluster}_k$ and 
and $\widetilde V^{\rm robust}_k$, with the weights selected to correct the bias of the cluster variance estimator as $k$ increases (see appendix for details). 

\begin{equation}\label{varfixedccv}
\widetilde{V}_k^{\rm CCV}=\widehat\lambda_k \widetilde V^{\rm cluster}_k 
+(1-\widehat\lambda_k) {\widetilde{V}}^{\rm robust}_k.
\end{equation}
where the estimated weight for the cluster variance is
\[
\widehat\lambda_k=1- q_k\,
\frac{\left(\displaystyle\frac{1}{M_k}\sum_{m=1}^{m_k} Q_{k,m}\mwidebar W_{k,m}(1-\mwidebar W_{k,m})\right)^2}{\displaystyle\frac{1}{M_k}\sum_{m=1}^{m_k} Q_{k,m}\mwidebar W_{k,m}^2(1-\mwidebar W_{k,m})^2},
\]
where
$Q_{k,m}$ is an indicator that takes value one if cluster $m$ of population $k$ is sampled, and $M_k=\sum_{m=1}^{m_k} Q_{k,m}$ is the total number of sampled clusters. 
The second factor in the second term approximately (that is, ignoring the variance of $\overline{W}_{k,m}$ conditional on $A_{k,m}]$) estimates the variance of $A_{k,m}(1-A_{k,m})$ divided by its second moment, so that
\[ \tilde\lambda\approx 1-q_k \frac{V(A_{k,m}(1-A_{k,m}))}{E[(A_{k,m}(1-A_{k,m}))^2]}.\]
If there is no variation in $W_{k,i}$ within any of the clusters the fixed effect estimator is not defined, and neither is this variance estimator. In all other cases the variance estimator is well-defined.

We also consider a bootstrap standard error, based on the same resampling procedure described in Section \ref{section:bootstrap}.

\section{Simulations}
\label{section:simulations}

We next report simulation results that illustrate  the performance of the proposed variance estimators relative to existing alternatives.
To operate in an empirically relevant setting, 
we create an artificial population based on the Census data briefly described in the introduction, which contains information on log earnings, an indicator for college attendance, and an indicator for state of residence for 2,632,838 individuals. 

For each individual in this population of 2,632,838 individuals, we define $m_{k,i}$ using state of residence (plus Washington, DC, and Puerto Rico), for a total of 52 clusters. 
We assign potential outcomes as $y_{k,i}(0)=Y_{k,i}-\widehat\tau_{k,m} W_{k,i}$ and
$y_{k,i}(1)=Y_{k,i}+\widehat\tau_{k,m} (1-W_{k,i})$, so treatment effects are constant within clusters.
We then repeatedly create samples from this population. Creating a sample requires fixing $p_k$, $q_k$, and fixing the distribution of $A_{k,m}$ and then drawing from the implied distribution for $R_{k,i}$ and $W_{k,i}$ to generate outcomes for all sampled units. In the baseline design, we set $p_k=q_k=1$, so we sample all $m_k=52$ clusters and all $n_k=2{,}632{,}838$ individuals in the population. For the assignment mechanism in the baseline design, we convert cluster means of the treatment variable into log-odds, $\widehat\ell_{k,m}=\ln(\mwidebar W_{k,m}/(1-\mwidebar W_{k,m}))$. Let $(\widehat\mu_\ell,\widehat\sigma_\ell)$ be the average and the sample standard deviation of $\widehat\ell_{k,m}$. We then draw $\ln (A_{k,m}/(1-A_{k,m}))$ for cluster $m$ from a normal distribution with expected value $\widehat\mu_\ell$ and standard deviation $\widehat\sigma_\ell$. Given the cluster assignment probability $A_{k,m}$, we assign the treatment in cluster $m$ by drawing from a binomial distribution with parameter $A_{k,m}$.

We calculate the standard deviation of the least squares and fixed effect estimators, normalized by the square root of the sample size, $N_k^{1/2} \mbox{s.d.}$, across 10,000 samples drawn according to the procedure outlined above. This is the benchmark against which we compare the various estimates of standard errors. For the least squares and the fixed effects estimators, respectively,
we first calculate the (infeasible) asymptotic standard errors $v_k^{1/2}$ and $\widetilde v_k^{1/2}$ to benchmark the performance of the feasible variance estimators. Next, we calculate the averages across 10,000 simulations of the robust, cluster, CCV, and TCSB standard errors, where we use 100 bootstrap replications in each simulation. 
Table \ref{table:standard_errors} reports the results. Table \ref{table:coverage_rates} reports coverage rates for 95 percent confidence intervals.
In the design column of the two tables $\sigma_{\tau_k}$ is the standard deviation of the cluster average treatment effect.

\begin{table}[t]
\begin{threeparttable}
    \centering
    \caption{Average standard errors across simulations}
\begin{tabular}{L{5cm}lccccccc}
\hline
&&&&&\multicolumn{4}{c}{normalized standard error}\\\cline{6-9}
&&\!$N_k^{1/2}\mbox{s.d.}$\!&\!$v_k^{1/2}$\!&\!$\widetilde v_k^{1/2}$\!&\!robust\!&\!cluster\!&\!CCV\!&\!TSCB\!\\ \hline
\hspace{-0.2cm}{\it }\\
\multirow{2}{5cm}{\hspace{-0.2cm}{\it Baseline design}: \\ $p_k=1$, $q_k=1$,\\
$\sigma_{\tau_k}=.120$,  
$\sigma_k=.057$}\\
                        &OLS\!&5.91 & 5.90 &  & 1.90 & 44.86 & 6.32 & 5.80\\
                        &FE\!&2.34 & & 2.32 & 1.90 & 44.63 & 2.31 & 2.29 \\[2ex]
\multirow{2}{5cm}{\hspace{-0.2cm}{\it Second Design}:\\$p_k=.1$, $q_k=1$,\\ $\sigma_{\tau_k}=.120$,   $\sigma_k=.057$}\\
                        &OLS\!& 2.61  & 2.59 &  & 1.90 & 14.28 & 3.78  & 2.60  \\
                        &FE\!& 1.95 &  & 1.95  & 1.90 & 14.21 & 1.95 & 1.94 \\[2ex] 
\multirow{2}{5cm}{\hspace{-0.2cm}{\it Third Design}:
\\
$p_k=.1$, $q_k=1$,\\  $\sigma_{\tau_k}=.480$,  $\sigma_k=.206$}\\
                        &OLS\!& 14.50 & 14.17 & & 1.98 & 56.46 & 13.70 & 14.33 \\
                        &FE\!& 12.14 & & 11.89 & 2.13 & 56.79 & 11.61 & 12.07  \\[2ex]
\multirow{2}{5cm}{\hspace{-0.2cm}{\it Fourth design}:\\$p_k=.1$, $q_k=1$,\\  $\sigma_{\tau_k}=0$,   $\sigma_k=.206$}\\
                        &OLS\!& 9.39 & 9.39 & & 1.90 & 8.20 & 9.19 & 9.37\\
                        &FE\!& 2.04 &  & 2.04 & 2.04 & 1.97 & 2.04 & 2.09\\[2ex] 
\multirow{2}{5cm}{\hspace{-0.2cm}{\it Fifth design}:\\ $p_k=.1$, $q_k=1$,\\  $\sigma_{\tau_k}=.480$,   $\sigma_k=0$}\\
                        &OLS\!& 1.95 & 1.97 & & 1.97 & 56.42 & 4.53 & 2.04 \\
                        &FE\!& 1.91 & & 1.94 & 1.94 & 56.42 & 1.96 & 1.90 \\[2ex]
\hline
\end{tabular}
\label{table:standard_errors}
\begin{tablenotes}[flushleft]%
\item\footnotesize{\!\!{\it Notes:} $N_k^{1/2}$s.d.\ is the standard deviation of the estimators over the simulations, multiplied by the square root of the sample size. $v_k^{1/2}$ is the square root of the asymptotic variance in equation (\ref{vk}). $\tilde v_k^{1/2}$ is the square root of the asymptotic variance of the fixed effect estimator in (\ref{vtildek}). The remaining four columns report average values of robust, cluster, CCV, and TSCB standard errors across simulations (multiplied by $N_k^{1/2}$). $p_k$ and $q_k$ are the unit and cluster sampling probabilities, respectively. $\sigma_{\tau_k}$ is the standard deviation of the cluster average treatment effect. $\sigma_k$ is the standard deviation across clusters of the treatment assignment probabilities.}
\end{tablenotes}
\end{threeparttable}
\end{table}

For the baseline design, the normalized standard deviation of the least squares estimator is 5.91. This is well approximated by the asymptotic standard error, 5.90. The robust standard error is on average over the simulations 1.90, less than one-third of the normalized standard deviation of the estimator. The cluster standard error is far too large, on average 44.86,  more than seven times the value of the normalized standard deviation. CCV improves considerably over robust and cluster. The average CCV standard error is 6.32, about 7 percent higher than the normalized standard deviation. The TSCB standard error is the most accurate, on average equal to 5.80. For the fixed effect estimator, the asymptotic standard error is again accurate. The robust standard error is about 16 percent too small, leading to a coverage rate for the nominal 95 percent confidence interval of 0.89 in Table \ref{table:coverage_rates}. The cluster standard error is too large by a factor of 20. CCV and TSCB standard errors closely approximate the normalized standard error. 

\begin{table}[t]
\begin{threeparttable}
    \centering
    \caption{Coverage rates across simulations}
\begin{tabular}{L{5cm}lcccccc}
\hline
&&\multicolumn{6}{c}{coverage of 95 percent confidence interval}\\\cline{3-8}
\vspace{.2cm}&&\!$v_k^{1/2}$\!&\!$\widetilde v_k^{1/2}$\!&\!robust\!&\!cluster\!&\!CCV\!&\!TSCB\!\\ \hline\hspace{-0.2cm}{\it }\\
\multirow{2}{5cm}{\hspace{-0.2cm}{\it Baseline design}:\\ $p_k=1$,  $q_k=1$,\\ $\sigma_{\tau_k}=.120$,   $\sigma_k=.057$}\\
                        &OLS\!& 0.949 & & 0.467 & 1.000 & 0.971 & 0.947    \\
                        &FE\! & & 0.950 & 0.893 & 1.000 & 0.947 & 0.942    \\[2ex]

\multirow{2}{5cm}{\hspace{-0.2cm}{\it Second design}:\\ $p_k=.1$,  $q_k=1$,\\  $\sigma_{\tau_k}=.120$,   $\sigma_k=.057$}\\
                        &OLS\!& 0.951 &  & 0.846 & 1.000 & 0.996 & 0.952    \\
                        &FE\!& & 0.950 & 0.944 & 1.000 & 0.950 & 0.948      \\[2ex] 

\multirow{2}{5cm}{\hspace{-0.2cm}{\it Third design}:\\ $p_k=.1$, $q_k=1$,\\ $\sigma_{\tau_k}=.480$,   $\sigma_k=.206$}\\
                        &OLS\!& 0.947 & & 0.208 & 1.000 & 0.960 & 0.950  \\ 
                        &FE\! &   & 0.941 & 0.284 & 1.000 & 0.918 & 0.948       \\[2ex]

\multirow{2}{5cm}{\hspace{-0.2cm}{\it Fourth design}:\\ $p_k=.1$, $q_k=1$, \\  $\sigma_{\tau_k}=0$,   $\sigma_k=.206$}\\
                        &OLS\!& 0.952 & & 0.308 & 0.905 & 0.966 & 0.952 \\  
                        &FE\! &   & 0.952 &  0.951 & 0.932 & 0.951 & 0.955 \\[2ex] 
\multirow{2}{5cm}{\hspace{-0.2cm}{\it Fifth design}:\\
 $p_k=.1$,  $q_k=1$,\\  $\sigma_{\tau_k}=.480$,   $\sigma_k=0$}\\
                        &OLS\!& 0.952 & & 0.953 &1.000 &1.000 & 0.959  \\ 
                        &FE\!&  & 0.954 & 0.955 & 1.000 & 0.957 &0.949   \\[2ex]
\hline
\end{tabular}
\label{table:coverage_rates}
\begin{tablenotes}[flushleft]
\item\footnotesize{\!\!{\it Notes:} See notes of Table \ref{table:standard_errors}.}
\end{tablenotes}
\end{threeparttable}
\end{table}

It is also interesting to consider the variation in the different variance estimators over the repeated samples relative to the true value of the standard deviation of the estimator. In the baseline design the normalized standard deviation is 5.91. The robust standard error is very precisely estimated, with  a standard deviation of the normalized robust standard over the 10,000 simulations equal to 0.005. The standard deviation of the cluster standard error is much larger, 1.48. For the CCV standard error the standard deviation is 1.21, and for the resampling-based TSCB the standard deviation is consideralby lower at 0.69.

We vary the design from the baseline case by changing $(i)$ the fraction of sampled units $p_k$, $(ii)$ the amount of treatment effect heterogeneity across  clusters, $\sigma_{\tau_k}$, and $(iii)$ the cross-cluster 
standard deviation
of the assignment probability, $\sigma_k$. In the second design, $p_k=0.1$ is the only change relative to the baseline design. This makes the robust standard errors less biased downward, and the cluster standard errors less biased upward. 
The result of decreasing the fraction of sampled units (and thus decreasing the sample size) is that the performance of the analytic CCV variance estimator declines, whereas the bootstraping vaiance estimator TSCB continues to perform well. We keep $p_k=0.1$ for the remaining three designs. In the third design, we increase both the treatment effect heterogeneity and the within-cluster correlation of the treatment by increasing the differences in treatment effects $\tau_{k,m}-\tau_k$ and the differences of the logs odds ratio $\ell_{k,m}-\ell_k$ by a factor of four. The resulting increase in $\sigma_k^2$ makes the performance of the robust standard error substantially worse, consistent with equation (\ref{equation:r2}). In this design, the bias of the robust standard error is substantial, also for the fixed effect estimator. 
The difference between the cluster variance and the true variance for the least squares estimator is proportional to the variation in the cluster average treatment effects, implying that the bias will increase for this design relative to the second design.  In the fourth design, we remove the heterogeneity in the treatment effect but keep the correlation in the treatment assignment the same as in the third design. Now, the cluster variance performs well, but the robust variance remains poor. In the fifth design, the assignment probabilities are identical in all clusters, and the treatment effect heterogeneity is the same as in the third and fourth designs. In this case the robust standard errors perform  well, but the cluster standard errors substantially over-estimate the uncertainty, as expected. In all designs, the CCV and especially the TSCB standard errors outperform the robust and cluster standard errors.

\section{Implications for Practice}

The analysis in this article has several implications for how to compute and, most importantly, interpret, standard errors in a variety of empirical settings. Some settings are clear cut and others are more subtle. First,
we discuss the case where there is no cluster sampling. If one has a random
sample of units from a large population with randomized treatment assignment
at the unit level, there is no reason to cluster the standard errors of the least squares estimator. Doing so can be harmful, resulting in unnecessarily wide confidence intervals. In this case, clustering is not appropriate
even if there is within-cluster correlation in outcomes (however those clusters are defined), and thus even if clustering makes a substantial difference in the magnitude of the standard errors. For example, if workers are sampled at random from a some population of interest and then randomly assigned to a job training program, clustering the standard errors at, say, the industry, county, or state level can result in standard errors that are unnecessarily conservative, often by a wide margin. Similarly, in a judge-leniency design---where defendants are randomly assigned to judges---standard errors should not be clustered at the level of the judge \citep*{frandsen}. If the sample represents a large fraction of the population and treatment effects are heterogeneous across units, robust standard errors are also conservative. If the data contains information on attributes of the units that are correlated with unit-level treatment effects, the methods in \cite{abadie2020sampling} can be applied to obtain less conservative standard errors. 

Next, consider the case of clustered assignment, and where we either have random sampling or we observe the entire population. This is one case where clustering becomes relevant, although conventional cluster standard errors can be extremely conservative. If assignment is perfectly clustered 
so that units that belong to the same cluster have the same treatment assignment, there is no improvement from using the CCV variance and the TSCB variance estimator is not applicable. If assignment is partially clustered---so there is variation in treatment assignment within clusters---and cluster sizes are large, the CCV and TSCB can be applied and can produce standard errors considerably smaller than the usual clustered standard errors. 

Another reason to cluster standard errors is cluster sampling. The case with $q_k$ close to zero is sometimes relevant, especially when the sample is a panel data on individuals or a cross-section of families, and the individuals or families in the sample are a small fraction of the population. Then, the clustered variance  estimator of the least squares estimator is asymptotically correct regardless of whether the treatment assignment is clustered or not. The same result holds when clusters are large ({\it e.g.}, states), $q_k$ is a substantial fraction of the clusters in the population, but $p_k$ is small---so the sample includes only a small number of units from each cluster. In other cases, cluster standard errors can be considerably larger than necessary. If cluster sizes are large and there is treatment variation within clusters, CCV and TSCB can substantially reduce the magnitude of standard errors.  

The insights in this article are relevant in other common settings of empirical economics. Consider a setting with unit-level panel data on outcomes and a treatment that is implemented on the same period for all units in the treatment group. The difference-in-difference estimator is in this case equal to the coefficient on the treatment variable in a cross-sectional regression of the change in unit-level average outcomes between the post-treatment and the pre-treatment periods on a constant and a treatment indicator that takes value one if the unit belongs to the treatment group. If treatment assignment is random across units, and the sample includes a random subset of the population or the entire population, robust standard error provide inference that is generally conservative if the sample is large relative to the population and treatment effects are heterogeneous. Here too, the methods in \cite{abadie2020sampling} can be applied to correct the bias of robust standard errors. With clustered assignment, one should cluster the standard errors at the level of assignment---for example, cluster at the village level if all farmers are assigned the same treatment status. Adding group-level fixed effects to this regression allows for group-specific linear trends in the underlying potential outcomes series but does not change the answer to the question whether one needs to adjust for clustering. Under partially clustered assignment, CCV and TSCB standard errors can continue to provide substantial improvements over conventional cluster standard errors for the fixed effect estimator.  

\section{Conclusion}

This article proposes a research framework aimed to address a question of central relevance for empirical practice: when and how we should cluster standard errors. Like in \cite{abadie2020sampling}, we shift the attention from  estimation of features of a data-generating process (i.e., infinite superpopulation) to estimation of average treatment effects of the finite population at hand. We show that, in this framework, the decision on when and how to cluster standard errors depends on the nature of the sampling and the assignment processes only, and not on the presence of within-cluster error components in the outcome variable. We derive expressions of the large sample variances of the OLS and FE estimators of the average treatment effect for a setting with clustered sampling and where assignment is random within clusters with assignment probabilities that may vary across clusters. For this setting, we demonstrate that robust standard errors can be too small and conventional cluster standard errors can be unnecessarily large. We propose two novel procedures, CCV and TSCB, that can be used to calculate more precise standard errors in settings with large clusters and where there is enough variation in treatment assignment within cluster (so that average treatment effects within clusters can be precisely estimated). While CCV and TSCB are designed for this particular setting, the general principles of the framework remain valid for other settings and estimators. If sampling is not clustered, standard errors should be clustered at the treatment assignment level because the estimand of interest depends on potential outcomes and the sampling of potential outcomes is determined only by the assignment mechanism. When the fraction of sampled clusters is non-negligible and there is variation in average treatment effects across clusters, conventional clustered standard errors may be off, and we provide an analytical framework that can be applied to derive appropriate standard errors. When sampling and assignment are random, clustering standard errors is not appropriate regardless of the structure of the covariance of the outcomes across the units in the population. In this setting, if there is substantial treatment effect heterogeneity and the sample represents a large fraction of the population of interest, robust standard errors are conservative in large samples. This bias can be corrected using the methods in \cite{abadie2020sampling}. Deriving standard error formulas for sampling and assignment processes other than the ones featured in this article is an important avenue for future research. \cite{rambachan2022design} is a recent contribution in this direction. In addition, in the present article we have restricted the analysis to linear estimators (least squares and fixed effects). \citet{xu2019asymptotic} extends the ideas and framework of this article to analyze the distribution of non-linear estimators. 

\bibliography{references}

\newpage

\appendix

\renewcommand{\thesection}{A.\arabic{section}}
\setcounter{section}{0}
\renewcommand{\thesubsection}{A.\arabic{section}.\arabic{subsection}}
\setcounter{page}{1}

\begin{center}
{\bf On-Line Appendix}\\ \medskip
{\bf When Should You Adjust Standard Errors for Clustering?}\\ \medskip
Alberto Abadie, Susan Athey, Guido W. Imbens, and Jeffrey M.\ Wooldridge\\ \medskip
Current version:~\Filemodtoday{appendix}
\end{center}
\section{Setting and notation}
\label{setting and notation}

We have a sequence of populations indexed by $k$. The $k$-th population has $n_k$ units, indexed by $i=1,\ldots, n_k$. The population is partitioned into $m_k$ strata or clusters. Let $m_{k,i}\in \{1,\ldots, m_k\}$ denote the stratum that unit $i$ of population $k$ belongs to. The number of units in cluster $m$ of population $k$ is $n_{k,m}\geq 1$. For each unit, $i$, there are two potential outcomes, $y_{k,i}(1)$ and $y_{k,i}(0)$, corresponding to treatment and no treatment. The parameter of interest is the population average treatment effect
\[
\tau_k = \frac{1}{n_k}\sum_{i=1}^{n_k} (y_{k,i}(1)-y_{k,i}(0)).
\]
The population treatment effect by cluster is
\[
\tau_{k,m} = \frac{1}{n_{k,m}}\sum_{i=1}^{n_k} 1\{m_{k,i}=m\}(y_{k,i}(1)-y_{k,i}(0)).
\]
Therefore, 
\[
\tau_k = \sum_{m=1}^{m_k}\frac{n_{k,m}}{n_k}\tau_{k,m}.
\]
We will assume  that potential outcomes, $y_{k,i}(1)$ and $y_{k,i}(0)$, are bounded in absolute value, uniformly for all $(k,i)$.

We next describe the two components of the stochastic nature of the sample. There is a stochastic binary treatment for each unit in each population, $W_{k,i}\in\{0,1\}$. The realized outcome for unit $i$ in population $k$ is $Y_{k,i}=y_{k,i}(W_{k,i})$. For a random sample of  the population, we observe the triple $(Y_{k,i},W_{k,i},m_{k,i})$. Inclusion in the sample is represented by the random variable $R_{k,i}$, which takes value one if unit $i$ belongs to the sample, and value zero if not. 

The sampling process that determines the values of $R_{k,i}$ is independent of the potential outcomes and the assignments. It consists of two stages. First, clusters are sampled with cluster sampling probability $q_k\in (0,1]$. Second, units are sampled from the subpopulation consisting of all the sampled clusters, with unit sampling probability equal to $p_k\in (0,1]$. Both $q_k$ and $p_k$ may be equal to one, or close to zero. If $q_k=1$, we sample all clusters. If $p_k=1$, we sample all units from the sampled clusters. If $q_k=p_k=1$, all units in the population are sampled.

The assignment process that determines the values of $W_{k,i}$ also consists of two stages. In the first stage, for cluster $m$ in population $k$, an assignment probability $A_{k,m}\in[0,1]$ is drawn randomly from a distribution with mean $\mu_k$, bounded away from zero and one uniformly in $k$, and variance $\sigma^2_k$, independently for each cluster. The variance $\sigma^2_k$ is key. If it is zero, we have random assignment across clusters. For positive values of $\sigma_k^2$ we have correlated assignment within the clusters. Because $A_{k,m}^2\leq A_{k,m}$, it follows that $\sigma_k^2$ is bounded above by $\mu_k(1-\mu_k)$ and that the bound is attained when $A_{k,m}$ can only take values zero or one (so all units within a  cluster have the same values for the treatment). In the second stage, each unit in cluster $m$ is assigned to the treatment independently, with cluster-specific probability $A_{k,m}$.

\section{Base case:~Difference in means}
\label{asection:base_case}

Let 
\[
N_{k,1}=\sum_{i=1}^{n_k} R_{k,i}W_{k,i}\quad\mbox{ and }\quad N_{k,0}=\sum_{i=1}^{n_k} R_{k,i}(1-W_{k,i})
\]
be the number of treated and untreated units in the sample, respectively. The total sample size is $N_k=N_{k,1}+N_{k,0}$.
We consider the simple difference of means between treated and non-treated, which is obtained as the coefficient on the treatment indicator in a regression of the outcome on a constant and the treatment,
\[
\widehat\tau_k = \frac{1}{N_{k,1}\vee 1}\sum_{i=1}^{n_k} R_{k,i}W_{k,i} Y_{k,i}
-
\frac{1}{N_{k,0}\vee 1}\sum_{i=1}^{n_k} R_{k,i}(1-W_{k,i}) Y_{k,i}.
\]
We make the following assumptions about the sampling process and the cluster sizes: {\em (i)} $q_km_k\rightarrow \infty$, {\em (ii)} $\liminf_{k\rightarrow\infty}p_k \min_m n_{k,m}>0$, and {\em (iii)} $\limsup_{k\rightarrow\infty} \max_m n_{k,m}/\min_m n_{k,m}<\infty$. The first assumption implies that the expected number of sampled clusters goes to infinity as $k$ increases. The second assumption implies that the average number of observations sampled per cluster, conditional on the cluster being sampled, does not go to zero. The third assumption restricts the imbalance between the number of units across clusters. Notice that assumptions {\em (i)} and {\em (ii)} imply $n_kp_kq_k\rightarrow\infty$, so the sample size becomes larger in expectation as $k$ increases. 

\subsection{Large $k$ distribution}

Let $\alpha_k = (1/n_k)\sum_{i=1}^{n_k} y_{k,i}(0)$ and $\tau_k = (1/n_k)\sum_{i=1}^{n_k} (y_{k,i}(1)-y_{i,k}(0))$, $u_{k,i}(1)=y_{k,i}(1)-(\alpha_k+\tau_k)$, and $u_{k,i}(0)=y_{k,i}(0)-\alpha_k$. Notice that,
\[
\sum_{i=1}^{n_k} u_{k,i}(1) =\sum_{i=1}^{n_k} u_{k,i}(0) =0.
\]
This implies
\begin{align*}
\sqrt{n_kp_kq_k}(\widehat\tau_k-\tau_k) &=\frac{b_{k,1}}{\widehat b_{k,1}}\widehat a_{k,1}
-\frac{b_{k,0}}{\widehat b_{k,0}}\widehat a_{k,0},
\end{align*}
where
\begin{align*}
\widehat a_{k,1} &= \frac{1}{\sqrt{n_kp_kq_k}\mu_k}\sum_{i=1}^{n_k} (R_{k,i}W_{k,i} -p_kq_k\mu_k)u_{k,i}(1),\\
\widehat a_{k,0} &= \frac{1}{\sqrt{n_kp_kq_k}(1-\mu_k)}\sum_{i=1}^{n_k} (R_{k,i}(1-W_{k,i}) -p_kq_k(1-\mu_k))u_{k,i}(0),
\end{align*}
$\widehat b_{k,1}=(N_{k,1}\vee 1)/n_k$, $\widehat b_{k,0}=(N_{k,0}\vee 1)/n_k$, $b_{k,1}=p_kq_k\mu_k$ and $b_{k,0}=p_kq_k(1-\mu_k)$.
We will first derive the large sample distribution of 
\begin{align*}
\widehat a_k&=\widehat a_{k,1}-\widehat a_{k,0}\\
&= \sum_{m=1}^{m_k} \big(\xi_{k,m,1}-\xi_{k,m,0}\big),
\end{align*}
where
\[
\xi_{k,m,1}=\frac{1}{\sqrt{n_kp_kq_k}\mu_k}\sum_{i=1}^{n_k}1\{m_{k,i}=m\}\big(R_{k,i}W_{k,i}-p_kq_k\mu_k\big)u_{k,i}(1),
\]
and 
\[
\xi_{k,m,0}=\frac{1}{\sqrt{n_kp_kq_k}(1-\mu_k)}\sum_{i=1}^{n_k}1\{m_{k,i}=m\}\big(R_{k,i}(1-W_{k,i})-p_kq_k(1-\mu_k)\big)u_{k,i}(0).
\]
Notice that $E[\xi_{k,m,1}]=E[\xi_{k,m,0}]=0$. Moreover, notice that the terms $\xi_{k,m,1}-\xi_{k,m,0}$ are independent across clusters, $m$. In addition,
\begin{align*}
E[\xi_{k,m,1}^2]=\frac{1}{n_k}\sum_{i=1}^{n_k}&1\{m_{k,i}=m\}\frac{1-p_kq_k\mu_k}{\mu_k}u_{k,i}^2(1)\\&+\frac{2}{n_k}\sum_{i=1}^{n_k-1}
\sum_{j=i+1}^{n_k }1\{m_{k,i}=m_{k,j}=m\}\frac{p_k(\sigma_k^2+\mu_k^2(1-q_k))}{\mu_k^2}u_{k,i}(1)u_{k,j}(1).
\end{align*}
\begin{align*}
E[\xi_{k,m,0}^2]=\frac{1}{n_k}\sum_{i=1}^{n_k }&1\{m_{k,i}=m\}\frac{1-p_kq_k(1-\mu_k)}{(1-\mu_k)}u_{k,i}^2(0)\\&+\frac{2}{n_k}\sum_{i=1}^{n_k-1}
\sum_{j=i+1}^{n_k }1\{m_{k,i}=m_{k,j}=m\}\frac{p_k(\sigma_k^2+(1-\mu_k)^2(1-q_k)) }{(1-\mu_k)^2}u_{k,i}(0)u_{k,j}(0),
\end{align*}
and 
\begin{align*}
E[&\xi_{k,m,1}\xi_{k,m,0}]=-\frac{1}{n_k}\sum_{i=1}^{n_k}1\{m_{k,i}=m\}p_kq_k u_{k,i}(1)u_{k,i}(0)\\
&+\frac{1}{n_k}\sum_{i=1}^{n_k-1}\sum_{j=i+1}^{n_k }1\{m_{k,i}=m_{k,j}=m\}\frac{p_k(\mu_k(1-\mu_k)(1-q_k)-\sigma_k^2)}{\mu_k(1-\mu_k)}\big(u_{k,i}(0)u_{k,j}(1) + u_{k,i}(1)u_{k,j}(0)\big).
\end{align*}
We obtain:
\begin{align*}
&n_kE[(\xi_{k,m,1}-\xi_{k,m,0})^2]\\&=\frac{1}{\mu_k}\sum_{i=1}^{n_k} 1\{m_{k,i}=m\}u^2_{k,i}(1)
+\frac{1}{1-\mu_k}\sum_{i=1}^{n_k} 1\{m_{k,i}=m\}u^2_{k,i}(0)\\
&+2p_k\sum_{i=1}^{n_k-1}\sum_{j=i+1}^{n_k }1\{m_{k,i}=m_{k,j}=m\}\Big(u_{k,i}(1)u_{k,j}(1)+u_{k,i}(0)u_{k,j}(0) - u_{k,i}(0)u_{k,j}(1)-u_{k,i}(1)u_{k,j}(0)\Big)\\
&-p_kq_k\Bigg(\sum_{i=1}^{n_k} 1\{m_{k,i}=m\}\big(u^2_{k,i}(1)+u^2_{k,i}(0)-2u_{k,i}(1)u_{k,i}(0)\big)\\
&\hspace*{1cm}+2\sum_{i=1}^{n_k-1}\sum_{j=i+1}^{n_k} 1\{m_{k,i}=m_{k,j}=m\}\big(u_{k,i}(1)u_{k,j}(1)
+u_{k,i}(0)u_{k,j}(0)-u_{k,i}(0)u_{k,j}(1)-u_{k,i}(1)u_{k,j}(0)\big)\Bigg)\\
&+2p_k\sigma^2_k\Bigg(\sum_{i=1}^{n_k-1}\sum_{j=i+1}^{n_k} 1\{m_{k,i}=m_{k,j}=m\}\bigg(\frac{u_{k,i}(1)u_{k,j}(1)}{\mu^2_k}+\frac{u_{k,i}(0)u_{k,j}(0)}{(1-\mu_k)^2}+\frac{u_{k,i}(0)u_{k,j}(1)}{\mu_k(1-\mu_k)}+\frac{u_{k,i}(1)u_{k,j}(0)}{\mu_k(1-\mu_k)}\Bigg).
\end{align*}
Therefore,
\begin{align*}
n_k&E[(\xi_{k,m,1}-\xi_{k,m,0})^2]\\&=\frac{1}{\mu_k}\sum_{i=1}^{n_k} 1\{m_{k,i}=m\}u^2_{k,i}(1)
+\frac{1}{1-\mu_k}\sum_{i=1}^{n_k} 1\{m_{k,i}=m\}u^2_{k,i}(0)\\
&+p_k\Bigg[\Bigg(\sum_{i=1}^{n_k}1\{m_{k,i}=m\}\big(u_{k,i}(1)-u_{k,i}(0)\big)\Bigg)^2-\sum_{i=1}^{n_k}1\{m_{k,i}=m\}\big(u_{k,i}(1)-u_{k,i}(0)\big)^2\Bigg]\\
&-p_kq_k\Bigg(\sum_{i=1}^{n_k} 1\{m_{k,i}=m\}\big(u_{k,i}(1)-u_{k,i}(0)\big)\Bigg)^2\\
&+p_k\sigma^2_k\Bigg[\Bigg(\sum_{i=1}^{n_k}1\{m_{k,i}=m\}\bigg(\frac{u_{k,i}(1)}{\mu_k}+\frac{u_{k,i}(0)}{1-\mu_k}\bigg)\Bigg)^2-\sum_{i=1}^{n_k}1\{m_{k,i}=m\}\bigg(\frac{u_{k,i}(1)}{\mu_k}+\frac{u_{k,i}(0)}{1-\mu_k}\bigg)^2\Bigg].
\end{align*}
Let $v_k=\sum_{m=1}^{m_k}E[(\xi_{k,m,1}-\xi_{k,m,0})^2]$, then
\begin{align*}
n_kv_k&=\sum_{i=1}^{n_k} \bigg(\frac{u^2_{k,i}(1)}{\mu_k}+\frac{u^2_{k,i}(0)}{1-\mu_k}\bigg)\\
&+p_k\sum_{m=1}^{m_k}\Bigg[\Bigg(\sum_{i=1}^{n_k}1\{m_{k,i}=m\}\big(u_{k,i}(1)-u_{k,i}(0)\big)\Bigg)^2-\sum_{i=1}^{n_k}1\{m_{k,i}=m\}\big(u_{k,i}(1)-u_{k,i}(0)\big)^2\Bigg]\\
&-p_kq_k\sum_{m=1}^{m_k}\Bigg(\sum_{i=1}^{n_k} 1\{m_{k,i}=m\}\big(u_{k,i}(1)-u_{k,i}(0)\big)\Bigg)^2\\
&+p_k\sigma^2_k\sum_{m=1}^{m_k}\Bigg[\Bigg(\sum_{i=1}^{n_k}1\{m_{k,i}=m\}\bigg(\frac{u_{k,i}(1)}{\mu_k}+\frac{u_{k,i}(0)}{1-\mu_k}\bigg)\Bigg)^2-\sum_{i=1}^{n_k}1\{m_{k,i}=m\}\bigg(\frac{u_{k,i}(1)}{\mu_k}+\frac{u_{k,i}(0)}{1-\mu_k}\bigg)^2\Bigg].
\end{align*}
Alternatively, we can write this expression as
\begin{align*}
n_kv_k&=\sum_{i=1}^{n_k} \bigg(\frac{u^2_{k,i}(1)}{\mu_k}+\frac{u^2_{k,i}(0)}{1-\mu_k}\bigg)\\
&-p_k\sum_{i=1}^{n_k}\big(u_{k,i}(1)-u_{k,i}(0)\big)^2-p_k\sigma_k^2\sum_{i=1}^{n_k}\bigg(\frac{u_{k,i}(1)}{\mu_k}+\frac{u_{k,i}(0)}{1-\mu_k}\bigg)^2\\
&+p_k(1-q_k)\sum_{m=1}^{m_k}\Bigg(\sum_{i=1}^{n_k} 1\{m_{k,i}=m\}\big(u_{k,i}(1)-u_{k,i}(0)\big)\Bigg)^2\\
&+p_k\sigma^2_k\sum_{m=1}^{m_k}\Bigg(\sum_{i=1}^{n_k}1\{m_{k,i}=m\}\bigg(\frac{u_{k,i}(1)}{\mu_k}+\frac{u_{k,i}(0)}{1-\mu_k}\bigg)\Bigg)^2.
\end{align*}
The sum of the first three terms is minimized for $p_k=1$ and $\sigma_k^2=\mu_k(1-\mu_k)$, in which case this sum is equal to zero. Therefore,
\begin{align}
\label{equation:vbound}
v_k&\geq (p_k\min_m n_{k,m})(1-q_k)\sum_{m=1}^{m_k}\frac{n_{k,m}}{n_k}\Bigg(\frac{1}{n_{k,m}}\sum_{i=1}^{n_k}1\{m_{k,i}=m\}\big(u_{k,i}(1)-u_{k,i}(0)\big)\Bigg)^2\nonumber\\
&+(p_k\min_m n_{k,m})\,\sigma^2_k\sum_{m=1}^{m_k}\frac{n_{k,m}}{n_k}\Bigg(\frac{1}{n_{k,m}}\sum_{i=1}^{n_k}1\{m_{k,i}=m\}\bigg(\frac{u_{k,i}(1)}{\mu_k}+\frac{u_{k,i}(0)}{1-\mu_k}\bigg)\Bigg)^2.
\end{align}
We will assume that $\liminf_{k\rightarrow\infty} ((1-q_k)\vee \sigma_k^2)>0$, so either sampling or assignment or both are correlated within cluster. (We study the case $q_k=1$ and $\sigma_k^2=0$ separately below.) In addition, assume {\em (i)} $\liminf_{k\rightarrow\infty} (1-q_k)>0$ and
\begin{equation}
\label{equation:tehet}
\liminf\limits_{k\rightarrow\infty}\sum_{m=1}^{m_k}\frac{n_{k,m}}{n_k}\Bigg(\frac{1}{n_{k,m}}\sum_{i=1}^{n_k}1\{m_{k,i}=m\}\big(u_{k,i}(1)-u_{k,i}(0)\big)\Bigg)^2>0,
\end{equation}
or {\em (ii)} $\liminf_{k\rightarrow\infty} \sigma_k^2>0$ and
\begin{equation}
\label{equation:pohet}
\liminf\limits_{k\rightarrow\infty}\sum_{m=1}^{m_k}\frac{n_{k,m}}{n_k}\Bigg(\frac{1}{n_{k,m}}\sum_{i=1}^{n_k}1\{m_{k,i}=m\}\bigg(\frac{u_{k,i}(1)}{\mu_k}+\frac{u_{k,i}(0)}{1-\mu_k}\bigg)\Bigg)^2>0.
\end{equation}
Equation (\ref{equation:tehet}) would be violated if, as $k$ increases, there is no variation in average treatment effects across clusters. Equation (\ref{equation:pohet}) would be violated if as $k$ increases there is no variation in average potential outcomes across clusters. If equations (\ref{equation:tehet}) and (\ref{equation:pohet}) hold, $v_k$ is bounded below by a term of order at least $p_k\min_m n_{k,m}$.
Recall our assumption, $\liminf_{k\rightarrow\infty}p_k \min_m n_{k,m}>0$, so the average number of observations sampled per cluster, conditional on the cluster being sampled, does not go to zero. Then,
\[
\liminf\limits_{k\rightarrow\infty} v_k>0.
\]

To obtain a CLT, we will check Lyapunov's condition,
\[
\lim_{k\rightarrow\infty }\sum_{m=1}^{m_k}\frac{1}{v_k^{1+\delta/2}}E[|\xi_{k,m,1}-\xi_{k,m,0}|^{2+\delta}] =0,
\]
for some $\delta>0$. Because potential outcomes are uniformly bounded and $\mu_k$ is uniformly bounded away from zero, we obtain 
\begin{align*}
|\xi_{k,m,1}|^{2+\delta}&\leq c\frac{n_{k,m}^{2+\delta}}{(n_kp_kq_k)^{1+\delta/2}}
\left|\frac{1}{n_{k,m}}\sum_{i=1}^{n_k}1\{m_{k,i}=m\} |R_{k,i}W_{k,i}-p_kq_k\mu_k|\right|^{2+\delta},
\end{align*}
where $c$ is some generic positive constant, whose value may change across equations. Consider $\delta=1$, and let
\begin{align*}
S&_{k,m,1}^3
= E\left[
\left|\frac{1}{n_{k,m}}\sum_{i=1}^{n_k}1\{m_{k,i}=m\} |R_{k,i}W_{k,i}-p_kq_k\mu_k|\right|^3\right]\\
&\leq \frac{1}{n_{k,m}^3}n_{k,m}E[|R_{k,i}W_{k,i}-p_kq_k\mu_k|^3]\\
&+\frac{3}{n_{k,m}^3}n_{k,m}(n_{k,m}-1)E\big[|R_{k,i}W_{k,i}-p_kq_k\mu_k|^2|R_{k,j}W_{k,j}-p_kq_k\mu_k|\big|m_{k,i}=m_{k,j}=m\big]\\
&+\frac{6}{n_{k,m}^3}\dbinom{n_{k,m}}{3}E\big[|R_{k,i}W_{k,i}-p_kq_k\mu_k||R_{k,j}W_{k,j}-p_kq_k\mu_k||R_{k,t}W_{k,t}-p_kq_k\mu_k|\big|m_{k,i}=m_{k,j}=m_{k,t}=m\big],
\end{align*}
for $i\neq j\neq t$. (The second and third terms on the left-hand side of last equation only appear when $n_{k,m}\geq 2$ and $n_{k,m}\geq 3$, respectively) As a result,
\begin{align*}
S_{k,m,1}^3
&\leq c\left(\frac{p_kq_k}{n_{k,m}^2}+ \frac{p_k^2q_k}{n_{k,m}}+p_k^3q_k\right)\\
&\leq c\, p_k^3q_k\left(\frac{1}{p_k^2\min_m n_{k,m}^2}+ \frac{1}{p_k\min_m n_{k,m}}+1\right).
\end{align*}
Because $\liminf_{k\rightarrow\infty} p_k\min_m n_{k,m}>0$, for large enough $k$ we obtain,
\begin{align*}
E[|\xi_{k,m,1}|^3]&\leq c\frac{p_k^3q_kn_{k,m}^3}{(n_kp_kq_k)^{3/2}},
\end{align*}
and the same bound applies for $E[|\xi_{k,m,0}|^3]$. Notice that 
\begin{align*}
\sum_{m=1}^{m_k} E[|\xi_{k,m,1}-\xi_{k,m,0}|^3]&\leq \sum_{m=1}^{m_k} E[(|\xi_{k,m,1}|+|\xi_{k,m,0}|)^3]\\&=\sum_{m=1}^{m_k} E[|\xi_{k,m,1}|^3] +
\sum_{m=1}^{m_k} E[|\xi_{k,m,0}|^3]\\
&+3\sum_{m=1}^{m_k} E[|\xi_{k,m,1}|^2|\xi_{k,m,0}|]+3\sum_{m=1}^{m_k} E[|\xi_{k,m,1}||\xi_{k,m,0}|^2].
\end{align*}
Now, H\"{o}lder's inequality implies that
\begin{equation}
\label{equation:lyapsuff}
\frac{p_k^3q_k\sum_{m=1}^{m_k}n_{k,m}^3}{v_k^{3/2}(n_kp_kq_k)^{3/2}}\longrightarrow 0,
\end{equation}
is sufficient for the Lyapunov condition to hold. 
Because $\max_m n_{k,m}/\min_m n_{k,m}$ is bounded asymptotically, we obtain, 
\begin{align*}
\limsup\limits_{k\rightarrow\infty}\frac{p_k^3q_k\sum_{m=1}^{m_k}n_{k,m}^{3}}{v_k^{3/2}(n_kp_kq_k)^{3/2}}
&\leq \limsup\limits_{k\rightarrow\infty}c\,\frac{p_k^3q_km_k\max_m n_{k,m}^{3}}{(p_k^2q_km_k\min_m n_{k,m}^2)^{3/2}}\\
&\leq \limsup\limits_{k\rightarrow\infty} \left(\frac{\max_m n_{k,m}}{\min_m n_{k,m}}\right)^3\frac{c}{\sqrt{q_km_k}}= 0,
\end{align*}
and so the Lyapunov condition holds.
As a result, we obtain
\[
\widehat a_k/\sqrt{v_k}\stackrel{d}{\longrightarrow} N(0,1).
\]
We will next prove that both $\widehat a_{k,1}/\sqrt{v_k}$ and $\widehat a_{k,0}/\sqrt{v_k}$ are $\mathcal O_p(1)$. 
\begin{align*}
E[\widehat a_{k,1}^2]&=\frac{1}{n_kp_kq_k}\frac{1}{\mu_k^2}
\sum_{m=1}^{m_k} E\left[\left(\sum_{i=1}^{n_k} 1\{m_{k,i}=m\}(R_{k,i}W_{k,i}-p_kq_k\mu_k)u_{k,i}(1)\right)^2\right]\\
&\leq c\,\frac{1}{n_kp_kq_k} \sum_{m=1}^{m_k} \big(n_{k,m}p_kq_k+n_{k,m}(n_{k,m}-1)p_k^2q_k\big)\\
&=c\, \Big(1+\sum_{m=1}^{m_k} \frac{n_{k,m}(n_{k,m}-1)p_k}{n_k}\Big).
\end{align*}
Therefore,
\begin{align*}
E[(\widehat a_{k,1}/\sqrt{v_k})^2]&\leq c \left(\frac{1}{p_k \min_m n_{k,m}}+\sum_{m=1}^{m_k} \frac{(\max_m n_{k,m})(n_{k,m}-1)p_k}{n_k p_k\min_m n_{k,m}}\right).
\end{align*}
Because $\limsup \max_m n_{k,m}/\min_m n_{k,m}<\infty$, we obtain 
$\limsup_{k\rightarrow\infty} E[(\widehat a_{k,1}/v_k)^2] <\infty$. As a result, $\widehat a_{k,1}/\sqrt{v_k}$ is $\mathcal O_p(1)$.

Let $\widetilde b_{k,1}=N_{k,1}/n_k$. Consider $k$ large enough, so $p_k\min_m n_{k,m}$ is bounded away from zero, making $\widetilde b_{k,1}/b_{k,1}$ well-defined. Notice that $E[\widetilde b_{k,1}/b_{k,1}]=1$ and 
\begin{align*}
\mbox{var}(\widetilde b_{k,1}/b_{k,1}) &= \frac{1}{(n_kp_kq_k\mu_k)^2}\sum_{m=1}^{m_k}
E\left[\left(\sum_{i=1}^{n_k} 1\{m_{k,i}=m\}(R_{k,i}W_{k,i}-n_kp_kq_k\mu_k)\right)^2\right]\\
&=\frac{n_kp_kq_k\mu_k(1-p_kq_k\mu_k)}{(n_kp_kq_k\mu_k)^2}
+\sum_{m=1}^{m_k} \frac{n_{k,m}(n_{k,m}-1)p_k^2q_k(\sigma_k^2+(1-q_k)\mu_k^2)}{(n_kp_kq_k\mu_k)^2}\\
&\leq \frac{1-p_kq_k\mu_k}{n_kp_kq_k\mu_k}
+c\,\frac{n_k(\max_m n_{k,m}-1)p_k^2q_k}{(n_kp_kq_k)^2}\\
&\leq \frac{1-p_kq_k\mu_k}{n_kp_kq_k\mu_k}
+c\,\frac{(\max_m n_{k,m}-1)}{\min n_{k,m}}\frac{1}{q_km_k}
\longrightarrow 0. 
\end{align*}
This implies $\widetilde b_{k,1}/b_{k,1}\stackrel{p}{\rightarrow} 1$. Analogous calculations yield $\widetilde b_{0,k}/b_{0,k}\stackrel{p}{\rightarrow} 1$. For large enough $k$, $\widetilde b_{k,1}/b_{k,1}=0$ if and only if $N_{k,1}=0$, which implies $\Pr(N_{k,1}=0)\rightarrow 0$. It follows that, for large enough $k$,
\[
\Pr(|\widetilde b_{k,1}/b_{k,1} -\widehat b_{k,1}/b_{k,1}|=0)=\Pr(N_{k,1}>0)\longrightarrow 1
\]
and $\widehat b_{k,1}/b_{k,1}\stackrel{p}{\rightarrow} 1$.
Using analogous calculations, we obtain $\widehat b_{k,0}/b_{k,0}\stackrel{p}{\rightarrow} 1$.
As a result,
\begin{align*}
\sqrt{n_kp_kq_k}(\widehat\tau_k-\tau_k)/v_k^{1/2} &= \frac{b_{k,1}}{\widehat b_{k,1}}\frac{\widehat a_{k,1}}{v_k^{1/2}} - \frac{b_{k,0}}{\widehat b_{k,0}}\frac{\widehat a_{k,0}}{v_k^{1/2}}\\
&=\frac{\widehat a_k}{v_k^{1/2}} + \left(\frac{b_{k,1}}{\widehat b_{k,1}}-1\right)\frac{\widehat a_{k,1}}{v_k^{1/2}} - \left(\frac{b_{k,0}}{\widehat b_{k,0}}-1\right)\frac{\widehat a_{k,0}}{v_k^{1/2}}\\\
&=\widehat a_k/\sqrt{v_k} + \scaleto{\mathcal{O}}{5pt}_p(1).
\end{align*}
Therefore, 
\[
\sqrt{n_kp_kq_k}(\widehat\tau_k-\tau_k)/v_k^{1/2}\stackrel{d}{\longrightarrow} N(0,1).
\]
Using $\widetilde b_{1,k}/b_{1,k}\stackrel{p}{\rightarrow} 1$ and $\widetilde b_{0,k}/b_{0,k}\stackrel{p}{\rightarrow} 1$, it is easy to show $N_k/(n_kp_kq_k)\stackrel{p}{\rightarrow}1$, which implies
\[
\sqrt{N_k}(\widehat\tau_k-\tau_k)/v_k^{1/2}\stackrel{d}{\longrightarrow} N(0,1).
\]

We will next consider the case of $q_k=1$ and $\sigma_k^2=0$, where no clustering is required. Consider 
\[
\vartheta_{k,i,1}=\frac{1}{\sqrt{n_kp_k}\mu_k} \big(R_{k,i}W_{k,i}-p_k\mu_k\big)u_{k,i}(1)
\]
and 
\[
\vartheta_{k,i,0}=\frac{1}{\sqrt{n_kp_k}(1-\mu_k)} \big(R_{k,i}(1-W_{k,i})-p_k(1-\mu_k)\big)u_{k,i}(0).
\]
Redefine now $v_k=\sum_{i=1}^{n_k} E\big[(\vartheta_{k,i,1}-\vartheta_{k,i,0})^2\big]$. Then,
\begin{align*}
v_k&= \frac{1}{n_k}\sum_{i=1}^{n_k} \bigg(\frac{u^2_{k,i}(1)}{\mu_k}+\frac{u^2_{k,i}(0)}{1-\mu_k}\bigg)- p_k \frac{1}{n_k}\sum_{i=1}^{n_k} \big(u_{k,i}(1)-u_{k,i}(0)\big)^2.
\end{align*}
Notice that $v_k$ is minimized for $p_k=1$, in which case
\begin{align*}
v_k&=\frac{1}{n_k}\sum_{i=1}^{n_k} \bigg(\frac{u^2_{k,i}(1)}{\mu_k}+\frac{u^2_{k,i}(0)}{1-\mu_k}\bigg)-\frac{1}{n_k}\sum_{i=1}^{n_k} \big(u_{k,i}(1)-u_{k,i}(0)\big)^2\\
&=\frac{1}{n_k}\sum_{i=1}^{n_k} \bigg(\frac{1-\mu_k}{\mu_k}u^2_{k,i}(1)+\frac{\mu_k}{1-\mu_k}u^2_{k,i}(0)+2u_{k,i}(1)u_{k,i}(0)\bigg)\\
&=\mu_k(1-\mu_k)\frac{1}{n_k}\sum_{i=1}^{n_k} \bigg(\frac{u^2_{k,i}(1)}{\mu^2_k}+\frac{u^2_{k,i}(0)}{(1-\mu_k)^2}+2\,\frac{u_{k,i}(1)u_{k,i}(0)}{\mu_k(1-\mu_k)}\bigg)\\
&=\mu_k(1-\mu_k)\frac{1}{n_k}\sum_{i=1}^{n_k} \bigg(\frac{u_{k,i}(1)}{\mu_k}+\frac{u_{k,i}(0)}{1-\mu_k}\bigg)^2.
\end{align*}
Therefore, the assumption
\[
\liminf\limits_{k\rightarrow\infty} \frac{1}{n_k}\sum_{i=1}^{n_k} \bigg(\frac{u_{k,i}(1)}{\mu_k}+\frac{u_{k,i}(0)}{1-\mu_k}\bigg)^2>0
\]
is enough for $\liminf_{k\rightarrow\infty} v_k>0$. Notice now that 
\begin{align*}
E[|\vartheta_{k,i,1}|^3]&\leq \frac{1}{(n_kp_k)^{3/2}} E[|R_{k,i}W_{k,i}-p_k\mu_k|^3]\\
&=\frac{1}{(n_kp_k)^{3/2}}(1-p_k\mu_k)^3p_k\mu_k+(p_k\mu_k)^3(1-p_k\mu_k)\\
&\leq c\frac{p_k}{(n_kp_k)^{3/2}},
\end{align*}
and the same bound holds for $E[|\vartheta_{k,i,0}|^3]$. Therefore, for the Lyapunov condition to hold, it is enough that
\[
\frac{n_kp_k}{(n_kp_k)^{3/2}}=\frac{1}{\sqrt{n_kp_k}}\longrightarrow 0,
\]
or $n_kp_k\rightarrow \infty$. That is, assumptions {\em (i)}-{\em (iii)}, which we used for the clustered case, are replaced by $n_kp_k\rightarrow \infty$.

\subsection{Estimation of the variance}

Let $\widehat U_{k,i}=Y_{k,i}-\widehat\alpha_k-\widehat\tau_kW_{k,i}$ be the residuals from the regression of $Y_{k,i}$ or a constant and $W_{k,i}$. Here, $\widehat\alpha_k$ is the coefficient on the constant regressor equal to one, and $\widehat\tau_k$ is the coefficient on $W_{k,i}$. We have already shown $v_k^{-1/2}(\widehat\tau_k-\tau_k)=\mathcal{O}_p(1/\sqrt{n_kp_kq_k})$. The same is true about $\widehat\alpha_k$ (e.g., apply the proof for $\widehat\tau_k$ after replacing each $y_{k,i}(1)$ with a zero). Define $\widehat\Sigma_k=\sum_{m=1}^{m_k} \widehat\Sigma_{k,m}$, where
\[
\widehat\Sigma_{k,m}= \left(\sum_{i=1}^{n_k}1\{m_{k,i}=m\}R_{k,i}\left(\begin{array}{c}\widehat U_{k,i}\\W_{k,i}\widehat U_{k,i}\end{array}\right)\right)\left(\sum_{i=1}^{n_k}1\{m_{k,i}=m\}R_{k,i}\left(\begin{array}{c}\widehat U_{k,i}\\W_{k,i}\widehat U_{k,i}\end{array}\right)\right)'.\]
Also, let 
\[
\widehat Q_k = \sum_{i=1}^{n_k} R_{k,i}\left(\begin{array}{c}1\\W_{k,i}\end{array}\right)\left(\begin{array}{c}1\\W_{k,i}\end{array}\right)',
\]
and $z= (0,1)'$. Then, the cluster estimator of the variance of 
$\sqrt{N_k}(\widehat\tau_k-\tau_k)$ is 
\begin{align*}
\widehat V_k^{\rm{cluster}} = N_kz'\widehat Q_k^{-1}\widehat\Sigma_k \widehat Q_k^{-1}z.
\end{align*}
Notice that
\[
(n_kp_kq_k)^{-1}E[\widehat Q_k]=\left(\begin{array}{cc}1 & \mu_k\\\mu_k&\mu_k\end{array}\right).
\]
In addition,
\[
\frac{1}{n_kp_kq_k}\widehat Q_k (2,2)= \frac{1}{n_kp_kq_k}\sum_{m=1}^{m_k}
\sum_{i=1}^{n_k} 1\{m_{k,i}=m\} R_{k,i}W_{k,i}.
\]
\begin{align*}
\mbox{var}\Bigg(\sum_{i=1}^{n_k} 1\{m_{k,i}=m\} R_{k,i}W_{k,i}\Bigg)&=n_{k,m} p_kq_k\mu_k(1-p_kq_k\mu_k)\\
&+ n_{k,m}(n_{k,m}-1)p_k^2q_k(\sigma_k^2+\mu_k^2(1-q_k)).
\end{align*}
Therefore, under conditions {\em (i)}-{\em (iii)}, we obtain
\begin{align*}
\mbox{var}\Bigg(\frac{1}{n_kp_kq_k}\widehat Q_k (2,2)\Bigg)
&\leq \frac{c}{n_kp_kq_k} \Bigg(1+p_k(\textstyle\max_m n_{k,m}-1)\Bigg)\\
&= c\, \frac{\textstyle\max_m n_{k,m}}{n_kq_k}+\scaleto{\mathcal{O}}{5pt}(1)\\
&\leq c\, \frac{\textstyle\max_m n_{k,m}}{\textstyle\min_m n_{k,m}}\,\frac{1}{q_km_k}+\scaleto{\mathcal{O}}{5pt}(1)\longrightarrow 0.
\end{align*}
Analogous calculations yield $\mbox{var}((n_kp_kq_k)^{-1}\widehat Q_k(1,1))\rightarrow 0$. Therefore, 
\[
\frac{1}{n_kp_kq_k}\widehat Q_k = \left(\begin{array}{cc}1 & \mu_k\\\mu_k&\mu_k\end{array}\right) + \scaleto{\mathcal{O}}{5pt}_p(1)
\]
and 
\[
n_kq_kp_k \widehat Q_k^{-1} =  H_k+\scaleto{\mathcal{O}}{5pt}_p(1),\quad\mbox{where}\quad H_k=\frac{1}{\mu_k(1-\mu_k)}\left(\begin{array}{rc}\mu_k & -\mu_k\\-\mu_k&\ \ 1\end{array}\right).
\]
Now, let $U_{k,i}=Y_{k,i}-\alpha_k-\tau_kW_{k,i}=W_{k,i}u_{k,i}(1)+(1-W_{k,i})u_{k,i}(0)$.  Notice that 
\[
v_k^{-1/2}\max_{i=1,\ldots, n_k} |\widehat U_{k,i}-U_{k,i}|\leq v_k^{-1/2}|\widehat\alpha_k-\alpha_k|+v_k^{-1/2}|\widehat\tau_k-\tau_k|=\mathcal{O}_p(1/\sqrt{n_kp_kq_k}).
\]
Define $\mwidebar\Sigma_k=\sum_{m=1}^{m_k} \mwidebar\Sigma_{k,m}$, where
\[
\mwidebar\Sigma_{k,m}= \left(\sum_{i=1}^{n_k}1\{m_{k,i}=m\}R_{k,i}\left(\begin{array}{c} U_{k,i}\\W_{k,i} U_{k,i}\end{array}\right)\right)\left(\sum_{i=1}^{n_k}1\{m_{k,i}=m\}R_{k,i}\left(\begin{array}{c} U_{k,i}\\W_{k,i} U_{k,i}\end{array}\right)\right)'.
\]
We will show 
\[
\frac{1}{n_kp_kq_kv_k}(\widehat\Sigma_k-\mwidebar\Sigma_k)\stackrel{p}{\longrightarrow} 0. 
\]
Notice that
\begin{align*}
\widehat\Sigma_{k,m}(2,2)-\mwidebar\Sigma_{k,m}(2,2)&= \Bigg(\sum_{i=1}^{n_k}1\{m_{k,i}=m\}R_{k,i}W_{k,i}(\widehat U_{k,i}-U_{k,i})\Bigg)^2\\
&+2\Bigg(\sum_{i=1}^{n_k}1\{m_{k,i}=m\}R_{k,i}W_{k,i}U_{k,i}\Bigg)\Bigg(\sum_{i=1}^{n_k}1\{m_{k,i}=m\}R_{k,i}W_{k,i}(\widehat U_{k,i}-U_{k,i})\Bigg).
\end{align*}
Therefore,
\begin{align*}
\frac{1}{n_kp_kq_kv_k}|\widehat\Sigma_k(2,2)-\mwidebar\Sigma_k(2,2)|
&\leq c \frac{1}{n_kp_kq_kv_k}\sum_{m=1}^{m_k}\Bigg(\sum_{i=1}^{n_k}1\{m_{k,i}=m\}R_{k,i}W_{k,i}\Bigg)^2\\
&\qquad\times\Bigg(\max_{i=1,\ldots, n_k} |\widehat U_{k,i}-U_{k,i}|^2 + \max_{i=1,\ldots, n_k} |\widehat U_{k,i}-U_{k,i}|
\Bigg).
\end{align*}
The same expression holds for the off-diagonal elements of $\widehat\Sigma_{k,m}-\mwidebar\Sigma_{k,m}$. For $\widehat\Sigma_{k,m}(1,1)-\mwidebar\Sigma_{k,m}(1,1)$, the expression holds once we replace each $W_{k,i}$ with a one. Let $\|\cdot\|$ be the Frobenius norm of a matrix. Then,
\begin{align*}
\frac{1}{n_kp_kq_kv_k}\|\widehat\Sigma_k-\mwidebar\Sigma_k\|
&\leq c \frac{1}{n_kp_kq_kv_k}\sum_{m=1}^{m_k}\Bigg(\sum_{i=1}^{n_k}1\{m_{k,i}=m\}R_{k,i}\Bigg)^2\\
&\qquad\times\Bigg(\max_{i=1,\ldots, n_k} |\widehat U_{k,i}-U_{k,i}|^2 + \max_{i=1,\ldots, n_k} |\widehat U_{k,i}-U_{k,i}|
\Bigg).
\end{align*}
We will prove that the right-hand side of the previous equation converges to zero in probability. We will factorize each term into a expression that is bounded in probability and one that converges to zero in $L_1$. 
\begin{align*}
E\Bigg[\sum_{m=1}^{m_k}\Bigg(\sum_{i=1}^{n_k}1\{m_{k,i}=m\}R_{k,i}\Bigg)^2\Bigg] 
&\leq n_kp_kq_k+n_k(\max_m n_{k,m}-1)p_k^2q_k. 
\end{align*}
For the first term, notice that
\begin{align*}
\max_{i=1,\ldots, n_k} |\widehat U_{k,i}&-U_{k,i}|^2\frac{n_kp_kq_k+n_k(\max_m n_{k,m}-1)p_k^2q_k}{n_kp_kq_kv_k}\\
&= \frac{n_kp_kq_k}{v_k}\max_{i=1,\ldots, n_k} |\widehat U_{k,i}-U_{k,i}|^2\Bigg(\frac{n_kp_kq_k+n_k(\max_m n_{k,m}-1)p_k^2q_k}{(n_kp_kq_k)^2}\Bigg)\\
&\leq\frac{n_kp_kq_k}{v_k}\max_{i=1,\ldots, n_k} |\widehat U_{k,i}-U_{k,i}|^2\Bigg(\frac{1}{n_kp_kq_k}+\frac{\max_m n_{k,m}-1}{\min_m n_{k,m}}\frac{1}{q_km_k}\Bigg)\\
&=\mathcal{O}_p(1)\,\scaleto{\mathcal{O}}{5pt}(1).
\end{align*}
For the second term, using the fact that $v_k$ is greater or equal to $ p_k\min_m n_{k,m}>0$ times a term with limit inferior that is bounded away from zero, we obtain
\begin{align*}
\max_{i=1,\ldots, n_k} |\widehat U_{k,i}&-U_{k,i}|\frac{n_kp_kq_k+n_k(\max_m n_{k,m}-1)p_k^2q_k}{n_kp_kq_kv_k}\\
&= \Big(\frac{n_kp_kq_k}{v_k}\Big)^{1/2}\!\!\max_{i=1,\ldots, n_k} |\widehat U_{k,i}-U_{k,i}|\Bigg(\frac{n_kp_kq_k+n_k(\max_m n_{k,m}-1)p_k^2q_k}{(n_kp_kq_k)^{3/2}v_k^{1/2}}\Bigg)\\
&\leq\Big(\frac{n_kp_kq_k}{v_k}\Big)^{1/2}\!\!\max_{i=1,\ldots, n_k} |\widehat U_{k,i}-U_{k,i}|\Bigg(\frac{1}{(n_kp_kq_kv_k)^{1/2}}+\frac{\max_m n_{k,m}-1}{\min_m n_{k,m}}\frac{1}{(q_km_k)^{1/2}}\Bigg)\\
&=\mathcal{O}_p(1)\,\scaleto{\mathcal{O}}{5pt}(1).
\end{align*}
As a result, we obtain 
\[
\frac{1}{n_kp_kq_kv_k}\|\widehat\Sigma_k-\mwidebar\Sigma_k\|=\scaleto{\mathcal{O}}{5pt}_p(1).
\]
Notice that
\begin{align*}
\frac{n_kp_kq_k}{v_k}&\widehat Q_k^{-1}\widehat\Sigma_k \widehat Q_k^{-1}-H_k\frac{\mwidebar\Sigma_k}{n_kp_kq_kv_k} H_k=H_k\frac{\mwidebar\Sigma_k}{n_kp_kq_kv_k}\Big(n_kp_kq_k\widehat Q_k^{-1}-H_k\Big)\\
&+\Big(n_kp_kq_k\widehat Q_k^{-1}-H_k\Big)\frac{\mwidebar\Sigma_k}{n_kp_kq_kv_k}\Big(n_kp_kq_k\widehat Q_k^{-1}\Big)+\Big(n_kp_kq_k\widehat Q_k^{-1}\Big) \frac{\widehat\Sigma_k-\mwidebar\Sigma_k}{n_kp_kq_kv_k}\Big(n_kp_kq_k\widehat Q_k^{-1}\Big). 
\end{align*}
Therefore, to show that the left-hand side of the last equation is $\scaleto{\mathcal{O}}{5pt}_p(1)$, it is only left to show that $\mwidebar\Sigma_k/(n_kp_kq_kv_k)$ is $\mathcal{O}_p(1)$. We will prove this next. Notice that
\begin{align*}
\frac{1}{n_kp_kq_kv_k}\|\mwidebar\Sigma_k\|
&\leq c \frac{1}{n_kp_kq_kv_k}\sum_{m=1}^{m_k}\Bigg(\sum_{i=1}^{n_k}1\{m_{k,i}=m\}R_{k,i}\Bigg)^2.
\end{align*}
Therefore,
\begin{align*}
E\Bigg[\frac{1}{n_kp_kq_kv_k}\|\mwidebar\Sigma_k\|\Bigg]
&\leq c \frac{1}{n_kp_kq_kv_k}\Bigg(n_kp_kq_k+n_k(\max_m n_{k,m}-1)p_k^2q_k\Bigg).
\end{align*}
Then,
\begin{align*}
E\Bigg[\frac{1}{n_kp_kq_kv_k}\|\mwidebar\Sigma_k\|\Bigg]
&\leq c \Bigg(\frac{1}{v_k}+\frac{p_k(\max_m n_{k,m}-1)}{p_k\min_m n_{k,m}}\Bigg)<\infty.
\end{align*}
We, therefore, obtain,
\[
\frac{n_kp_kq_k}{v_k}\widehat Q_k^{-1}\widehat\Sigma_k \widehat Q_k^{-1}-H_k\frac{\mwidebar\Sigma_k}{n_kp_kq_kv_k} H_k\stackrel{p}{\longrightarrow} 0.
\]
Because $N_k/(n_kp_kq_k)\stackrel{p}{\rightarrow} 1$, we obtain
\begin{align*}
\widehat V_k^{\rm{cluster}}/v_k &= z'H_k\frac{\mwidebar\Sigma_k}{n_kp_kq_kv_k} H_kz+\scaleto{\mathcal{O}}{5pt}_p(1)\\
&=\frac{1}{n_kp_kq_kv_k}\left(\frac{1}{\mu_k(1-\mu_k)}\right)^2
\sum_{m=1}^{m_k} \bigg(\sum_{i=1}^{n_k} 1\{m_{k,i}=m\}R_{k,i}(W_{k,i}-\mu_k)U_{k,i}\bigg)^2+\scaleto{\mathcal{O}}{5pt}_p(1).
\end{align*}
Recall that $U^2_{k,i}=u^2_{k,i}(1)W_{k,i}+u^2_{k,i}(0)(1-W_{k,i})$. Notice that
\begin{align*}
E\bigg[\bigg(\sum_{i=1}^{n_k}& 1\{m_{k,i}=m\}R_{k,i}(W_{k,i}-\mu_k)U_{k,i}\bigg)^2\bigg]\\&=\sum_{i=1}^{n_k} 1\{m_{k,i}=m\} p_kq_k\mu_k(1-\mu_k)\Big((1-\mu_k)u^2_{k,i}(1)+\mu_k u^2_{k,i}(0)\Big)\\
&+2\sum_{i=1}^{n_k-1}\sum_{j=i+1}^{n_k}1\{m_{k,i}=m_{k,j}=m\}p_k^2q_k\Big[(\sigma_k^2 + \mu_k^2)(1-\mu_k)^2u_{k,i}(1)u_{k,j}(1)\\
&\qquad+\mu_k(1-\mu_k)(\sigma_k^2-\mu_k(1-\mu_k))(u_{k,i}(0)u_{k,j}(1)+u_{k,i}(1)u_{k,j}(0))\\
&\qquad+(\sigma_k^2+(1-\mu_k)^2)\mu_k^2u_{k,i}(0)u_{k,j}(0)\Big].
\end{align*}
Let
\[
v_k^{\rm cluster} = \frac{1}{n_kp_kq_k}\left(\frac{1}{\mu_k(1-\mu_k)}\right)^2\sum_{m=1}^{m_k}E\bigg[\bigg(\sum_{i=1}^{n_k} 1\{m_{k,i}=m\}R_{k,i}(W_{k,i}-\mu_k)U_{k,i}\bigg)^2\bigg].
\]
Then,
\begin{align*}
n_kv_k^{\rm{cluster}}&=\sum_{i=1}^{n_k} \bigg(\frac{u^2_{k,i}(1)}{\mu_k}+\frac{u^2_{k,i}(0)}{1-\mu_k}\bigg)\\
&+p_k\sum_{m=1}^{m_k}\Bigg[\Bigg(\sum_{i=1}^{n_k}1\{m_{k,i}=m\}\big(u_{k,i}(1)-u_{k,i}(0)\big)\Bigg)^2-\sum_{i=1}^{n_k}1\{m_{k,i}=m\}\big(u_{k,i}(1)-u_{k,i}(0)\big)^2\Bigg]\\
&+p_k\sigma^2_k\sum_{m=1}^{m_k}\Bigg[\Bigg(\sum_{i=1}^{n_k}1\{m_{k,i}=m\}\bigg(\frac{u_{k,i}(1)}{\mu_k}+\frac{u_{k,i}(0)}{1-\mu_k}\bigg)\Bigg)^2-\sum_{i=1}^{n_k}1\{m_{k,i}=m\}\bigg(\frac{u_{k,i}(1)}{\mu_k}+\frac{u_{k,i}(0)}{1-\mu_k}\bigg)^2\Bigg].
\end{align*}
Alternatively, we can write
\begin{align*}
n_kv_k^{\rm cluster}&=\sum_{i=1}^{n_k} \bigg(\frac{u^2_{k,i}(1)}{\mu_k}+\frac{u^2_{k,i}(0)}{1-\mu_k}\bigg)\\
&-p_k\sum_{i=1}^{n_k}\big(u_{k,i}(1)-u_{k,i}(0)\big)^2-p_k\sigma_k^2\sum_{i=1}^{n_k}\bigg(\frac{u_{k,i}(1)}{\mu_k}+\frac{u_{k,i}(0)}{1-\mu_k}\bigg)^2\\
&+p_k\sum_{m=1}^{m_k}\Bigg(\sum_{i=1}^{n_k} 1\{m_{k,i}=m\}\big(u_{k,i}(1)-u_{k,i}(0)\big)\Bigg)^2\\
&+p_k\sigma^2_k\sum_{m=1}^{m_k}\Bigg(\sum_{i=1}^{n_k}1\{m_{k,i}=m\}\bigg(\frac{u_{k,i}(1)}{\mu_k}+\frac{u_{k,i}(0)}{1-\mu_k}\bigg)\Bigg)^2.
\end{align*}
We will next show that
\[
z'H_k\frac{\mwidebar\Sigma_k}{n_kp_kq_kv_k} H_kz -\frac{v_k^{\rm cluster}}{v_k}\stackrel{p}{\longrightarrow} 0.
\]
Given the $\mu_k(1-\mu_k)$ is bounded away from zero, by the weak law of large numbers for arrays, it is enough to show
\[
\frac{1}{(n_kp_kq_kv_k)^2}\sum_{m=1}^{m_k}
E\bigg[\bigg(\sum_{i=1}^{n_k} 1\{m_{k,i}=m\}R_{k,i}(W_{k,i}-\mu_k)U_{k,i}\bigg)^4\bigg]
\longrightarrow 0. 
\]
Applying the multinomial theorem and the fact that all moments of $W_{k,i}$ as well as all potential outcomes are bounded, we obtain:
\begin{align*}
\frac{1}{(n_kp_kq_kv_k)^2}&\sum_{m=1}^{m_k}
E\bigg[\bigg(\sum_{i=1}^{n_k} 1\{m_{k,i}=m\}R_{k,i}(W_{k,i}-\mu_k)U_{k,i}\bigg)^4\bigg]\\
&\leq \frac{c}{(n_kp_kq_kv_k)^2}\Big(n_kp_kq_k + n_kp_k^2q_k \max_m n_{k,m}+n_kp_k^3q_k\max_m n^2_{k,m} + n_kp_k^4q_k\max_m n^3_{k,m}\Big).
\end{align*}
Now, using $\limsup_{k\rightarrow\infty} \max_m n_{k,m}/\min_m n_{k,m}<\infty$,
 $\limsup_{k\rightarrow\infty} p_k \min_m n_{k,m}/v_k<\infty$, and $q_km_k\rightarrow\infty$ we obtain
\begin{align*}
\frac{1}{(n_kp_kq_kv_k)^2}&\sum_{m=1}^{m_k}
E\bigg[\bigg(\sum_{i=1}^{n_k} 1\{m_{k,i}=m\}R_{k,i}(W_{k,i}-\mu_k)U_{k,i}\bigg)^4\bigg]\\
&\leq c\left(\frac{1}{n_kp_kq_kv_k^2} + \frac{\max_m n_{k,m}}{\min_m n_{k,m}}\frac{1}{q_km_kv_k^2}+\frac{p_k\max_m n^2_{k,m}}{v_k\min_m n_{k,m}}\frac{1}{q_km_kv_k}
+\frac{p_k^2\max_m n^3_{k,m}}{v^2_k\min_m n_{k,m}}\frac{1}{q_km_k}\right)\\
&\longrightarrow 0.
\end{align*}
As a result,
\[
\frac{\widehat V_k^{\rm{cluster}}}{v_k}=\frac{v_k^{\rm cluster}}{v_k}+\scaleto{\mathcal{O}}{5pt}_p(1).
\]
The robust (sandwich) estimator of the variance of $\sqrt{N_k}(\widehat\tau_k-\tau_k)$ is given by
\[
\widehat V_k^{\rm robust}=N_kz'\widehat Q_k^{-1}\widehat\Omega_k \widehat Q_k^{-1}z.
\]
where
\[
\widehat\Omega_k = \sum_{i=1}^{n_k}R_{k,i}\left(\begin{array}{c}\widehat U_{k,i}\\W_{k,i}\widehat U_{k,i}\end{array}\right)\left(\begin{array}{c}\widehat U_{k,i}\\W_{k,i}\widehat U_{k,i}\end{array}\right)'.
\]
We will derive the limit of $\widehat V_k^{\rm robust}/v_k$. Let
\[
\mwidebar\Omega_k = \sum_{i=1}^{n_k}R_{k,i}\left(\begin{array}{c}U_{k,i}\\W_{k,i} U_{k,i}\end{array}\right)\left(\begin{array}{c}U_{k,i}\\W_{k,i}U_{k,i}\end{array}\right)'.
\]
Because potential outcomes (and $W_{k,i}$) are bounded, we obtain
\begin{align*}
\frac{1}{n_kp_kq_kv_k}\|\widehat\Omega_k-\mwidebar\Omega_k\|\leq
c
\left(\frac{1}{n_kp_kq_kv_k}\sum_{i=1}^{n_k} R_{k,i}\right) \max_{i=1,\ldots, n_k} |\widehat U_{k,i}^2-U_{k,i}^2|.
\end{align*}
Because the limsup of the expectation of the first factor (which is non-negative) is bounded and the second factor converges to zero in probability as proved above, we obtain 
\[
\frac{1}{n_kp_kq_kv_k}\|\widehat\Omega_k-\mwidebar\Omega_k\|=\scaleto{\mathcal{O}}{5pt}_p(1).
\]
Notice that
\begin{align*}
\frac{1}{n_kp_kq_kv_k}\|\mwidebar\Omega_k\|&\leq c
\left(\frac{1}{n_kp_kq_kv_k}\sum_{i=1}^{n_k} R_{k,i}\right) 
\end{align*}
Again, the limsup of the expectation of the right-hand side of this equation is non-negative and bounded. As a result, we obtain $\|\mwidebar\Omega_k\|/(n_kp_kq_kv_k)=\mathcal O_p(1)$.
\begin{align*}
\widehat V_k^{\rm{robust}}/v_k &= z'H_k\frac{\mwidebar\Omega_k}{n_kp_kq_kv_k} H_kz+\scaleto{\mathcal{O}}{5pt}_p(1)\\
&=\frac{1}{n_kp_kq_kv_k}\left(\frac{1}{\mu_k(1-\mu_k)}\right)^2 \sum_{i=1}^{n_k} R_{k,i}(W_{k,i}-\mu_k)^2U_{k,i}^2+\scaleto{\mathcal{O}}{5pt}_p(1).
\end{align*}
Notice that
\begin{align*}
E\bigg[\sum_{i=1}^{n_k}& R_{k,i}(W_{k,i}-\mu_k)^2U_{k,i}^2\bigg]=\sum_{i=1}^{n_k}  p_kq_k\mu_k(1-\mu_k)\Big((1-\mu_k)u^2_{k,i}(1)+\mu_k u^2_{k,i}(0)\Big).
\end{align*} 
Finally, notice that
\begin{align*}
\frac{1}{(n_kp_kq_kv_k)^2}\sum_{m=1}^{m_k} E\Bigg[\Bigg(\sum_{i=1}^{n_k} R_{k,i}(W_{k,i}-\mu_k)^2U_{k,i}^2\Bigg)^2\Bigg]&\leq c \frac{n_kp_kq_k+n_kp_k^2q_k\max_m n_{k,m}}{(n_kp_kq_kv_k)^2}\\
&\leq c \Bigg(\frac{1}{n_kp_kq_kv_k^2} + \frac{\max_m n_{k,m}}{\min_m n_{k,m}}\frac{1}{q_km_kv_k^2}\Bigg)\\&\longrightarrow 0.
\end{align*}
Therefore, by the weak law of large numbers for arrays, we obtain
\[
\frac{\widehat V_k^{\rm{robust}}}{v_k}=\frac{v_k^{\rm robust}}{v_k}+\scaleto{\mathcal{O}}{5pt}_p(1),
\]
where
\[
v_k^{\rm robust} = \frac{1}{n_k}\sum_{i=1}^{n_k} \bigg(\frac{u^2_{k,i}(1)}{\mu_k}+\frac{u^2_{k,i}(0)}{1-\mu_k}\bigg).
\]

\section{Fixed effects}
\label{asection:fixed_effects} 

\subsection{Large $k$ distribution}

Let
\[
\mwidebar N_{k,m} = \sum_{i=1}^{n_k} 1\{m_{k,i}=m\}R_{k,i}
\]
and
\begin{align}
\label{equation:fixed}
\widehat\tau_k^{\rm{\,fixed}}&= \frac{\displaystyle\sum_{m=1}^{m_k}\sum_{i=1}^{n_k} 1\{m_{k,i}=m\}R_{k,i}Y_{k,i}(W_{k,i}- \mwidebar W_{k,m})}{\displaystyle\sum_{m=1}^{m_k}
\sum_{i=1}^{n_k} 1\{m_{k,i}=m\}R_{k,i}W_{k,i}(W_{k,i}- \mwidebar W_{k,m})},
\end{align}
where 
\[
\mwidebar W_{k,m}=\frac{1}{\mwidebar N_{k,m}\vee 1}\sum_{i=1}^{n_k} 1\{m_{k,i}=m\}R_{k,i}W_{k,i}.
\]
Notice that we need $\liminf_{k\rightarrow\infty} \mu_k(1-\mu_k)-\sigma^2_k=\liminf_{k\rightarrow\infty}  E[A_{k,m}(1-A_{k,m})]>0$  for this estimator to be well-defined in large samples (otherwise, the denominator in the formula for $\widehat\tau_k^{\rm{\,fixed}}$ could be equal to zero). Although it is not strictly necessary, and because it entails little loss of generality and simplifies the exposition, we will assume that the supports of the cluster probabilities, $A_{k,m}$, are bounded away from zero and one (uniformly in $k$ and $m$). In finite samples we assign $\widehat\tau_k^{\rm{\,fixed}}=0$ to the cases when the denominator of $\widehat\tau_k^{\rm{\,fixed}}$ in equation (\ref{equation:fixed}) is equal to zero. Notice that
\[
\sum_{i=1}^{n_k} 1\{m_{k,i}=m\}R_{k,i}(W_{k,i}- \mwidebar W_{k,m}) =0.
\]
Let
\[
\alpha_{k,m} = \frac{1}{n_{k,m}}\sum_{i=1}^{n_k} 1\{m_{k,i}=m\} y_{k,i}(0),\quad
\tau_{k,m} = \frac{1}{n_{k,m}}\sum_{i=1}^{n_k} 1\{m_{k,i}=m\} (y_{k,i}(1)-y_{k,i}(0)),
\]
$e_{k,i}(0)=y_{k,i}(0)-\alpha_{k,m_{k,i}}$, and $e_{k,i}(1)=y_{k,i}(1)-\alpha_{k,m_{k,i}}-\tau_{k,m_{k,i}}$. It follows that
\[
\sum_{i=1}^{n_k} 1\{m_{k,i}=m\}e_{k,i}(1)=\sum_{i=1}^{n_k} 1\{m_{k,i}=m\}e_{k,i}(0)=0.
\]
Now, $Y_{k,i}=e_{k,i}(1)W_{k,i}+e_{k,i}(0)(1-W_{k,i})+\alpha_{k,m_{k,i}}+\tau_{k,m_{k,i}}W_{k,i}$. Then,
\begin{align*}
\widehat\tau_k^{\rm{\,fixed}}&= \frac{\displaystyle\sum_{m=1}^{m_k}\sum_{i=1}^{n_k} 1\{m_{k,i}=m\}R_{k,i}((e_{k,i}(1)+\tau_{k,m})W_{k,i}+e_{k,i}(0)(1-W_{k,i}))(W_{k,i}- \mwidebar W_{k,m})}{\displaystyle\sum_{m=1}^{m_k}
\sum_{i=1}^{n_k} 1\{m_{k,i}=m\}R_{k,i}W_{k,i}(W_{k,i}- \mwidebar W_{k,m})}.
\end{align*}
Let 
\begin{align}
\label{equation:taubar}
\mwidebar\tau_k &= \frac{\displaystyle\sum_{m=1}^{m_k}\tau_{k,m}\sum_{i=1}^{n_k} 1\{m_{k,i}=m\}R_{k,i}W_{k,i}(W_{k,i}- \mwidebar W_{k,m})}{\displaystyle\sum_{m=1}^{m_k}
\sum_{i=1}^{n_k} 1\{m_{k,i}=m\}R_{k,i}W_{k,i}(W_{k,i}- \mwidebar W_{k,m})},
\end{align}
where, as before, we make $\mwidebar\tau_k=0$ if the denominator on the right-hand side of (\ref{equation:taubar}) is equal to zero. Now, $\widehat\tau_k^{\rm{\,fixed}}-\tau_k=(\widehat\tau_k^{\rm{\,fixed}}-\mwidebar\tau_k)+(\mwidebar\tau_k-\tau_k)$, where
\[
\widehat\tau_k^{\rm{\,fixed}}-\mwidebar\tau_k= \frac{\displaystyle\sum_{m=1}^{m_k}\sum_{i=1}^{n_k} 1\{m_{k,i}=m\}R_{k,i}(e_{k,i}(1)W_{k,i}+e_{k,i}(0)(1-W_{k,i}))(W_{k,i}- \mwidebar W_{k,m})}{\displaystyle\sum_{m=1}^{m_k}
\sum_{i=1}^{n_k} 1\{m_{k,i}=m\}R_{k,i}W_{k,i}(W_{k,i}- \mwidebar W_{k,m})}
\]
and
\[
\mwidebar\tau_k -\tau_k= \frac{\displaystyle\sum_{m=1}^{m_k}(\tau_{k,m}-\tau_k)\sum_{i=1}^{n_k} 1\{m_{k,i}=m\}R_{k,i}W_{k,i}(W_{k,i}- \mwidebar W_{k,m})}{\displaystyle\sum_{m=1}^{m_k}
\sum_{i=1}^{n_k} 1\{m_{k,i}=m\}R_{k,i}W_{k,i}(W_{k,i}- \mwidebar W_{k,m})}.
\]
Notice that outcomes enter the term $\widehat\tau_k^{\rm{\,fixed}}-\mwidebar\tau_k$ only through the intra-cluster errors, $e_{k,i}(1)$ and $e_{k,i}(0)$. In contrast, the term $\mwidebar\tau_k-\tau_k$ depends on outcomes only through inter-cluster variability in treatment effects, $\tau_{k,m}-\tau_k$. The numerator in the expression for $\mwidebar\tau_k-\tau_k$ in the last displayed equation does not have mean zero in general, and this will be reflected in a bias term, $B_k$, which we define next.
Let,
\[
D_k =\frac{1}{n_kp_kq_k}\sum_{m=1}^{m_k}\sum_{i=1}^{n_k} 1\{m_{k,i}=m\}R_{k,i}W_{k,i}(W_{k,i}- \mwidebar W_{k,m}),
\]
and
\[
B_k = -\frac{\displaystyle\frac{1}{n_kp_k}E[A_{k,m}(1-A_{k,m})]\sum_{m=1}^{m_k} (\tau_{k,m}-\tau_k)(1-(1-p_k)^{n_{k,m}})}{\displaystyle\frac{1}{n_kp_kq_k}\sum_{m=1}^{m_k}\sum_{i=1}^{n_k} 1\{m_{k,i}=m\}R_{k,i}W_{k,i}(W_{k,i}- \mwidebar W_{k,m})}.
\]
Then, $\sqrt{n_kp_kq_k}(\widehat\tau_k^{\rm{\,fixed}} -\tau_k-B_k)=F_k/D_k$, where
\[
F_k = \sum_{m=1}^{m_k} (\psi_{k,m}-\lowwidebar\psi_{k,m}) + (\varphi_{k,m}-\mwidebar\varphi_{k,m}),
\]
\begin{align*}
\psi_{k,m} &= \frac{1}{\sqrt{n_kp_kq_k}}\sum_{i=1}^{n_k} 1\{m_{k,i}=m\}R_{k,i}(e_{k,i}(1)W_{k,i}+e_{k,i}(0)(1-W_{k,i}))(W_{k,i}- A_{k,m}),\\
\lowwidebar\psi_{k,m} & = \frac{1}{\sqrt{n_kp_kq_k}}\sum_{i=1}^{n_k} 1\{m_{k,i}=m\}R_{k,i}(e_{k,i}(1)W_{k,i}+e_{k,i}(0)(1-W_{k,i}))(\mwidebar W_{k,m}-A_{k,m}),\\
\varphi_{k,m} &= \frac{1}{\sqrt{n_kp_kq_k}}(\tau_{k,m}-\tau_k)\sum_{i=1}^{n_k} 1\{m_{k,i}=m\}\big(R_{k,i}W_{k,i}(W_{k,i}- A_{k,m})-p_kq_kE[A_{k,m}(1-A_{k,m})]\big),\\
\shortintertext{and}
\mwidebar\varphi_{k,m} &= \frac{1}{\sqrt{n_kp_kq_k}}(\tau_{k,m}-\tau_k)\Bigg(\sum_{i=1}^{n_k} 1\{m_{k,i}=m\}R_{k,i}W_{k,i}(\mwidebar W_{k,m}-A_{k,m})\\
&\qquad\qquad\qquad\qquad\qquad\qquad\qquad\qquad\qquad\qquad-q_kE[A_{k,m}(1-A_{k,m})](1-(1-p_k)^{n_{k,m}})\Bigg).
\end{align*}
The terms $\psi_{k,m}$ and $\lowwidebar\psi_{k,m}$ depend on the within-cluster errors $e_{k,i}(1)$ and $e_{k,i}(0)$. The terms $\varphi_{k,m}$ and $\mwidebar\varphi_{k,m}$ depend on the inter-clusters errors $\tau_{k,m}-\tau_k$. $\psi_{k,m}$ and $\varphi_{k,m}$ replace $\mwidebar W_{k,m}$ with $A_{k,m}$, while $\lowwidebar\psi_{k,m}$  and $\mwidebar\varphi_{k,m}$ correct for the difference,  $\mwidebar W_{k,m}-A_{k,m}$.

It can be seen (in intermediate calculations below) that
\[
E\Bigg[\sum_{i=1}^{n_k} 1\{m_{k,i}=m\}R_{k,i}W_{k,i}(W_{k,i}- A_{k,m})\Bigg]=n_{k,m}p_kq_kE[A_{k,m}(1-A_{k,m})]
\]
and
\begin{align*}
E\Bigg[\sum_{i=1}^{n_k} 1\{m_{k,i}=m\}R_{k,i}W_{k,i}(\mwidebar W_{k,m}-A_{k,m})\Bigg]
&=q_kE[A_{k,m}(1-A_{k,m})](1-(1-p_k)^{n_{k,m}}).
\end{align*}
These two expectations are substracted in $\varphi_{k,m}$ and $\mwidebar\varphi_{k,m}$ , so $\varphi_{k,m}$  and $\mwidebar\varphi_{k,m}$ have mean zero. Doing so for  $\varphi_{k,m}$ does not require adjustments elsewhere. Because
\[
\sum_{m=1}^{m_k} (\tau_{k,m}-\tau_k) n_{k,m} =0, 
\]
the $n_{k,m}p_kq_kE[A_{k,m}(1-A_{k,m})]$ terms do not change the sum $F_k$. In contrast, demeaning $\mwidebar\varphi_{k,m}$ creates the bias term $B_k$.
If the size of the clusters $n_{k,m}$ does not vary across clusters, then $B_k$ is equal to zero. More generally, $\sqrt{n_kp_kq_k}D_kB_k = \mathcal O (m_k\sqrt{q_k/(n_kp_k)})$. Therefore, if
\begin{equation}
\label{equation:clobsratio}
\frac{m_kq_k}{p_k (n_k/m_k)}\longrightarrow 0,
\end{equation}
(that is, if the expected number of sampled clusters is small relative to the expected number of sampled observations per sampled cluster) then $\sqrt{n_kp_kq_k}D_kB_k$ converges to zero. As a result, $\sqrt{n_kp_kq_k}B_k$ converges in probability to zero, because, as we will show later, $D_k$ converges in probability to $\mu_k(1-\mu_k)-\sigma_k^2$, which is bounded away from zero. In our large sample analysis, we will assume that the expected number of sampled clusters grows to infinity, $m_kq_k\rightarrow\infty$. Then, equation (\ref{equation:clobsratio}) implies that the expected number of observations per sampled cluster goes to infinity, $p_k (n_k/m_k)\rightarrow\infty$. Notice also that $n_kp_kq_k=(n_kp_k/m_k)(m_kq_k)\rightarrow\infty$.

We summarize now the assumptions we made thus far. We first assumed that the supports of the cluster probabilities, $A_{k,m}$, are bounded away from zero and one (uniformly in $k$ and $m$), and that potential outcomes are bounded. Moreover, we assumed $m_kq_k\rightarrow \infty$ and $(m_kq_k)/((p_k n_k)/m_k)\rightarrow 0$.
These imply $(p_k n_k)/m_k\rightarrow\infty$ and $n_kp_kq_k\rightarrow\infty$. 
We will add the assumption that the ratio between maximum and minimum cluster size is bounded, $\limsup_{k\rightarrow\infty} \max_m n_{k,m}/\min_m n_{k,m}<\infty$. This assumption implies $p_k\min_m n_{k,m}\rightarrow\infty$ and $(m_kq_k)/(p_k \min_m n_{k,m})\rightarrow 0$.

We will now study the behavior of $D_k$. Notice that
\begin{align*}
E\Bigg[\sum_{m=1}^{m_k}\sum_{i=1}^{n_k} &1\{m_{k,i}=m\}R_{k,i}W_{k,i}(W_{k,i}- \mwidebar W_{k,m})\Bigg]\\
&=E\Bigg[\sum_{m=1}^{m_k}\sum_{i=1}^{n_k} 1\{m_{k,i}=m\}R_{k,i}W_{k,i}(W_{k,i}- A_{k,m})\Bigg]\\&-E\Bigg[\sum_{m=1}^{m_k}\sum_{i=1}^{n_k} 1\{m_{k,i}=m\}R_{k,i}W_{k,i}(\mwidebar W_{k,m}- A_{k,m})\Bigg]\\
&=n_kp_kq_kE[A_{k,m}(1-A_{k,m})] - q_kE[A_{k,m}(1-A_{k,m})]\sum_{m=1}^{m_k}(1-(1-p_k)^{n_{k,m}}).
\end{align*}
In addition, 
\begin{align*}
\frac{1}{(n_kp_kq_k)^2}\sum_{m=1}^{m_k}E\Bigg[\Bigg(&\sum_{i=1}^{n_k} 1\{m_{k,i}=m\}R_{k,i}W_{k,i}(W_{k,i}- \mwidebar W_{k,m})\Bigg)^2\Bigg]\\
&\leq c\,\frac{n_kp_kq_k+n_kp_k^2q_k\max_m n_{k,m}}{(n_kp_kq_k)^2}\\
&=c\,\Bigg(\frac{1}{n_kp_kq_k}+\frac{\max_m n_{k,m}}{\min_m n_{k,m}}\frac{1}{m_kq_k}\Bigg)\longrightarrow 0.
\end{align*}
The weak law of large numbers for arrays implies
\[
D_k - E[A_{k,m}(1-A_{k,m})] + \frac{1}{n_kp_k}E[A_{k,m}(1-A_{k,m})]\sum_{m=1}^{m_k} (1-(1-p_k)^{n_{k,m}}) \stackrel{p}{\longrightarrow} 0. 
\]
Because $m_k/(n_kp_k)\rightarrow 0$ and $E[A_{k,m}(1-A_{k,m})]=\mu_k(1-\mu_k)-\sigma_k^2$, we obtain
\[
D_k - (\mu_k(1-\mu_k)-\sigma_k^2)\stackrel{p}{\longrightarrow} 0. 
\]

We now turn our attention to $F_k$. 
We will first calculate the variance of $\psi_{k,m}$. Let $Q_{k,m}$ be a binary variable that takes value one if cluster $m$ in population $k$ is sampled, and zero otherwise. Notice that 
\[
E[R_{k,i}W_{k,i}(W_{k,i}-A_{k,m})|A_{k,m}, Q_{k,m}=1,m_{k,i}=m] = p_kA_{k,m}(1-A_{k,m}),
\]
and
\[
E[R_{k,i}(1-W_{k,i})(W_{k,i}-A_{k,m})|A_{k,m}, Q_{k,m}=1,m_{k,i}=m] = -p_kA_{k,m}(1-A_{k,m}).
\]
Consider now 
\begin{align*}
\psi_{k,m,1} &= \frac{1}{\sqrt{n_kp_kq_k}}\sum_{i=1}^{n_k} 1\{m_{k,i}=m\}R_{k,i}W_{k,i}(W_{k,i}-A_{k,m})e_{k,i}(1)\\
&=\frac{Q_{k,m}}{\sqrt{n_kp_kq_k}}\sum_{i=1}^{n_k} 1\{m_{k,i}=m\}\Big(R_{k,i}W_{k,i}(W_{k,i}-A_{k,m})- p_kA_{k,m}(1-A_{k,m})\Big)e_{k,i}(1),
\end{align*}
and
\begin{align*}
\psi_{k,m,0}&=\frac{1}{\sqrt{n_kp_kq_k}}\sum_{i=1}^{n_k} 1\{m_{k,i}=m\}R_{k,i}(1-W_{k,i})(W_{k,i}-A_{k,m})e_{k,i}(0)\\
&=\frac{Q_{k,m}}{\sqrt{n_kp_kq_k}}\sum_{i=1}^{n_k} 1\{m_{k,i}=m\}\Big(R_{k,i}(1-W_{k,i})(W_{k,i}-A_{k,m})+p_kA_{k,m}(1-A_{k,m})\Big)e_{k,i}(0).
\end{align*}
It holds that $\psi_{k,m}=\psi_{k,m,1}+\psi_{k,m,0}$ and $E[\psi_{k,m}]=0$. Now, notice that
\begin{align*}
E[\psi_{k,m,1}^2]&=\frac{1}{n_k}E[A_{k,m}(1-A_{k,m})^2-p_kA_{k,m}^2(1-A_{k,m})^2]\sum_{i=1}^{n_k} 1\{m_{k,i}=m\}e_{k,i}^2(1),\\
E[\psi_{k,m,0}^2]&=\frac{1}{n_k}E[A^2_{k,m}(1-A_{k,m})-p_kA_{k,m}^2(1-A_{k,m})^2]\sum_{i=1}^{n_k} 1\{m_{k,i}=m\}e_{k,i}^2(0),\\
\shortintertext{and}
E[\psi_{k,m,1}\psi_{k,m,0}]&=\frac{1}{n_k} p_k E[A_{k,m}^2(1-A_{k,m})^2]\sum_{i=1}^{n_k} 1\{m_{k,i}=m\}e_{k,i}(1)e_{k,i}(0).
\end{align*}
Therefore,
\begin{align*}
E[(\psi_{k,m,1}+\psi_{k,m,0})^2]&=\frac{1}{n_k}E[A_{k,m}(1-A_{k,m})^2]\sum_{i=1}^{n_k} 1\{m_{k,i}=m\}e_{k,i}^2(1)\\
&+\frac{1}{n_k}E[A^2_{k,m}(1-A_{k,m})]\sum_{i=1}^{n_k} 1\{m_{k,i}=m\}e_{k,i}^2(0)\\
&- \frac{1}{n_k}p_kE[A^2_{k,m}(1-A_{k,m})^2]\sum_{i=1}^{n_k} 1\{m_{k,i}=m\}(e_{k,i}(1)-e_{k,i}(0))^2,
\end{align*}
and
\begin{align}
\label{equation:psi2}
\sum_{m=1}^{m_k}E[(\psi_{k,m,1}+\psi_{k,m,0})^2]&=E[A_{k,m}(1-A_{k,m})^2]\frac{1}{n_k}\sum_{i=1}^{n_k}  e_{k,i}^2(1)\nonumber\\
&+E[A^2_{k,m}(1-A_{k,m})]\frac{1}{n_k}\sum_{i=1}^{n_k} e_{k,i}^2(0)\nonumber\\
&- p_k E[A^2_{k,m}(1-A_{k,m})^2]\frac{1}{n_k}\sum_{i=1}^{n_k}(e_{k,i}(1)-e_{k,i}(0))^2.
\end{align} 
We will next show that the terms $\lowwidebar\psi_{k,m}$ do not matter for the asymptotic distribution of $\sqrt{n_kp_kq_k}(\widehat\tau_k-\tau_k)$. Notice that, because the cluster sum of $e_{k,i}(1)$ is equal to zero, we obtain $E[\mwidebar\psi_{k,m}]=0$ and, therefore,
\[
\sum_{m=1}^{m_k} E\Big[\mwidebar\psi_{k,m}\Big]=0.
\]
Moreover
\[
2\sum_{i=1}^{n_k-1}\sum_{j=i+1}^{n_k}1\{m_{k,i}=m_{k,j}=m\}e_{k,i}(1)e_{k,j}(1) = -\sum_{i=1}^{n_k}1\{m_{k,i}=m\}e_{k,i}^2(1) \leq 0.
\]
In addition, $E[R_{k,i}W_{k,i}(\mwidebar W_{k,m}-A_{k,m})^2|m_{k,i}=m]\leq q_kE[A_{k,m}(1-A_{k,m})]/n_{k,m}$ (see intermediate calculations). Therefore, 
\begin{align*}
E\Bigg[\Bigg(&\sum_{i=1}^{n_k}1\{m_{k,i}=m\}R_{k,i}W_{k,i}(\mwidebar W_{k,m}-A_{k,m})e_{k,i}(1)\Bigg)^2\Bigg]\\
&=\sum_{i=1}^{n_k}1\{m_{k,i}=m\}E\Big[R_{k,i}W_{k,i}(\mwidebar W_{k,m}-A_{k,m})^2|m_{k,i}=m\Big]e_{k,i}^2(1)\\
&+ 2\sum_{i=1}^{n_k-1}\sum_{j=i+1}^{n_k}E\Big[1\{m_{k,i}=m_{k,j}=m\}R_{k,i}R_{k,j}W_{k,i}W_{k,j}(\mwidebar W_{k,m}-A_{k,m})^2\Big]e_{k,i}(1)e_{k,j}(1)\\
&\leq q_kE[A_{k,m}(1-A_{k,m})]\frac{1}{n_{k,m}}\sum_{i=1}^{n_k}1\{m_{k,i}=m\}e_{k,i}^2(1).
\end{align*}
Now, because errors are bounded, we obtain
\begin{equation}
\label{equation:sqsume1}
\sum_{m=1}^{m_k}E\Bigg[\Bigg(\frac{1}{\sqrt{n_kp_kq_k}}\sum_{i=1}^{n_k}1\{m_{k,i}=m\}R_{k,i}W_{k,i}(\mwidebar W_{k,m}-A_{k,m})e_{k,i}(1)\Bigg)^2\Bigg]\leq c\, \frac{m_k}{n_kp_k}.  
\end{equation}
Because $m_k/(n_kp_k)\rightarrow 0$, the weak law of large numbers for arrays, implies,
\[
\frac{1}{\sqrt{n_kp_kq_k}}\sum_{m=1}^{m_k}\sum_{i=1}^{n_k} 1\{m_{k,i}=m\}R_{k,i}W_{k,i}(\mwidebar W_{k,m}-A_{k,m})e_{k,i}(1)\stackrel{p}{\longrightarrow} 0. 
\]
with the analogous result involving the errors $e_{k,i}(0)$. If follows that
\[
\sum_{m=1}^{m_k} \lowwidebar\psi_{k,m}\stackrel{p}{\longrightarrow} 0. 
\]
Consider now $\varphi_{k,m}$. Notice that
\begin{align*}
E\Big[\Big(R_{k,i}W_{k,i}&(W_{k,i}-A_{k,m})-p_kq_kE[A_{k,m}(1-A_{k,m})])\Big)^2\Big]\\
&= p_kq_kE[A_{k,m}(1-A_{k,m})^2]-p_k^2q_k^2\big(E[A_{k,m}(1-A_{k,m})]\big)^2,
\end{align*}
and
\begin{align*}
E\Big[\Big(R_{k,i}W_{k,i}&(W_{k,i}-A_{k,m})-p_kq_kE[A_{k,m}(1-A_{k,m})])\Big)\\
&\times\Big(R_{k,j}W_{k,j}(W_{k,j}-A_{k,m})-p_kq_kE[A_{k,m}(1-A_{k,m})])\Big)\big|m_{k,i}=m_{k,j}=m\Big]\\
&\qquad =p_k^2q_kE[A^2_{k,m}(1-A_{k,m})^2]-p_k^2q_k^2\big(E[A_{k,m}(1-A_{k,m})]\big)^2.
\end{align*}
Therefore, 
\begin{align*}
E[\varphi_{k,m}^2]&= \Big(E[A_{k,m}(1-A_{k,m})^2]-p_kq_k(E[A_{k,m}(1-A_{k,m})])^2\Big)\frac{n_{k,m}}{n_k}(\tau_{k,m}-\tau_k)^2\\
&+\Big(p_kE[A^2_{k,m}(1-A_{k,m})^2]-p_kq_k(E[A_{k,m}(1-A_{k,m})])^2\Big)\frac{n_{k,m}(n_{k,m}-1)}{n_k}(\tau_{k,m}-\tau_k)^2,
\end{align*}
and
\begin{align*}
\sum_{m=1}^{m_k}E[\varphi_{k,m}^2]&= \Big(E[A_{k,m}(1-A_{k,m})^2]-p_kq_k(E[A_{k,m}(1-A_{k,m})])^2\Big)\sum_{m=1}^{m_k}\frac{n_{k,m}}{n_k}(\tau_{k,m}-\tau_k)^2\\
&+\Big(p_kE[A^2_{k,m}(1-A_{k,m})^2]-p_kq_k(E[A_{k,m}(1-A_{k,m})])^2\Big)\sum_{m=1}^{m_k}\frac{n_{k,m}(n_{k,m}-1)}{n_k}(\tau_{k,m}-\tau_k)^2.
\end{align*}

Next, we calculate the variance of $\mwidebar\varphi_{k,m}$. Using results on the moments of a Binomial distribution, we obtain, for $n\geq 1$,
\begin{align*}
E\Bigg[\Bigg(\sum_{i=1}^{n_k}&1\{m_{k,i}=m\}R_{k,i}W_{k,i}(\mwidebar W_{k,m}-A_{k,m})\Bigg)^2\Big| \begin{array}{l}Q_{k,m}=1,\\\mwidebar N_{k,m}=n\end{array}\Bigg]\\
&=\frac{1}{n^2}E\Bigg[\Bigg(\sum_{i=1}^{n_k} 1\{m_{k,i}=m\}R_{k,i}W_{k,i}\Big(\sum_{i=1}^{n_k} 1\{m_{k,i}=m\}R_{k,i}W_{k,i}-nA_{k,m}\Big)\Bigg)^2\Big| \begin{array}{l}Q_{k,m}=1,\\\mwidebar N_{k,m}=n\end{array}\Bigg]\\
&=n E[A_{k,m}^3(1-A_{k,m})]+E[A^2_{k,m}(1-A_{k,m})(5-7A_{k,m})]\\&+
\frac{1}{n}E[A_{k,m}(1-A_{k,m})(6A^2_{k,m} - 6 A_{k,m} + 1)].
\end{align*}
Therefore,
\begin{align*}
E\Bigg[\Bigg(\sum_{i=1}^{n_k}&1\{m_{k,i}=m\}R_{k,i}W_{k,i}(\mwidebar W_{k,m}-A_{k,m})\Bigg)^2\Bigg]\\
&=n_{k,m}p_kq_k E[A_{k,m}^3(1-A_{k,m})]+q_kE[A^2_{k,m}(1-A_{k,m})(5-7A_{k,m})](1-(1-p_k)^{n_{k,m}})\\&+
q_kE[A_{k,m}(1-A_{k,m})(6A^2_{k,m} - 6 A_{k,m} + 1)]r_{k,m}, 
\end{align*}
where 
\[
r_{k,m} = \sum_{n=1}^{n_{k,m}}\frac{1}{n} 
\Pr(\mwidebar N_{k,m}=n|Q_{k,m}=1)\leq \sum_{n=1}^{n_{k,m}}
\Pr(\mwidebar N_{k,m}=n|Q_{k,m}=1)\leq 1.
\]
It follows that,
\begin{align*}
E[\mwidebar\varphi^2_{k,m}]&=(\tau_{k,m}-\tau_k)^2\Big(\frac{n_{k,m}}{n_k}E[A_{k,m}^3(1-A_{k,m})]+\frac{1}{n_kp_k}E[A^2_{k,m}(1-A_{k,m})(5-7A_{k,m})](1-(1-p_k)^{n_{k,m}})\\&+
\frac{1}{n_kp_k}E[A_{k,m}(1-A_{k,m})(6A^2_{k,m} - 6 A_{k,m} + 1)]r_{k,m}
\\&-\frac{q_k}{n_kp_k}(E[A_{k,m}(1-A_{k,m})])^2(1-(1-p_k)^{n_{k,m}})^2\Big).
\end{align*}
Therefore,
\[
\sum_{m=1}^{m_k}E[\mwidebar\varphi^2_{k,m}] = \sum_{m=1}^{m_k}
(\tau_{k,m}-\tau_k)^2\Big(\frac{n_{k,m}}{n_k}\Big)E[A_{k,m}^3(1-A_{k,m})] + \scaleto{\mathcal{O}}{5pt}(1).
\]

We will now study the covariance between $\varphi_{k,m}$ and $\lowwidebar\varphi_{k,m}$.  Using results on the moments of a Binomial distribution, we obtain, for $n\geq 1$,
\begin{align*}
E\Bigg[\Bigg(&\sum_{i=1}^{n_k}1\{m_{k,i}=m\}R_{k,i}W_{k,i}(W_{k,i}-A_{k,m})\Bigg)\Bigg(\sum_{i=1}^{n_k}1\{m_{k,i}=m\}R_{k,i}W_{k,i}(\mwidebar W_{k,m}-A_{k,m})\Bigg)\Big| \begin{array}{l}Q_{k,m}=1,\\\mwidebar N_{k,m}=n\end{array}\Bigg]\\
&=E\Bigg[\frac{1-A_{k,m}}{n}\Bigg(\sum_{i=1}^{n_k}1\{m_{k,i}=m\}R_{k,i}W_{k,i}\Bigg)^2\Bigg(\sum_{i=1}^{n_k}1\{m_{k,i}=m\}R_{k,i}W_{k,i}-nA_{k,m})\Bigg)\Big| \begin{array}{l}Q_{k,m}=1,\\\mwidebar N_{k,m}=n\end{array}\Bigg]\\
&=2nE[A^2_{k,m}(1-A_{k,m})^2]+E[A_{k,m}(1-A_{k,m})^2(1-2A_{k,m})].
\end{align*}
Therefore,
\begin{align*}
E\Bigg[\Bigg(&\sum_{i=1}^{n_k}1\{m_{k,i}=m\}R_{k,i}W_{k,i}(W_{k,i}-A_{k,m})\Bigg)\Bigg(\sum_{i=1}^{n_k}1\{m_{k,i}=m\}R_{k,i}W_{k,i}(\mwidebar W_{k,m}-A_{k,m})\Bigg)\Bigg]\\
&=2n_{k,m}p_kq_kE[A^2_{k,m}(1-A_{k,m})^2]+q_kE[A_{k,m}(1-A_{k,m})^2(1-2A_{k,m})] \Pr(\mwidebar N_{k,m}\geq 1|Q_{k,m}=1).
\end{align*}
In addition,
\begin{align*}
E\Bigg[&\sum_{i=1}^{n_k}1\{m_{k,i}=m\}R_{k,i}W_{k,i}(W_{k,i}-A_{k,m})\Bigg]E\Bigg[\sum_{i=1}^{n_k}1\{m_{k,i}=m\}R_{k,i}W_{k,i}(\mwidebar W_{k,m}-A_{k,m})\Bigg]\\
&=n_{k,m}p_kq^2_k(E[A_{k,m}(1-A_{k,m})])^2
\Pr(\mwidebar N_{k,m}\geq 1|Q_{k,m}=1).
\end{align*}
As a result,
\begin{align*}
E[\varphi_{k,m}\lowwidebar\varphi_{k,m}]&=\Big(2E[A^2_{k,m}(1-A_{k,m})^2]-q_k(E[A_{k,m}(1-A_{k,m})])^2\Big)(\tau_{k,m}-\tau_k)^2\left(\frac{n_{k,m}}{n_k}\right)+\mathcal O\left(\frac{1}{n_kp_k}\right)\\
&+\mathcal O\left(\frac{q_k}{p_k\min_m n_{k,m}}\big(p_k\min_m n_{k,m}(1-p_k)^{\min_m n_{k,m}}\big)\right). 
\end{align*}
Notice that $m_k/(n_kp_k)\rightarrow 0$. In addition, $m_kq_k/(p_k\min_m n_{k,m})\rightarrow 0$ and
\begin{align*}
p_k\min_m n_{k,m} (1-p_k)^{\min_m n_{k,m}}&= p_k\min_m n_{k,m} \left(1-\frac{p_k\min_m n_{k,m}}{\min_m n_{k,m}}\right)^{\min_m n_{k,m}}\\
&< p_k\min_m n_{k,m} e^{-p_k\min_m n_{k,m}}\longrightarrow 0.
\end{align*}
Therefore,
\[
\sum_{m=1}^{m_k} E[\varphi_{k,m}\lowwidebar\varphi_{k,m}]= \Big(2E[A^2_{k,m}(1-A_{k,m})^2]-q_k(E[A_{k,m}(1-A_{k,m})])^2\Big)\sum_{m=1}^{m_k}(\tau_{k,m}-\tau_k)^2\left(\frac{n_{k,m}}{n_k}\right)+\ \scaleto{\mathcal{O}}{5pt}(1).
\]
Next, we will study the remaining covariances between $\psi_{k,m}$, $\varphi_{k,m}$, $\lowwidebar\psi_{k,m}$, and $\mwidebar\varphi_{k,m}$. Because the intra-cluster errors, $e_{k,i}(1)$ and $e_{k,i}(0)$ sum to zero, it can be easily seen that $E[\psi_{k,m}\varphi_{k,m}]=E[\psi_{k,m}\mwidebar\varphi_{k,m}]=0$. It can also be seen that the inter-clusters sums of covariances between $\lowwidebar\psi_{k,m}$ and any of the other terms go to zero. To prove this for the covariance with $\psi_{k,m}$, 
we have
\begin{align*}
\left(\sum_{m=1}^{m_k} E[|\psi_{k,m}\lowwidebar\psi_{k,m}|]\right)^2&\leq \left(\sum_{m=1}^{m_k} (E[\psi^2_{k,m}]E[\lowwidebar\psi^2_{k,m}])^{1/2}\right)^2\\
&\leq \sum_{m=1}^{m_k} E[\psi^2_{k,m}]\sum_{m=1}^{m_k} E[\lowwidebar\psi^2_{k,m}]\\
&=\mathcal O(1)\scaleto{\mathcal{O}}{5pt}(1)=\scaleto{\mathcal{O}}{5pt}(1).
\end{align*}
The same argument and result applies to $E[\lowwidebar\psi_{k,m}\varphi_{k,m}]$ and $E[\lowwidebar\psi_{k,m}\mwidebar\varphi_{k,m}]$.
Putting all the pieces together, we obtain
\begin{align*}
n_kp_kq_kE[D_k^2(\widehat\tau_k^{\rm{\,fixed}}-\tau_k)^2] =  f_k +\scaleto{\mathcal{O}}{5pt}(1),
\end{align*}
where
\begin{align*}
f_k&=E[A_{k,m}(1-A_{k,m})^2]\frac{1}{n_k}\sum_{i=1}^{n_k} e_{k,i}^2(1)
+E[A^2_{k,m}(1-A_{k,m})]\frac{1}{n_k}\sum_{i=1}^{n_k} e_{k,i}^2(0)\\
&- p_k E[A^2_{k,m}(1-A_{k,m})^2]\frac{1}{n_k}\sum_{i=1}^{n_k}(e_{k,i}(1)-e_{k,i}(0))^2\\
&+\Big(E[A_{k,m}(1-A_{k,m})^2]-p_kq_k(E[A_{k,m}(1-A_{k,m})])^2\Big)\sum_{m=1}^{m_k}\frac{n_{k,m}}{n_k}(\tau_{k,m}-\tau_k)^2\\
&+\Big(p_kE[A^2_{k,m}(1-A_{k,m})^2]-p_kq_k(E[A_{k,m}(1-A_{k,m})])^2\Big)\sum_{m=1}^{m_k}\frac{n_{k,m}(n_{k,m}-1)}{n_k}(\tau_{k,m}-\tau_k)^2\\
&+E[A_{k,m}^3(1-A_{k,m})]\sum_{m=1}^{m_k}
(\tau_{k,m}-\tau_k)^2\Big(\frac{n_{k,m}}{n_k}\Big)\\
&-2\Big(2E[A^2_{k,m}(1-A_{k,m})^2]-q_k(E[A_{k,m}(1-A_{k,m})])^2\Big)\sum_{m=1}^{m_k}(\tau_{k,m}-\tau_k)^2\left(\frac{n_{k,m}}{n_k}\right).
\end{align*}
Collecting terms with identical factors, we obtain
\begin{align*}
f_k&=E[A_{k,m}(1-A_{k,m})^2]\frac{1}{n_k}\sum_{i=1}^{n_k} e_{k,i}^2(1)
+E[A^2_{k,m}(1-A_{k,m})]\frac{1}{n_k}\sum_{i=1}^{n_k} e_{k,i}^2(0)\\
&- p_k E[A^2_{k,m}(1-A_{k,m})^2]\frac{1}{n_k}\sum_{i=1}^{n_k}(e_{k,i}(1)-e_{k,i}(0))^2\\
&+\Big(E[A_{k,m}(1-A_{k,m})^2]-(4+p_k)E[A^2_{k,m}(1-A_{k,m})^2]\\
&\qquad\qquad +E[A_{k,m}^3(1-A_{k,m})] + 2q_k(E[A_{k,m}(1-A_{k,m})])^2 \Big)\sum_{m=1}^{m_k}\frac{n_{k,m}}{n_k}(\tau_{k,m}-\tau_k)^2\\
&+\Big(p_kE[A^2_{k,m}(1-A_{k,m})^2]-p_kq_k(E[A_{k,m}(1-A_{k,m})])^2\Big)\sum_{m=1}^{m_k}\frac{n^2_{k,m}}{n_k}(\tau_{k,m}-\tau_k)^2.
\end{align*}
The first three terms in the expression above depend on intra-cluster heterogeneity in potential outcomes and treatment effects. The last two terms depend on inter-cluster variation in average treatment effects. 

A more compact expression for $f_k$ is
\begin{align}
\label{equation:fformula}
f_k&=E[A_{k,m}(1-A_{k,m})^2]\frac{1}{n_k}\sum_{i=1}^{n_k} e_{k,i}^2(1)
+E[A^2_{k,m}(1-A_{k,m})]\frac{1}{n_k}\sum_{i=1}^{n_k} e_{k,i}^2(0)\nonumber\\
&- p_k E[A^2_{k,m}(1-A_{k,m})^2]\frac{1}{n_k}\sum_{i=1}^{n_k}(e_{k,i}(1)-e_{k,i}(0))^2\nonumber\\
&+\Big(E[A_{k,m}(1-A_{k,m})]-(5+p_k)E[A^2_{k,m}(1-A_{k,m})^2]\nonumber\\
&\qquad\qquad + 2q_k(E[A_{k,m}(1-A_{k,m})])^2 \Big)\sum_{m=1}^{m_k}\frac{n_{k,m}}{n_k}(\tau_{k,m}-\tau_k)^2\nonumber\\
&+\Big(p_kE[A^2_{k,m}(1-A_{k,m})^2]-p_kq_k(E[A_{k,m}(1-A_{k,m})])^2\Big)\sum_{m=1}^{m_k}\frac{n^2_{k,m}}{n_k}(\tau_{k,m}-\tau_k)^2.
\end{align}
Notice that the first four terms in (\ref{equation:fformula}) are bounded, and that
\[
E[A^2_{k,m}(1-A_{k,m})^2]-q_k(E[A_{k,m}(1-A_{k,m})])^2=\mbox{var}(A_{k,m}(1-A_{k,m}))+(1-q_k)(E[A_{k,m}(1-A_{k,m})])^2.
\] 
Assume that 
\begin{equation}
\label{equation:chte}
\liminf_{k\rightarrow\infty} \sum_{m=1}^{m_k}\frac{n_{k,m}}{n_k}(\tau_{k,m}-\tau_k)^2>0,
\end{equation}
and 
\begin{equation}
\label{equation:varAq}
\liminf_{k\rightarrow\infty} \mbox{var}(A_{k,m}(1-A_{k,m}))\vee (1-q_k)>0.
\end{equation}
The last term in equation (\ref{equation:fformula}) is greater than
\[
p_k\min_m n_{k,m}\big(E[A^2_{k,m}(1-A_{k,m})^2]-(E[A_{k,m}(1-A_{k,m})])^2\big)\sum_{m=1}^{m_k}\frac{n_{k,m}}{n_k}(\tau_{k,m}-\tau_k)^2,
\]
which converges to infinity because $p_k\min_m n_{k,m}\rightarrow \infty$. That is, the last term dominates the variance in large samples provided that \eqref{equation:chte} and \eqref{equation:varAq} hold.

We will now derive the large sample distribution of $\widehat\tau_k^{\rm{\,fixed}}$. To show that Lyapunov's condition holds for $F_k$, notice that 
\begin{align*}
|(&\psi_{k,m}-\lowwidebar\psi_{k,m}) + (\varphi_{k,m}-\mwidebar\varphi_{k,m})|^3\\
&=\frac{1}{(n_kp_kq_k)^{3/2}}\Bigg|\sum_{i=1}^{n_k} 1\{m_{k,i}=m\}R_{k,i}((e_{k,i}(1)+\tau_{k,m}-\tau_k)W_{k,i}+e_{k,i}(0)(1-W_{k,i}))(W_{k,i}- \mwidebar W_{k,m})\\
&\qquad\qquad\qquad\qquad\qquad - (\tau_{k,m}-\tau)q_kE[A_{k,m}(1-A_{k,m})](1-(1-p_k)^{n_{k,m}})\Bigg|^3,
\end{align*}
where the last term inside the absolute value comes from the bias correction. Notice that, 
\begin{align*}
\Bigg|\sum_{i=1}^{n_k} 1\{m_{k,i}=m\}&R_{k,i}(e_{k,i}(1)+\tau_{k,m}-\tau_k)W_{k,i}(W_{k,i}- \mwidebar W_{k,m})\Bigg|^3\\
&=\Bigg|(1- \mwidebar W_{k,m})\sum_{i=1}^{n_k} 1\{m_{k,i}=m\}R_{k,i}(e_{k,i}(1)+\tau_{k,m}-\tau_k)W_{k,i}\Bigg|^3\\
&\leq c\,\Bigg|\sum_{i=1}^{n_k} 1\{m_{k,i}=m\}R_{k,i}W_{k,i}\Bigg|^3\\
&\leq c\, \mwidebar N_{k,m}^3. 
\end{align*}
From the formula of the third moment of a binomial random variable, we obtain
\begin{align*}
E[\mwidebar N_{k,m}^3]&=q_k E[\mwidebar N_{k,m}^3|Q_{k,m}=1]\\
&= n_{k,m}^3p_k^3q_k + \scaleto{\mathcal{O}}{5pt}(n_{k,m}^3p_k^3q_k),
\end{align*}
as $p_kn_{k,m}\rightarrow \infty$.
 Now,
\begin{align*}
\frac{1}{f_k^{3/2}}\sum_{k=1}^{m_k}E\Bigg[\Bigg|&\frac{1}{\sqrt{n_kp_kq_k}}\sum_{i=1}^{n_k} 1\{m_{k,i}=m\}R_{k,i}(e_{k,i}(1)+\tau_{k,m}-\tau_k)W_{k,i}(W_{k,i}- \mwidebar W_{k,m})\Bigg|^3\Bigg]\\
&\leq c\,\frac{n_k\max_n n_{k,m}^2 p^3_kq_k}{(n_kp_kq_k)^{3/2}(p_k\min_m n_{k,m})^{3/2}}
=c\left(\frac{\max_m n_{k,m}}{\min_m n_{k,m}}\right)^2\frac{1}{(m_kq_k)^{1/2}}\longrightarrow 0. 
\end{align*} 

Similar calculations deliver the analogous result for the term involving $e_{k,i}(0)$, and proving the result for the bias term is straightforward. Therefore, we obtain
\[
\frac{1}{f_k^{3/2}}\sum_{m=1}^{m_k}|(\psi_{k,m}-\lowwidebar\psi_{k,m}) + (\varphi_{k,m}-\mwidebar\varphi_{k,m})|^3\longrightarrow 0.
\]
By the Central Limit Theorem for arrays, this implies
\[
\sqrt{n_kp_kq_k}F_k/f_k^{1/2}\stackrel{d}{\longrightarrow} N(0,1).
\]
Let $\tilde v_k = f_k/(\mu_k(1-\mu_k)-\sigma^2_k)^2$. Then,
\[
\sqrt{n_kp_kq_k}(\widehat\tau_k^{\rm{\,fixed}} -\tau_k)/\tilde v_k^{1/2}\stackrel{d}{\longrightarrow} N(0,1).
\]
As a result, 
\[
\sqrt{N_k}(\widehat\tau_k^{\rm{\,fixed}} -\tau_k)/\tilde v_k^{1/2}\stackrel{d}{\longrightarrow} N(0,1).
\]

\subsection{Estimation of the variance}
Let
\begin{align*}
N_{k,m,0}&=\sum_{i=1}^{n_k} 1\{m_{k,i}=m\}R_{k,i} (1-W_{k,i})\\
\shortintertext{and}
N_{k,m,1}&=\sum_{i=1}^{n_k} 1\{m_{k,i}=m\}R_{k,i} W_{k,i}.
\end{align*}
Let
\[
\mwidebar Y_{k,m}=\frac{1}{\mwidebar N_{k,m}\vee 1}\sum_{i=1}^{n_k} 1\{m_{k,i}=m\}R_{k,i} Y_{k,i}.
\]
Then,
\[
\mwidebar Y_{k,m}=\widehat\alpha_{k,m}+\widehat\tau_{k,m}\mwidebar W_{k,m},
\]
where
\begin{gather*}
\widehat\alpha_{k,m}=\frac{1}{N_{k,m,0}\vee 1}\sum_{i=1}^{n_k} 1\{m_{k,i}=m\}R_{k,i} (1-W_{k,i})Y_{k,i},\\
\widehat\tau_{k,m}=\frac{1}{N_{k,m,1}\vee 1}\sum_{i=1}^{n_k} 1\{m_{k,i}=m\}R_{k,i} W_{k,i}Y_{k,i}-\frac{1}{N_{k,m,0}\vee 1}\sum_{i=1}^{n_k} 1\{m_{k,i}=m\}R_{k,i} (1-W_{k,i})Y_{k,i},\\
\shortintertext{and, as before,}
\mwidebar W_{k,m}=\frac{1}{\mwidebar N_{k,m}\vee 1}\sum_{i=1}^{n_k} 1\{m_{k,i}=m\}R_{k,i} W_{k,i}.
\end{gather*}
Let $\widetilde U_{k,i}=\widetilde Y_{k,i} - \widehat\tau_k^{\,{\rm fixed}} \widetilde W_{k,i}$, where $\widetilde Y_{k,i}=Y_{k,i}-\mwidebar Y_{k,m_{k,i}}$, $\widetilde W_{k,i}=(W_{k,i}-\mwidebar W_{k,m_{k,i}})$, and $\widehat\tau_k^{\,\rm fixed}$ is the within estimator of $\tau_k$. 
Let $\widetilde\Sigma_k=\sum_{m=1}^{m_k} \widetilde\Sigma_{k,m}$, where
\begin{align*}
\widetilde\Sigma_{k,m}&= \left(\sum_{i=1}^{n_k}1\{m_{k,i}=m\}R_{k,i}\widetilde W_{k,i}\widetilde U_{k,i}\right)^2.
\end{align*}
Also, let 
\[
\widetilde Q_k = \sum_{i=1}^{n_k} R_{k,i}\widetilde W_{k,i}^2.
\]
Then, the cluster estimator of the variance of 
$\sqrt{N_k}(\widehat\tau^{\,{\rm fixed}}_k-\tau_k)$ is 
\begin{align*}
\widetilde V_k^{\rm{cluster}} = N_k\widetilde Q_k^{-1}\widetilde\Sigma_k \widetilde Q_k^{-1}.
\end{align*}
We know already that
\[
\frac{1}{n_kp_kq_k} \widetilde Q_k-(\mu_k(1-\mu_k)-\sigma_k^2)\stackrel{p}{\longrightarrow} 0,
\]
with $\mu_k(1-\mu_k)-\sigma_k^2$ bounded away from zero. To establish convergence of $\widetilde\Sigma_k/(n_kp_kq_kf_k)$, first notice that, for $m_{k,i}=m$, we have
\begin{align*}
\widetilde U_{k,i} & = Y_{k,i} - (\widehat\alpha_{k,m}+\widehat\tau_{k,m}\mwidebar W_{k,m})-\widehat\tau_k^{\,{\rm fixed}}(W_{k,i}-\mwidebar W_{k,m})\\
&=y_{k,i}(1)W_{k,i}+y_{k,i}(0)(1-W_{k,i})-(\alpha_{k,m}+\tau_{k,m}\mwidebar W_{k,m})-\widehat\tau_k^{\,{\rm fixed}}(W_{k,i}-\mwidebar W_{k,m})\\
&-(\widehat\alpha_{k,m}-\alpha_{k,m}) - (\widehat\tau_{k,m}-\tau_{k,m})\mwidebar W_{k,m}\\
&=e_{k,i}(1)W_{k,i}+e_{k,i}(0)(1-W_{k,i})+(\tau_{k,m}-\widehat\tau_k^{\,{\rm fixed}})(W_{k,i}-\mwidebar W_{k,m})\\
&-(\widehat\alpha_{k,m}-\alpha_{k,m}) - (\widehat\tau_{k,m}-\tau_{k,m})\mwidebar W_{k,m}\\
&=e_{k,i}(1)W_{k,i}+e_{k,i}(0)(1-W_{k,i})+(\tau_{k,m}-\tau_k)(W_{k,i}-\mwidebar W_{k,m})\\
&-(\widehat\tau_k^{\,{\rm fixed}}-\tau_k)(W_{k,i}-\mwidebar W_{k,m})-(\widehat\alpha_{k,m}-\alpha_{k,m}) - (\widehat\tau_{k,m}-\tau_{k,m})\mwidebar W_{k,m}.
\end{align*}
For $m_{k,i}=m$ and $N_{k,m,0}, N_{k,m,1}\geq 1$, let
\[
\mwidebar U_{k,i} = e_{k,i}(1)W_{k,i}+e_{k,i}(0)(1-W_{k,i})+(\tau_{k,m}-\tau_k)(W_{k,i}-\mwidebar W_{k,m}), 
\]
and let $\mwidebar U_{k,i} = 0$ for $m_{k,i}=m$ and $N_{k,m,0}N_{k,m,1}=0$.
Then, for $m_{k,i}=m$ and $N_{k,m,0}N_{k,m,1}\geq 1$, we have
\[
\widetilde U_{k,i} - \mwidebar U_{k,i} = -(\widehat\tau_k^{\,{\rm fixed}}-\tau_k)(W_{k,i}-\mwidebar W_{k,m})-(\widehat\alpha_{k,m}-\alpha_{k,m}) - (\widehat\tau_{k,m}-\tau_{k,m})\mwidebar W_{k,m}.
\]
Then,
\begin{align*}
\Bigg(&\sum_{i=1}^{n_k}1\{m_{k,i}=m\}R_{k,i}\widetilde W_{k,i}\widetilde U_{k,i}\Bigg)^2\\
&= \Bigg(\sum_{i=1}^{n_k}1\{m_{k,i}=m\}R_{k,i}\widetilde W_{k,i}\Big(\mwidebar U_{k,i}+\big(\widetilde U_{k,i}-\mwidebar U_{k,i}\big)\Big)\Bigg)^2\\
&=\Bigg(\sum_{i=1}^{n_k}1\{m_{k,i}=m\}R_{k,i}\widetilde W_{k,i}\Big(\mwidebar U_{k,i} -(\widehat\tau_k^{\,{\rm fixed}}-\tau_k)(W_{k,i}-\mwidebar W_{k,m})\Big)\Bigg)^2\\
&=\Bigg(\sum_{i=1}^{n_k}1\{m_{k,i}=m\}R_{k,i}\widetilde W_{k,i}\mwidebar U_{k,i} -(\widehat\tau_k^{\,{\rm fixed}}-\tau_k)\sum_{i=1}^{n_k}1\{m_{k,i}=m\}R_{k,i}W_{k,i}(W_{k,i}-\mwidebar W_{k,m})\Bigg)^2.
\end{align*}
Using the formula for the second moment of a binomial distribution and $n\geq 1$, we obtain,
\begin{align*}
E\Bigg[\Bigg(\sum_{i=1}^{n_k}1\{m_{k,i}=m\}R_{k,i}&W_{k,i}(W_{k,i}-\mwidebar W_{k,m})\Bigg)^2\Big|\mwidebar N_{k,m}=n\Bigg]\\
&=E\Bigg[\Bigg(\sum_{i=1}^{n_k}1\{m_{k,i}=m\}(1-\mwidebar W_{k,m})R_{k,i}W_{k,i}\Bigg)^2\Big|\mwidebar N_{k,m}=n\Bigg]\\
&\leq E\Bigg[\Bigg(\sum_{i=1}^{n_k}1\{m_{k,i}=m\}R_{k,i}W_{k,i}\Bigg)^2\Big|\mwidebar N_{k,m}=n\Bigg]\\
&\leq n^2 + n.
\end{align*}
From the formula of the sum of the first two moments of a binomial distribution, we obtain
\begin{align*}
\sum_{m=1}^{m_k}E\Bigg[\Bigg(\sum_{i=1}^{n_k}1\{m_{k,i}=m\}R_{k,i}W_{k,i}(W_{k,i}-\mwidebar W_{k,m})\Bigg)^2\Bigg]\leq \sum_{m=1}^{m_k} (n_{k,m}^2p_k^2q_k
+ n_{k,m}p_kq_k).
\end{align*}
Therefore, 
\begin{align*}
\frac{1}{n_kp_kq_kf_k}(\widehat\tau_k^{\,{\rm fixed}}-\tau_k)^2&\sum_{m=1}^{m_k}E\Bigg[\Bigg(\sum_{i=1}^{n_k}1\{m_{k,i}=m\}R_{k,i}W_{k,i}(W_{k,i}-\mwidebar W_{k,m})\Bigg)^2\Bigg]\\
&\leq \frac{n_kp_kq_k}{f_k} (\widehat\tau_k^{\,{\rm fixed}}-\tau_k)^2
\frac{1}{(n_kp_kq_k)^2} \sum_{m=1}^{m_k} (n_{k,m}^2p_k^2q_k
+ n_{k,m}p_kq_k)\\
&=\mathcal O_p(1) \left(
\frac{\max_m n_{k,m}}{\min_m n_{k,m}}\frac{1}{m_kq_k}+ \frac{1}{n_kp_kq_k}\right) \stackrel{p}{\longrightarrow} 0.
\end{align*}
Now, notice that 
\begin{align*}
\frac{1}{n_kp_kq_kf_k}&\sum_{m=1}^{m_k}\left(\sum_{i=1}^{n_k}1\{m_{k,i}=m\}R_{k,i}\widetilde W_{k,i}\mwidebar U_{k,i}\right)^2\\
&=\frac{1}{n_kp_kq_kf_k}\sum_{m=1}^{m_k}\Bigg(\sum_{i=1}^{n_k}1\{m_{k,i}=m\}R_{k,i}\big(e_{k,i}(1)W_{k,i}+e_{k,i}(0)(1-W_{k,i})\big)(W_{k,i}-\mwidebar W_{k,m})\\
&\hspace*{6cm}+(\tau_{k,m}-\tau_k)\sum_{i=1}^{n_k}1\{m_{k,i}=m\}R_{k,i}(W_{k,i}-\mwidebar W_{k,m})^2\Bigg)^2.
\end{align*}
Equation \eqref{equation:sqsume1} (and the analogous result for the sum involving terms with $e_{k,i}(0)$), implies
\[
\frac{1}{n_kp_kq_kf_k}\sum_{m=1}^{m_k}\Bigg(\sum_{i=1}^{n_k}1\{m_{k,i}=m\}R_{k,i}\big(e_{k,i}(1)W_{k,i}+e_{k,i}(0)(1-W_{k,i})\big)(\mwidebar W_{k,m}-A_{k,m})\Bigg)^2\stackrel{p}{\longrightarrow} 0.
\]
As a result, it is enough to establish convergence of $\mwidebar\Sigma_k/(n_kp_kq_kf_k)$, where
\begin{align*}
\mwidebar\Sigma_k &=\sum_{m=1}^{m_k}\Bigg(\sum_{i=1}^{n_k}1\{m_{k,i}=m\}R_{k,i}\big(e_{k,i}(1)W_{k,i}+e_{k,i}(0)(1-W_{k,i})\big)(W_{k,i}-A_{k,m})\\&
+(\tau_{k,m}-\tau_k)\sum_{i=1}^{n_k}1\{m_{k,i}=m\}R_{k,i}(W_{k,i}-\mwidebar W_{k,m})^2\Bigg)^2\\
&=\sum_{m=1}^{m_k}\Bigg(\sum_{i=1}^{n_k}1\{m_{k,i}=m\}\Big(R_{k,i}W_{k,i}(W_{k,i}-A_{k,m})- p_kq_kA_{k,m}(1-A_{k,m})\Big)e_{k,i}(1)\\
&+\sum_{i=1}^{n_k}1\{m_{k,i}=m\}\Big(R_{k,i}(1-W_{k,i})(W_{k,i}-A_{k,m})+ p_kq_kA_{k,m}(1-A_{k,m})\Big)e_{k,i}(0)\\
&+(\tau_{k,m}-\tau_k)\sum_{i=1}^{n_k}1\{m_{k,i}=m\}R_{k,i}(W_{k,i}-\mwidebar W_{k,m})^2\Bigg)^2.
\end{align*}
We will next show that
\begin{equation}
\label{equation:convsigmabar}
\frac{1}{n_kp_kq_kf_k}\mwidebar\Sigma_k-\frac{f_k^{\rm cluster}}{f_k}\stackrel{p}{\longrightarrow} 0,
\end{equation}
where
\begin{align*}
f_k^{\rm cluster}&=\frac{1}{n_k}E[A_{k,m}(1-A_{k,m})^2]\sum_{i=1}^{n_k} e_{k,i}^2(1)\\
&+\frac{1}{n_k}E[A^2_{k,m}(1-A_{k,m})]\sum_{i=1}^{n_k} e_{k,i}^2(0)\\
&- \frac{1}{n_k}p_kE[A^2_{k,m}(1-A_{k,m})^2]\sum_{i=1}^{n_k}(e_{k,i}(1)-e_{k,i}(0))^2\\
&+ (E[A_{k,m}(1-A_{k,m})]-(5+p_k)E[A^2_{k,m}(1-A_{k,m})^2])\sum_{m=1}^{m_k}\frac{n_{k,m}}{n_k}(\tau_{k,m}-\tau_k)^2\\
&+p_kE[A^2_{k,m}(1-A_{k,m})^2]\sum_{m=1}^{m_k}\frac{n^2_{k,m}}{n_k}(\tau_{k,m}-\tau_k)^2.
\end{align*}
Let
\begin{align*}
X_{k,m}&=\frac{1}{n_kp_kq_k}\Bigg(\sum_{i=1}^{n_k}1\{m_{k,i}=m\}R_{k,i}\big(e_{k,i}(0)(1-W_{k,i})\\&+e_{k,i}(1)W_{k,i}\big)(W_{k,i}-A_{k,m})
+(\tau_{k,m}-\tau_k)\sum_{i=1}^{n_k}1\{m_{k,i}=m\}R_{k,i}(W_{k,i}-\mwidebar W_{k,m})^2\Bigg)^2
\end{align*}
Using the result in equation \eqref{equation:psi2} and results on the moments of the binomial distribution (see intermediate calculations in section \ref{section:interfe}), we obtain 
\begin{align*}
\frac{1}{n_kp_kq_k}E[\mwidebar\Sigma_k] &=\sum_{m=1}^{m_k} E[X_{k,m}]\\&=f_k^{\rm cluster} +  \scaleto{\mathcal{O}}{5pt}(1).
\end{align*}
Therefore, to show that equation (\ref{equation:convsigmabar}) holds, we will show
\begin{equation}
\label{equation:conv4m}
\frac{1}{f_k^2}\sum_{m=1}^{m_k}E[X_{k,m}^2]\longrightarrow 0.
\end{equation}
Let
\begin{align*}
\theta_k&=E[(R_{k,i}W_{k,i}(W_{k,i}-A_{k,m})-p_kA_{k,m}(1-A_{k,m}))^2|m_{k,i}=m,Q_{k,m}=1]\\
&=p_k\left(E[A_{k,m}(1-A_{k,m})^2]-p_kE[A_{k,m}^2(1-A_{k,m})^2]\right),
\end{align*}
and
\begin{align*}
\pi_k&=E[(R_{k,i}W_{k,i}(W_{k,i}-A_{k,m})-p_kA_{k,m}(1-A_{k,m}))^4|m_{k,i}=m,Q_{k,m}=1]\\
&=p_k E[(W_{k,i}(W_{k,i}-A_{k,m})-p_kA_{k,m}(1-A_{k,m}))^4|m_{k,i}=m]+p_k^4(1-p_k)
E[A^4_{k,m}(1-A_{k,m})^4].
\end{align*}
Let
\begin{align*}
X_{k,m,1}&=\frac{1}{n_kp_kq_k}\Bigg(\sum_{i=1}^{n_k}1\{m_{k,i}=m\}R_{k,i}W_{k,i}(W_{k,i}-A_{k,m})e_{k,i}(1)\Bigg)^2\\
&=\frac{Q_{k,m}}{n_kp_kq_k}\Bigg(\sum_{i=1}^{n_k}1\{m_{k,i}=m\}(R_{k,i}W_{k,i}(W_{k,i}-A_{k,m})-p_kA_{k,m}(1-A_{k,m}))e_{k,i}(1)\Bigg)^2.
\end{align*}
Then,
\begin{align*}
E[X_{k,m,1}^2]&=q_kE[X_{k,m,1}^2|Q_{k,m}=1]\\
&=\frac{\pi_k}{n^2_kp^2_kq_k}\sum_{i=1}^{n_k}1\{m_{k,i}=m\} e_{k,i}^4(1)\\
&+\frac{6\theta^2_k}{n^2_kp^2_kq_k}\sum_{i=1}^{n_k-1}\sum_{j=i+1}^{n_k}1\{m_{k,i}=m_{k,j}=m\} e_{k,i}^2(1)e_{k,j}^2(1) .
\end{align*}
Therefore, because $n_kp_kq_k\rightarrow\infty$ and $m_kq_k\rightarrow\infty$, we obtain
\begin{align}
\label{equation:x1conv}
\sum_{m=1}^{m_k}E[X_{k,m,1}^2]
&\leq \frac{c}{n_kp_kq_k}\left(\frac{1}{n_k}\sum_{i=1}^{n_k} e_{k,i}^4(1)\right)
+\frac{c}{m_kq_k}\frac{\max_m n^2_{k,m}}{\min_m n_{k,m}^2}\nonumber\\&\hspace*{1cm}\times\left(\frac{1}{m_k}\sum_{m=1}^{m_k}\frac{1}{\max_m n^2_{k,m}}\sum_{i=1}^{n_k-1}\sum_{j=i+1}^{n_k}1\{m_{k,i}=m_{k,j}=m\} e_{k,i}^2(1)e_{k,j}^2(1)\right)\nonumber\\
&\longrightarrow 0.
\end{align}
Using the same argument, we obtain
\begin{equation}
\label{equation:x2conv}
\sum_{m=1}^{m_k}E[X_{k,m,2}^2]\longrightarrow 0,
\end{equation}
where 
\[
X_{k,m,2}=\frac{1}{n_kp_kq_k}\Bigg(\sum_{i=1}^{n_k}1\{m_{k,i}=m\}R_{k,i}(1-W_{k,i})(W_{k,i}-A_{k,m})e_{k,i}(0)\Bigg)^2.
\]
Notice that equations (\ref{equation:x1conv}) and (\ref{equation:x2conv}) imply
\[
\frac{1}{f_k^2}\sum_{m=1}^{m_k}E[X_{k,m,1}^2]\longrightarrow 0
\]
and
\[
\frac{1}{f_k^2}\sum_{m=1}^{m_k}E[X_{k,m,2}^2]\longrightarrow 0.
\]
Notice that the last two equations hold even if $f_k$ is bounded (e.g., when $\tau_{k,m}-\tau_k=0$ for all $k$ and $m$), as long as $f_k$ is bounded away from zero in large samples. In section \ref{section:fecte} we derive conditions so that  $f_k$ is bounded away from zero in large samples even if $\tau_{k,m}-\tau_k=0$ for all $k$ and $m$. 
Now, let 
\[
X_{k,m,3}=\frac{1}{n_kp_kq_k}\Bigg((\tau_{k,m}-\tau_k)\sum_{i=1}^{n_k}1\{m_{k,i}=m\}R_{k,i}W_{k,i}(W_{k,i}-\mwidebar W_{k,m})\Bigg)^2.
\]
Recall that, under the conditions in \eqref{equation:chte} and \eqref{equation:varAq}, $f_k\rightarrow\infty$ and $p_k\min n_{k,m}/f_k$ is bounded for large $k$ and, therefore, $p_k\max n_{k,m}/f_k$ is bounded for large $k$. Then (see intermediate calculations at the end of this document), for large $k$,
\begin{align*}
\frac{1}{f_k^2}\sum_{m=1}^{m_k} E[X_{k,m,3}^2]&= \frac{1}{(n_kp_kq_kf_k)^2}\sum_{m=1}^{m_k}n_{k,m}^4p_k^4q_k(\tau_{k,m}-\tau_k)^4\left(1+\mathcal O\left(\frac{1}{p_k\min_{m}n_{k,m}}\right)\right)\\
&= \frac{p_k\max_m n^2_{k,m}}{m_kq_kf_k\min_m n_{k,m}}\left(\frac{p_k}{f_k}\sum_{m=1}^{m_k}\frac{n_{k,m}^2}{n_k}(\tau_{k,m}-\tau_k)^4\right)\left(1+\mathcal O\left(\frac{1}{p_k\min_{m}n_{k,m}}\right)\right)\\
&=\mathcal O\left(\frac{1}{m_kq_k}\right)\left(1+\mathcal O\left(\frac{1}{p_k\min_{m}n_{k,m}}\right)\right)\rightarrow 0.
\end{align*}
Now, H\"{o}lder's inequality implies that equation (\ref{equation:conv4m}) holds (see intermediate calculations).

Now let,
\[
\tilde v_k^{\rm cluster} = f_k^{\rm cluster}/(\mu_k(1-\mu_k)-\sigma^2_k)^2.
\]
We obtain,
\[
\frac{\widetilde V_k^{\rm{cluster}}}{\tilde v_k}=\frac{\tilde v_k^{\rm cluster}}{\tilde v_k}+\scaleto{\mathcal{O}}{5pt}_p(1).
\]
We will next establish the analogous result for the heteroskedaticity-robust variance estimator.
Let
\begin{align*}
\widetilde\Sigma_{k}^{\rm robust}&=\sum_{i=1}^{n_k}R_{k,i}\widetilde W_{k,i}^2\widetilde U_{k,i}^2.
\end{align*}
Then, the heteroskedasticity-robust estimator of the variance of 
$\sqrt{N_k}(\widehat\tau^{\,{\rm fixed}}_k-\tau_k)$ is 
\begin{align*}
\widetilde V_k^{\rm{robust}} = N_k\widetilde Q_k^{-1}\widetilde\Sigma_k^{\rm robust} \widetilde Q_k^{-1}.
\end{align*}
As we have established before,
\begin{align*}
\widetilde U_{k,i} & =e_{k,i}(1)W_{k,i}+e_{k,i}(0)(1-W_{k,i})+(\tau_{k,m}-\tau_k)(W_{k,i}-\mwidebar W_{k,m})\\
&-(\widehat\tau_k^{\,{\rm fixed}}-\tau_k)(W_{k,i}-\mwidebar W_{k,m})-(\widehat\alpha_{k,m}-\alpha_{k,m}) - (\widehat\tau_{k,m}-\tau_{k,m})\mwidebar W_{k,m}.
\end{align*}
For $m_{k,i}=m$ and $N_{k,m,0}N_{k,m,1}\geq 1$, let
\[
\mwidebar U_{k,i} = e_{k,i}(1)W_{k,i}+e_{k,i}(0)(1-W_{k,i})+(\tau_{k,m}-\tau_k)(W_{k,i}-\mwidebar W_{k,m}), 
\]
and let $\mwidebar U_{k,i} = 0$ for $m_{k,i}=m$ and $N_{k,m,0}N_{k,m,1}=0$.
Then, for $m_{k,i}=m$ and $N_{k,m,0}N_{k,m,1}\geq 1$, we have
\[
\widetilde U_{k,i} - \mwidebar U_{k,i} = -(\widehat\tau_k^{\,{\rm fixed}}-\tau_k)(W_{k,i}-\mwidebar W_{k,m})-(\widehat\alpha_{k,m}-\alpha_{k,m}) - (\widehat\tau_{k,m}-\tau_{k,m})\mwidebar W_{k,m},
\]
and
\begin{equation}
\label{equation:Us}
\frac{1}{n_kp_kq_k}\sum_{i=1}^{n_k}R_{k,i}\widetilde W_{k,i}^2\widetilde U_{k,i}^2= \frac{1}{n_kp_kq_k}\sum_{i=1}^{n_k}R_{k,i}\widetilde W_{k,i}^2\Big(\mwidebar U_{k,i}+\big(\widetilde U_{k,i}-\mwidebar U_{k,i}\big)\Big)^2.
\end{equation}
Focusing on the part of the right hand side of last equation that depends on the first term of $\widetilde U_{k,i}-\mwidebar U_{k,i}$, we obtain 
\begin{align*}
\frac{1}{n_kp_kq_k}\sum_{i=1}^{n_k}R_{k,i}\widetilde W_{k,i}^4(\widehat\tau_k^{\,{\rm fixed}}-\tau_k)^2\leq (\widehat\tau_k^{\,{\rm fixed}}-\tau_k)^2\frac{1}{n_kp_kq_k}\sum_{i=1}^{n_k}R_{k,i}\widetilde W_{k,i}^2\stackrel{p}{\longrightarrow} 0.
\end{align*}
We will focus now on the part of the right-hand side of equation \eqref{equation:Us} that that depends on the second term of $\widetilde U_{k,i}-\mwidebar U_{k,i}$,
\begin{align*}
\frac{1}{n_kp_kq_k}&\sum_{m=1}^{m_k}\sum_{i=1}^{n_k}1\{m_{k,i}=m\}R_{k,i}\widetilde W_{k,i}^2(\widehat\alpha_{k,m}-\alpha_{k,m})^2.
\end{align*}
Using the formula for the variance of a sample mean under sampling without replacement \citep[e.g., in the supplement of][]{abadie2020sampling}, we obtain for $1\leq n\leq n_{k,m}-1$,
\begin{align}
\label{equation:bound_alpha}
E\Bigg[(\widehat\alpha_{k,m}-\alpha_{k,m})^2\sum_{i=1}^{n_k}&1\{m_{k,i}=m\}R_{k,i}\widetilde W_{k,i}^2\Big| N_{k,m,0}=n\Bigg]\nonumber\\
&=E\Bigg[(\widehat\alpha_{k,m}-\alpha_{k,m})^2\mwidebar N_{k,m}\mwidebar W_{k,m}(1-\mwidebar W_{k,m})\Big| N_{k,m,0}=n\Bigg]\nonumber\\
&\leq E\Big[n(\widehat\alpha_{k,m}-\alpha_{k,m})^2\big| N_{k,m,0}=n\Big]\nonumber\\
&=n\, \mbox{var}(\widehat\alpha_{k,m}|N_{k,m,0}=n)\nonumber\\
&=s^2_{k,m,0}\Big(1-\frac{n}{n_{k,m}}\Big),
\end{align}
where
\[
s^2_{k,m,0} = \frac{1}{n_{k,m}-1}\sum_{i=1}^{n_k}1\{m_{k,i}=m\}(y_{k,i}(0)-\alpha_{k,m})^2.
\]
Because $s^2_{k,m,0}$ is bounded, so is the right-hand side of equation \eqref{equation:bound_alpha}. As a result
\[
E\Bigg[\frac{1}{n_kp_kq_k}\sum_{m=1}^{m_k}\sum_{i=1}^{n_k}1\{m_{k,i}=m\}R_{k,i}\widetilde W_{k,i}^2(\widehat\alpha_{k,m}-\alpha_{k,m})^2\Bigg]\leq c\,\frac{m_k}{n_kp_k}\longrightarrow 0.
\]
An analogous derivation applies to the part of the right-hand side of equation \eqref{equation:Us} that depends on the third term of $\widetilde U_{k,i}-\mwidebar U_{k,i}$. (Notice that $\mwidebar W_{k,m}\leq 1$ and that $\widehat\tau_{k,m}-\tau_{k,m}$ is equal to minus the difference between $\widehat\alpha_{k,m}-\alpha_{k,m}$ and the analogous difference for the treated. 

Therefore, we will study the behavior of 
\begin{equation}
\label{equation:ubarsum}
\frac{1}{n_kp_kq_k}\sum_{i=1}^{n_k}R_{k,i}\widetilde W_{k,i}^2\mwidebar U_{k,i}^2.
\end{equation}
First, notice that 
\begin{align}
\label{equation:diffUs}
\frac{1}{n_kp_kq_k}&\sum_{m=1}^{m_k}\sum_{i=1}^{n_k} 1\{m_{k,i}=m\}R_{k,i} \Big|(W_{k,i}-\mwidebar W_{k,m})^2-(W_{k,i}-A_{k,m})^2\Big|W_{k,i} e^2_{k,i}(1)\nonumber\\
&\leq c\frac{1}{n_kp_kq_k}\sum_{m=1}^{m_k}\sum_{i=1}^{n_k} 1\{m_{k,i}=m\}R_{k,i} \Big|(W_{k,i}-\mwidebar W_{k,m})+(W_{k,i}-A_{k,m})\Big| \big|\mwidebar W_{k,m}-A_{k,m}\big|W_{k,i}\nonumber\\
&\leq c 
\Bigg(\frac{1}{n_kp_kq_k}\sum_{m=1}^{m_k}\sum_{i=1}^{n_k} 1\{m_{k,i}=m\}R_{k,i} \Big((W_{k,i}-\mwidebar W_{k,m})+(W_{k,i}-A_{k,m})\Big)^2W_{k,i} \Bigg)^{1/2}\nonumber\\
&\quad\times \Bigg(\frac{1}{n_kp_kq_k}\sum_{m=1}^{m_k}\sum_{i=1}^{n_k} 1\{m_{k,i}=m\}R_{k,i}\big(\mwidebar W_{k,m}-A_{k,m}\big)^2\Bigg)^{1/2}.
\end{align}
The inside of the first square root in equation \eqref{equation:diffUs} is bounded by a constant times
\[
\frac{1}{n_kp_kq_k}\sum_{m=1}^{m_k}\sum_{i=1}^{n_k} 1\{m_{k,i}=m\}R_{k,i},
\]
which converges in probability to one. The expectation of the inside of the second square root in equation \eqref{equation:diffUs} is
\begin{align*}
\frac{1}{n_kp_kq_k}\sum_{m=1}^{m_k}E\big[\mwidebar N_{k,m}\big(\mwidebar W_{k,m}-A_{k,m}\big)^2\big]\leq c \frac{m_k}{n_kp_k}\longrightarrow 0. 
\end{align*}
As a result, the right-hand side of equation \eqref{equation:diffUs} converges to zero in probability. The derivation with $(1-W_{k,i})e^2_{k,i}(0)$ replacing $W_{k,i}e^2_{k,i}(1)$ in equation \eqref{equation:diffUs} is analogous. Now, notice that
\begin{align*}
(W_{k,i}&-\mwidebar W_{k,m})^4-(W_{k,i}-\mwidebar A_{k,m})^4\\
&=-\big((W_{k,i}-\mwidebar W_{k,m})^2+(W_{k,i}- A_{k,m})^2\big)\big((W_{k,i}-\mwidebar W_{k,m})+(W_{k,i}- A_{k,m})\big)(\mwidebar W_{k,m}-A_{k,m}).
\end{align*}
Because the first factor of the expression above is bounded, we obtain
\begin{align}
\label{equation:diffUs2}
\frac{1}{n_kp_kq_k}&\sum_{m=1}^{m_k}\sum_{i=1}^{n_k} 1\{m_{k,i}=m\}R_{k,i} \Big|(W_{k,i}-\mwidebar W_{k,m})^4-(W_{k,i}-A_{k,m})^4\Big|(\tau_{k,m}-\tau_k)^2\nonumber\\
&\leq c 
\Bigg(\frac{1}{n_kp_kq_k}\sum_{m=1}^{m_k}\sum_{i=1}^{n_k} 1\{m_{k,i}=m\}R_{k,i} \Bigg)^{1/2}\nonumber\\
&\quad\times \Bigg(\frac{1}{n_kp_kq_k}\sum_{m=1}^{m_k}\sum_{i=1}^{n_k} 1\{m_{k,i}=m\}R_{k,i}\big(\mwidebar W_{k,m}-A_{k,m}\big)^2\Bigg)^{1/2}.
\end{align}
Now, the right-hand side of equation \eqref{equation:diffUs2} converges to zero in probability by the same argument as for equation \eqref{equation:diffUs}.
Cauchy-Schwarz inequality implies,
\begin{equation*}
\frac{1}{n_kp_kq_k}\sum_{i=1}^{n_k}R_{k,i}\widetilde W_{k,i}^2\mwidebar U_{k,i}^2=
\frac{1}{n_kp_kq_k}\sum_{i=1}^{n_k}1\{m_{k,i}=m\}R_{k,i}(W_{k,i}-A_{k,m})^2\breve U_{k,i}^2+
\scaleto{\mathcal{O}}{5pt}_p(1),
\end{equation*}
where
\begin{equation}
\label{equation:Ubreve}
\breve U_{k,i}=e_{k,i}(1)W_{k,i}+e_{k,i}(0)(1-W_{k,i})+(\tau_{k,m}-\tau_k)(W_{k,i}-A_{k,m}),
\end{equation}
for $m_{k,i}=m$ and $N_{k,m,0}N_{k,m,1}\geq 1$, and $\breve U_{k,i}=0$ for $N_{k,m,0}N_{k,m,1}\geq 0$.
Therefore, we will study the behavior of 
\[
\frac{1}{n_kp_kq_k}\sum_{m=1}^{m_k}\sum_{i=1}^{n_k} 1\{m_{k,i}=m\}R_{k,i} (W_{k,i}-A_{k,m})^2\breve U_{k,i}^2.
\]
We know,
\begin{align*}
\frac{1}{n_kp_kq_k}\sum_{m=1}^{m_k}\sum_{i=1}^{n_k} 1\{m_{k,i}=m\}R_{k,i}& (W_{k,i}-A_{k,m})^2W_{k,i}e^2_{k,i}(1)\\
& -
E[A_{k,m}(1-A_{k,m})^2]\frac{1}{n_k}\sum_{i=1}^{n_k} e^2_{k,i}(1)
\stackrel{p}{\longrightarrow} 0,
\end{align*}
and 
\begin{align*}
\frac{1}{n_kp_kq_k}\sum_{m=1}^{m_k}\sum_{i=1}^{n_k} 1\{m_{k,i}=m\}R_{k,i}& (W_{k,i}-A_{k,m})^2(1-W_{k,i})e^2_{k,i}(0)\\
\\&-
E[A^2_{k,m}(1-A_{k,m})]\frac{1}{n_k}\sum_{i=1}^{n_k} e^2_{k,i}(0)
\stackrel{p}{\longrightarrow} 0.
\end{align*}
Now, notice that
\begin{align*}
E[(W_{k,i}-A_{k,m})^4|m_{k,i}=m, R_{k,i}=1, A_{k,m}=a]&=(1-a)^4a+a^4(1-a)\\
&=a(1-a)[(1-a)^3+a^3]\\
&=a(1-a)[1-3a(1-a)],
\end{align*}
which implies
\begin{align*}
E\Bigg[\sum_{i=1}^{n_k} 1\{m_{k,i}=m\}R_{k,i}& (W_{k,i}-A_{k,m})^4(\tau_{k,m}-\tau_k)^2\Bigg]\\&=n_{k,m}p_kq_kE[A_{k,m}(1-A_{k,m})(1-3A_{k,m}(1-A_{k,m}))](\tau_{k,m}-\tau_k)^2,
\end{align*}
and
\begin{align*}
E\Bigg[\frac{1}{n_kp_kq_k}\sum_{m=1}^{m_k}&\sum_{i=1}^{n_k} 1\{m_{k,i}=m\}R_{k,i}(W_{k,i}-A_{k,m})^4(\tau_{k,m}-\tau_k)^2\Bigg]\nonumber\\
&=E[A_{k,m}(1-A_{k,m})(1-3A_{k,m}(1-A_{k,m}))]\sum_{m=1}^{m_k}\frac{n_{k,m}}{n_k}
(\tau_{k,m}-\tau_k)^2.
\end{align*}
Notice now that
\begin{align*}
\frac{1}{(n_kp_kq_k)^2}\sum_{m=1}^{m_k}&E\Bigg[\Bigg(\sum_{i=1}^{n_k} 1\{m_{k,i}=m\}R_{k,i}(W_{k,i}-A_{k,m})^4(\tau_{k,m}-\tau_k)^2\Bigg)^2\Bigg]\\
&\leq c\, \frac{1}{(n_kp_kq_k)^2}\sum_{m=1}^{m_k}E\Bigg[\Bigg(\sum_{i=1}^{n_k} 1\{m_{k,i}=m\}R_{k,i}\Bigg)^2\Bigg]\\
&\leq c\, \frac{q_k}{(n_kp_kq_k)^2}\sum_{m=1}^{m_k}(n_{k,m}p_k+n_{k,m}^2p_k^2)\\
&=c\Bigg(\frac{1}{n_kp_kq_k}+\frac{\max_m n_{k,m}}{\min_m n_{k,m}}\frac{1}{m_kq_k}\Bigg)
\stackrel{p}{\longrightarrow} 0.
\end{align*}

Notice also that expectations of the sums of products of the terms on the right-hand side of equation \eqref{equation:Ubreve} are equal to zero. Then, 
\begin{align*}
\frac{1}{n_kp_kq_k}\widetilde\Sigma_{k}^{\rm robust}- f_k^{\rm robust}\stackrel{p}{\longrightarrow} 0, 
\end{align*}
where
\begin{align*}
f_k^{\rm robust}&=E[A_{k,m}(1-A_{k,m})^2]\frac{1}{n_k}\sum_{i=1}^{n_k} e^2_{k,i}(1) + E[A^2_{k,m}(1-A_{k,m})]\frac{1}{n_k}\sum_{i=1}^{n_k} e^2_{k,i}(0)\\
&+E[A_{k,m}(1-A_{k,m})(1-3A_{k,m}(1-A_{k,m}))]\sum_{m=1}^{m_k}\frac{n_{k,m}}{n_k}
(\tau_{k,m}-\tau_k)^2.
\end{align*}
Now let,
\[
\tilde v_k^{\rm robust} = f_k^{\rm robust}/(\mu_k(1-\mu_k)-\sigma^2_k)^2.
\]
We obtain,
\[
\widetilde V_k^{\rm{robust}}=\tilde v_k^{\rm robust}+\scaleto{\mathcal{O}}{5pt}_p(1).
\]
\subsection{Large $k$ results the fixed effects case under homogeneous average treatment effects across clusters}
\label{section:fecte}

We will now study the Lyapounov's condition for the case $\tau_{k,m}=\tau_k$ for all $k$ and $m=1, \ldots, m_k$, so
\[
f_k = \sum_{m=1}^{m_k}E[\psi_{k,m}^2].
\]
Notice that 
\begin{align*}
\sum_{m=1}^{m_k}E[\psi_{k,m}^2]&\geq\frac{1}{n_k}E[A_{k,m}(1-A_{k,m})^2]\sum_{i=1}^{n_k} 1\{m_{k,i}=m\}e_{k,i}^2(1)\\
&+\frac{1}{n_k}E[A^2_{k,m}(1-A_{k,m})]\sum_{i=1}^{n_k} 1\{m_{k,i}=m\}e_{k,i}^2(0)\\
&- \frac{1}{n_k}E[A^2_{k,m}(1-A_{k,m})^2]\sum_{i=1}^{n_k} 1\{m_{k,i}=m\}(e_{k,i}(1)-e_{k,i}(0))^2\\
&=\frac{1}{n_k}E[A_{k,m}(1-A_{k,m})^3]\sum_{i=1}^{n_k} 1\{m_{k,i}=m\}e_{k,i}^2(1)\\
&+\frac{1}{n_k}E[A^3_{k,m}(1-A_{k,m})]\sum_{i=1}^{n_k} 1\{m_{k,i}=m\}e_{k,i}^2(0)\\
&+\frac{2}{n_k}E[A^2_{k,m}(1-A_{k,m})^2]\sum_{i=1}^{n_k} 1\{m_{k,i}=m\}e_{k,i}(1)e_{k,i}(0)\\
&=E\left[\frac{1}{n_k}\sum_{m=1}^{m_k}A^3_{k,m}(1-A_{k,m})^3\sum_{i=1}^{n_k} 1\{m_{k,i}=m\}\left(\frac{e_{k,i}(1)}{A_{k,m}}+\frac{e_{k,i}(0)}{1-A_{k,m}}\right)^2\right].
\end{align*}
Therefore,
\[
\liminf_{k\rightarrow\infty}E\left[\frac{1}{n_k}\sum_{m=1}^{m_k}A^3_{k,m}(1-A_{k,m})^3\sum_{i=1}^{n_k} 1\{m_{k,i}=m\}\left(\frac{e_{k,i}(1)}{A_{k,m}}+\frac{e_{k,i}(0)}{1-A_{k,m}}\right)^2\right]>0
\]
is sufficient for $\liminf_{k\rightarrow\infty} f_k>0$ (even if condition \eqref{equation:chte} does not hold). Given our assumption that the supports of the cluster probabilities, $A_{k,m}$, are bounded away from zero and one (uniformly in $k$ and $m$), then
\begin{equation}
\label{equation:boundforcte}
\liminf_{k\rightarrow\infty}E\left[\frac{1}{n_k}\sum_{m=1}^{m_k}\sum_{i=1}^{n_k} 1\{m_{k,i}=m\}\left(\frac{e_{k,i}(1)}{A_{k,m}}+\frac{e_{k,i}(0)}{1-A_{k,m}}\right)^2\right]>0
\end{equation}
is sufficient for $\liminf_{k\rightarrow\infty} f_k>0$. Assume that \eqref{equation:boundforcte} holds, so $\liminf_{k\rightarrow\infty} f_k>0$. We now obtain,
\begin{align*}
E\Bigg[\Bigg|\sum_{i=1}^{n_k}1\{m_{k,i}=m\}R_{k,i}&W_{k,i}(W_{k,i}-\mwidebar W_{k,m})e_{k,i}(1)\Bigg|^4\,\Big|\,Q_{k,m}=1, A_{k,m}\Bigg]\\
&=E\Bigg[(1-\mwidebar W_{k,m})^4\Bigg|\sum_{i=1}^{n_k}1\{m_{k,i}=m\}R_{k,i}W_{k,i}e_{k,i}(1)\Bigg|^4\,\Big|\,Q_{k,m}=1, A_{k,m}\Bigg]\\
&\leq E\Bigg[\Bigg|\sum_{i=1}^{n_k}1\{m_{k,i}=m\}R_{k,i}W_{k,i}e_{k,i}(1)\Bigg|^4\,\Big|\,Q_{k,m}=1, A_{k,m}\Bigg],
\end{align*}
and
\begin{align*}
E\Bigg[\Bigg|\sum_{i=1}^{n_k}1\{m_{k,i}=m\}R_{k,i}&W_{k,i}e_{k,i}(1)\Bigg|^4\,\Big|\,Q_{k,m}=1, A_{k,m}\Bigg]\\
&=E\Bigg[\Bigg|\sum_{i=1}^{n_k}1\{m_{k,i}=m\}(R_{k,i}W_{k,i}-p_kA_{k,m})e_{k,i}(1)\Bigg|^4\,\Big|\,Q_{k,m}=1, A_{k,m}\Bigg]\\
&=n_{k,m} E[(R_{k,i}W_{k,i}-p_kA_{k,m})^4|Q_{k,m}=1, A_{k,m}]\\&+3
n_{k,m}(n_{k,m}-1)(E[(R_{k,i}W_{k,i}-p_kA_{k,m})^2|Q_{k,m}=1, A_{k,m}])^2.
\end{align*}
The first equality holds because the terms $e_{k,i}(1)$ sum to zero within clusters. The second equality holds because, if $m_{k,i}=m_{k,j}=m$, with $i\neq j$, then $R_{k,i}W_{k,i}$ and $R_{k,i}W_{k,i}$ are independent conditional on $Q_{k,m}=1, A_{k,m}$, and $E[R_{k,i}W_{k,i}-p_kA_{k,m}|Q_{k,m}=1, A_{k,m}]=0$. Notice that 
\[
E[(R_{k,i}W_{k,i}-p_kA_{k,m})^2|Q_{k,m}=1, A_{k,m}]=p_kA_{k,m}(1-p_kA_{k,m})
\leq p_k,
\]
which also implies $E[(R_{k,i}W_{k,i}-p_kA_{k,m})^4|Q_{k,m}=1, A_{k,m}]
\leq p_k$. As a result,
\begin{align*}
\sum_{m=1}^{m_k}E\Bigg[\Bigg|\frac{1}{\sqrt{n_kp_kq_k}}\sum_{i=1}^{n_k}1\{m_{k,i}=m\}&R_{k,i}W_{k,i}e_{k,i}(1)\Bigg|^4\Bigg]\\&\leq \frac{1}{n_k p_k q_k}+ 3\,\frac{\max_m n_{k,m}}{\min_m n_{k,m}}\frac{1}{m_kq_k}\rightarrow 0.
\end{align*}

\section{Derivations of the variance estimators}
\label{section:new_variance}

In this section, we derive the adjustments in the CCV variance. (We do this under the assumption that the $Z_i$ are independent. In our simulations we actually use a slightly different sampling scheme for the $Z_i$ where the average $\overline{Z}_{k,m}$ is identical and fixed in each cluster.) To derive the CCV variance of the least squares estimator, consider first a variance estimator of the form
\[ 
\left(\sum_{i=1}^n V_i\right)^2.\]
We aim, however, to design an estimator based on a subsample consisting of units with $Z_i=1$, where $Z_i\in\{0,1\}$ is i.i.d.\ binary
with $\Pr(Z_i=1)=p_Z$ and independent of $V_i$.
First, notice that
\[
E\left[\left(\sum_{i=1}^n V_i\right)^2\right]=\sum_{i=1}^n E[V_i^2]+2\sum_{i=1}^{n-1}\sum_{j=i+1}^n E[V_iV_j],
\]
and
\[
E\left[\left(\sum_{i=1}^n Z_iV_i\right)^2\right]=p_Z\sum_{i=1}^n E[V_i^2]+2p_Z^2\sum_{i=1}^{n-1}\sum_{j=i+1}^n E[V_iV_j].
\]
Therefore,
\[
E\left[\frac{1}{p_Z}\left(\sum_{i=1}^n Z_iV_i\right)^2\right]=\sum_{i=1}^n E[V_i^2]+2p_Z\sum_{i=1}^{n-1}\sum_{j=i+1}^n E[V_iV_j],
\]
and
\[
\frac{(1-p_Z)}{p_Z^2}\left(E\left[\left(\sum_{i=1}^n Z_iV_i\right)^2\right]-p_Z\sum_{i=1}^n E[V_i^2]\right)=2(1-p_Z)\sum_{i=1}^{n-1}\sum_{j=i+1}^n E[V_iV_j].
\]
Adding the last two equations,
\begin{align}
E\left[\left(\sum_{i=1}^n V_i\right)^2\right]&=\frac{1}{p_Z^2}E\left[\left(\sum_{i=1}^n Z_iV_i\right)^2\right]-\frac{(1-p_Z)}{p_Z}\sum_{i=1}^n E[V_i^2]\nonumber\\
&=\frac{1}{p_Z^2}E\left[\left(\sum_{i=1}^n Z_iV_i\right)^2\right]-\frac{(1-p_Z)}{p_Z^2}\sum_{i=1}^n E[Z_iV_i^2].
\label{equation:deriv_ccv_ols}
\end{align}
The first term of the CCV variance estimator for least squares is based on the sample counterpart of the right-hand side of equation (\ref{equation:deriv_ccv_ols}), with $1\{m_{k,i}=m\}R_{k,i}((W_{k,i}-\mwidebar W_k)\widehat U_{k,i}-(\widehat\tau_{k,m}-\widehat\tau_k) \mwidebar W_k(1-\mwidebar W_k))$ in the role of $V_i$.

To derive the CCV variance estimator for the fixed effect case, consider
\[
\lambda_k = 1-q_k\frac{(E[A_{k,m}(1-A_{k,m})])^2}{E[A_{k,m}^2(1-A_{k,m})^2]},
\]
and let $f_k^{\rm CCV}=\lambda_k f_k^{\rm cluster}+(1-\lambda_k)f_k^{\rm robust}$. 
This transformation is designed to reproduce the terms in $f_k$ with factor
\[
\sum_{m=1}^{m_k}\frac{n^2_{k,m}}{n_k}
(\tau_{k,m}-\tau_k)^2.
\]
These terms dominate $f_k$ as $k$ increases. It also reproduces several lower order terms. 

Notice that 
\begin{align*}
f_k^{\rm robust}&=E[A_{k,m}(1-A_{k,m})^2]\frac{1}{n_k}\sum_{i=1}^{n_k} e^2_{k,i}(1) + E[A^2_{k,m}(1-A_{k,m})]\frac{1}{n_k}\sum_{i=1}^{n_k} e^2_{k,i}(0)\\
&+\Big(E[A_{k,m}(1-A_{k,m})]-(5+p_k)E[A^2_{k,m}(1-A_{k,m})^2]\Big)\sum_{m=1}^{m_k}\frac{n_{k,m}}{n_k}
(\tau_{k,m}-\tau_k)^2\\
&+(2+p_k)E[A^2_{k,m}(1-A_{k,m})^2]\sum_{m=1}^{m_k}\frac{n_{k,m}}{n_k}
(\tau_{k,m}-\tau_k)^2.
\end{align*}
Then,
\begin{align*}
    f_k^{\rm CCV}-f_k&=(1-\lambda_k)p_k E[A^2_{k,m}(1-A_{k,m})^2]\left(\sum_{m=1}^{m_k}\frac{n_{k,m}}{n_k}(\tau_{k,m}-\tau_k)^2+\frac{1}{n_k}\sum_{i=1}^{n_k}(e_{k,i}(1)-e_{k,i}(0))^2\right)\\
    &=p_kq_k(E[A_{k,m}(1-A_{k,m})])^2\left(\sum_{m=1}^{m_k}\frac{n_{k,m}}{n_k}(\tau_{k,m}-\tau_k)^2+\frac{1}{n_k}\sum_{i=1}^{n_k}(e_{k,i}(1)-e_{k,i}(0))^2\right).
\end{align*}
For $\tilde v_k^{\rm CCV}=f_k^{\rm CCV}/(\mu_k(1-\mu_k)-\sigma^2_k)^2$, we obtain,
\begin{equation}
\tilde v_k^{\rm CCV}-\tilde v_k = p_kq_k\sum_{m=1}^{m_k}\frac{n_{k,m}}{n_k}(\tau_{k,m}-\tau_k)^2+p_kq_k\frac{1}{n_k}\sum_{i=1}^{n_k}(e_{k,i}(1)-e_{k,i}(0))^2.
\label{equation:vkmvccv}
\end{equation}
The difference $\tilde v_k^{\rm CCV}-\tilde v_k$ is non-negative and of smaller order than $\tilde v_k$. Therefore, $\tilde v_k^{\rm CCV}/\tilde v_k\rightarrow 1$ (even if $\tilde v_k^{\rm CCV}-\tilde v_k$ is bounded away from zero). The first term on the right-hand side of (\ref{equation:vkmvccv}) could be estimated to further correct the difference between the CCV estimator and the variance of $\widehat\tau_k^{\rm fixed}$.  

\section{Limit results}
\label{section:Limit results}

Let $X_{k,m}$ be an infinite array of random variables, with rows indexed by $k=1, 2, \ldots$, and the columns of the $k$-th row indexed by $m=1, \ldots, m_k$. Let
\[
S_k =\sum_{m=1}^{m_k} X_{k,m},
\]
and  $a_k=E[S_k]$.

{\em A Weak Law of Large Numbers for Arrays:} For each $k=1,2,\ldots$ , suppose that $X_{k,1}, \ldots,$ $X_{k,m_k}$ are independent and have finite second moments. In addition, let $b_k$ be a sequence of positive constants such that 
\[
\frac{1}{b_k^2}\sum_{m=1}^{m_k} E[X_{k,m}^2]\longrightarrow 0.
\]
Then,
\[
\frac{S_k-a_k}{b_k}\stackrel{p}{\longrightarrow} 0.
\]
Proof: By Chebyshev's inequality, for any $\varepsilon>0$
\begin{align*}
\Pr\left(\left|\frac{S_k-a_k}{b_k}\right|>\varepsilon\right)&\leq \frac{1}{b_k^2\varepsilon^2}\mbox{var}(S_k)\\&=\frac{1}{b_k^2\varepsilon^2}\sum_{m=1}^{m_k}\mbox{var}(X_{k,m})\\
&\leq \frac{1}{b_k^2\varepsilon^2}\sum_{m=1}^{m_k}E[X_{k,m}^2]\longrightarrow 0.
\end{align*}
\hfill$\square$

{\em A Central Limit Theorem for Arrays:}  For each $k=1,2,\ldots$ , suppose that $X_{k,1}, \ldots,$ $X_{k,m_k}$ are independent, with zero means, $E[X_{k,m}]=0$, and finite variances,
$\sigma^2_{k,m}=E[X^2_{k,m}]$,  for $m=1,\ldots, m_k$. Let
\[
s_k^2 = \sum_{m=1}^{m_k} \sigma^2_{k,m}.
\]
Assume also that Lyapounov's condition holds, 
\[
\lim_{k\rightarrow\infty} \frac{1}{s_k^{2+\delta}}\sum_{m=1}^{m_k} E[|X_{k,m}|^{2+\delta}] = 0,
\]
for some $\delta>0$. Then,
\[
\frac{S_k}{s_k}\stackrel{d}{\longrightarrow} N(0,1).
\]
Proof: \cite{billingsley}, Chapter 27.

\section{Intermediate calculations for Section \ref{asection:base_case}}

The calculation of $v_k$ uses the following results.
\[
E[(R_{k,i}W_{k,i}-p_kq_k\mu_k)^2]=p_kq_k\mu_k(1-p_kq_k\mu_k),
\]
\[
E[(R_{k,i}(1-W_{k,i})-p_kq_k(1-\mu_k))^2]=p_kq_k(1-\mu_k)(1-p_kq_k(1-\mu_k)),
\]
\[
E[(R_{k,i}W_{k,i}-p_kq_k\mu_k)(R_{k,i}(1-W_{k,i})-p_kq_k(1-\mu_k))]= -p_k^2q_k^2\mu_k(1-\mu_k),
\]
\[
E[R_{k,i}W_{k,i}R_{k,j}W_{k,j}|m_{k,i}=m_{k,j}]=E[p_k^2q_kA_{k,m}^2]=p_k^2q_k(\sigma_k^2+\mu_k^2),
\]
and
\begin{align*}
E[(R_{k,i}W_{k,i}-p_kq_k\mu_k)(R_{k,j}W_{k,j}-p_kq_k\mu_k)|m_{k,i}=m_{k,j}]&=p_k^2q_k(\sigma_k^2+\mu_k^2)-(p_kq_k\mu_k)^2\\
&=p_k^2q_k(\sigma_k^2+(1-q_k)\mu_k^2).
\end{align*}
Similarly,
\begin{align*}
E[(R_{k,i}(1-W_{k,i})&-p_kq_k(1-\mu_k))(R_{k,j}(1-W_{k,j})-p_kq_k(1-\mu_k))|m_{k,i}
=m_{k,j}]\\
&=p_k^2q_k(\sigma_k^2+(1-q_k)(1-\mu_k)^2).
\end{align*}
Notice also that 
\begin{align*}
E[R_{k,i}W_{k,i}R_{k,j}(1-W_{k,j})|m_{k,i}=m_{k,j}]&=E[p_k^2q_kA_{k,m}(1-A_{k,m})]\\
&=p_k^2q_k(\mu_k(1-\mu_k)-\sigma_k^2),
\end{align*}
and
\begin{align*}
E[(R_{k,i}W_{k,i}-p_kq_k\mu_k)&(R_{k,j}(1-W_{k,j})-p_kq_k(1-\mu_k))|m_{k,i}=m_{k,j}]\\&=p_k^2q_k(\mu_k(1-\mu_k)-\sigma_k^2)-p_k^2q_k^2\mu_k(1-\mu_k)\\
&=p_k^2q_k(\mu_k(1-\mu_k)(1-q_k)-\sigma_k^2).
\end{align*}
\noindent\makebox[\linewidth]{\rule{\textwidth}{0.4pt}}
The following bounds are useful to prove Lyapunov's condition.
\begin{align*}
E[|R_{k,i}W_{k,i}-p_kq_k\mu_k|^3]&=(1-p_kq_k\mu_k)^3p_kq_k\mu_k+(p_kq_k\mu_k)^3(1-p_kq_k\mu_k)\\
&\leq c\, p_kq_k.
\end{align*}
Let $Q_{k,m}$ be a binary indicator that takes value one if cluster $m$ of population $k$ is sampled. 
\begin{align*}
E\big[&|R_{k,i}W_{k,i}-p_kq_k\mu_k|^2|R_{k,j}W_{k,j}-p_kq_k\mu_k|\big|m_{k,i}=m_{k,j}=m\big]\\
&=E\big[\big((1-p_kq_k\mu_k)^2p_kA_{k,m}+(p_kq_k\mu_k)^2(1-p_kA_{k,m})\big)\\
&\hspace*{2cm}\times\big((1-p_kq_k\mu_k)p_kA_{k,m}+(p_kq_k\mu_k)(1-p_kA_{k,m})\big)\big|m_{k,i}=m_{k,j}=m, Q_{k,m}=1\big]q_k\\
&+E\big[\big(p_kq_k\mu_k\big)^3\big|m_{k,i}=m_{k,j}=m, Q_{k,m}=0\big](1-q_k)\\
&\leq c p_k^2q_k.
\end{align*}
\begin{align*}
E\big[&|R_{k,i}W_{k,i}-p_kq_k\mu_k||R_{k,j}W_{k,j}-p_kq_k\mu_k||R_{k,t}W_{k,t}-p_kq_k\mu_k|\big|m_{k,i}=m_{k,j}=m_{k,t}=m\big]\\
&=E\big[\big((1-p_kq_k\mu_k)p_kA_{k,m}+(p_kq_k\mu_k)(1-p_kA_{k,m})\big)^3\big|m_{k,i}=m_{k,j}=m_{k,t}=m,Q_{k,m}=1\Big]q_k\\
&+E\big[\big(p_kq_k\mu_k\big)^3\big|m_{k,i}=m_{k,j}=m_{k,t}=m,Q_{k,m}=1\Big](1-q_k)\\
&\leq c p^3_kq_k.
\end{align*}
\noindent\makebox[\linewidth]{\rule{\textwidth}{0.4pt}}

Other useful intermediate calculations. 

For the moments of treatment indicators, notice that
$E[(W_{k,i}-\mu_k)^2W_{k,i}]=\mu_k(1-\mu_k)^2$, and $E[(W_{k,i}-\mu_k)^2(1-W_{k,i})]=(1-\mu_k)\mu_k^2$. In addition,
\begin{align*}
E[W_{k,i}W_{k,j}|m_{k,i}=m_{k,j}]&=E[A_{k,m}^2]\quad \mbox{(for $m\in\{1,\ldots, m_k\}$)}\\
& = \sigma_k^2+\mu_k^2.
\end{align*}
Similarly, $E[(1-W_{k,i})(1-W_{k,j})|m_{k,i}=m_{k,j}]=\sigma_k^2+(1-\mu_k)^2$. Therefore, $E[(W_{k,i}-\mu_k)W_{k,j}|m_{k,i}=m_{k,j}]=\sigma^2_k$ and $E[(W_{k,i}-\mu_k)(1-W_{k,j})|m_{k,i}=m_{k,j}]=-\sigma^2_k$. In addition, 
\begin{align*}
E[(W_{k,i}-\mu_k)(W_{k,j}&-\mu_k)W_{k,i}W_{k,j}|m_{k,i}=m_{k,j}]\\
&=E[A^2_{k,m}](1-\mu_k)^2\quad \mbox{(for $m\in\{1,\ldots, m_k\}$)}\\
&=(\sigma_k^2+\mu_k^2)(1-\mu_k)^2.
\end{align*}
Similarly,
\begin{align*}
E[(W_{k,i}-\mu_k)(W_{k,j}&-\mu_k)(1-W_{k,i})(1-W_{k,j})|m_{k,i}=m_{k,j}]=(\sigma_k^2+(1-\mu_k)^2)\mu_k^2,
\end{align*}
and
\begin{align*}
E[(W_{k,i}-\mu_k)(W_{k,j}-\mu_k)W_{k,i}(1-W_{k,j})|m_{k,i}=m_{k,j}]=\mu_k(1-\mu_k)(\sigma^2_k-\mu_k(1-\mu_k)).
\end{align*}

$\mbox{var}(R_{k,i}W_{k,i})=p_kq_k\mu_k(1-p_kq_k\mu_k)$, $\mbox{var}(R_{k,i}(1-W_{k,i}))=p_kq_k(1-\mu_k)(1-p_kq_k(1-\mu_k))$. Moreover,
\begin{align*}
\mbox{cov}(R_{k,i}W_{k,i},R_{k,i}(1-W_{k,i}))&=E[R_{k,i}W_{k,i}R_{k,i}(1-W_{k,i})]-E[R_{k,i}W_{k,i}]E[R_{k,i}(1-W_{k,i})]\\
&=-p_k^2q_k^2\mu_k(1-\mu_k).
\end{align*}
Recall that $E[W_{k,i}W_{k,j}|m_{k,i}=m_{k,j}] = \sigma_k^2+\mu_k^2$. Therefore, $\mbox{cov}(W_{k,i},W_{k,j}|m_{k,i}=m_{k,j})=\sigma_k^2$. Also,
\[
E[W_{k,i}(1-W_{k,j})|m_{k,i}=m_{k,j}]=\mu_k(1-\mu_k)-\sigma_k^2.
\]
\begin{align*}
E[R_{k,i}W_{k,i}R_{k,j}W_{k,j}|m_{k,i}=m_{k,j}]&=E[R_{k,i}R_{k,j}|m_{k,i}=m_{k,j}]E[W_{k,i}W_{k,j}|m_{k,i}=m_{k,j}]\\
&=p_k^2q_k(\sigma_k^2+\mu_k^2).
\end{align*}
Similarly,
\[
E[R_{k,i}(1-W_{k,i})R_{k,j}(1-W_{k,j})|m_{k,i}=m_{k,j}]=p_k^2q_k(\sigma_k^2+(1-\mu_k)^2).
\]
Therefore, 
\begin{align*}
\mbox{cov}(R_{k,i}W_{k,i},R_{k,j}W_{k,j}|m_{k,i}=m_{k,j})&=p_k^2q_k(\sigma_k^2+\mu_k^2)
-p_k^2q_k^2\mu_k^2\\
&=p_k^2q_k(\sigma_k^2+\mu_k^2(1-q_k)),\\ \shortintertext{and}
\mbox{cov}(R_{k,i}(1-W_{k,i}),R_{k,j}(1-W_{k,j})|m_{k,i}=m_{k,j})&=p_k^2q_k(\sigma_k^2+(1-\mu_k)^2)-p^2_kq^2_k(1-\mu_k)^2\\&=p_k^2q_k(\sigma_k^2+(1-\mu_k)^2(1-q_k)).
\end{align*}

In addition,
\begin{align*}
\mbox{cov}(R_{k,i}W_{k,i},R_{k,j}(1-W_{k,j})|m_{k,i}=m_{k,j})&=E[R_{k,i}W_{k,i}R_{k,j}(1-W_{k,j})|m_{k,i}=m_{k,j}]\\&-E[R_{k,i}W_{k,i}|m_{k,i}=m_{k,j}]E[R_{k,j}(1-W_{k,j})|m_{k,i}=m_{k,j}]\\
&=E[R_{k,i}R_{k,j}|m_{k,i}=m_{k,j}]E[W_{k,i}(1-W_{k,j})|m_{k,i}=m_{k,j}]\\&-E[R_{k,i}W_{k,i}|m_{k,i}=m_{k,j}]E[R_{k,j}(1-W_{k,j})|m_{k,i}=m_{k,j}]\\
&=p_k^2q_k(\mu_k(1-\mu_k)-\sigma_k^2)
-p_k^2q_k^2\mu_k(1-\mu_k)\\
&=p_k^2q_k(\mu_k(1-\mu_k)(1-q_k)-\sigma_k^2).
\end{align*}

\section{Intermediate calculations for Section \ref{asection:fixed_effects}}
\label{section:interfe}

\noindent\makebox[\linewidth]{\rule{\textwidth}{0.4pt}}
\[
E[R_{k,i}W_{k,i}(W_{k,i}-A_{k,m})|A_{k,m}, Q_{k,m}=1,m_{k,i}=m]=p_kA_{k,m}(1-A_{k,m}).
\]
This implies
\[
E[R_{k,i}W_{k,i}(W_{k,i}-A_{k,m})|m_{k,i}=m]=p_kq_kE[A_{k,m}(1-A_{k,m})].
\]
Therefore,
\[
E\Bigg[\sum_{i=1}^{n_k} 1\{m_{k,i}=m\}R_{k,i}W_{k,i}(W_{k,i}- A_{k,m})\Bigg]=n_{k,m}p_kq_kE[A_{k,m}(1-A_{k,m})].
\]

\noindent\makebox[\linewidth]{\rule{\textwidth}{0.4pt}}
For $n\geq 1$,
\begin{align*}
E\Bigg[\sum_{i=1}^{n_k} &1\{m_{k,i}=m\}R_{k,i}W_{k,i}(\mwidebar W_{k,m}-A_{k,m})\Big| \mwidebar N_{k,m}=n\Bigg]\\
&=\frac{1}{n}E\Bigg[\sum_{i=1}^{n_k} 1\{m_{k,i}=m\}R_{k,i}W_{k,i}\Bigg(\sum_{i=1}^{n_k} 1\{m_{k,i}=m\}R_{k,i}W_{k,i}-nA_{k,m}\Bigg)\Big| \mwidebar N_{k,m}=n\Bigg]\\
&=E[A_{k,m}(1-A_{k,m})].
\end{align*}
Therefore,
\begin{align*}
E\Bigg[\sum_{i=1}^{n_k} 1\{m_{k,i}=m\}R_{k,i}W_{k,i}(\mwidebar W_{k,m}-A_{k,m})\Bigg]&=E[A_{k,m}(1-A_{k,m})]\Pr(\mwidebar N_{k,m}\geq 1)\\
&=q_kE[A_{k,m}(1-A_{k,m})](1-(1-p_k)^{n_{k,m}}).
\end{align*}

\noindent\makebox[\linewidth]{\rule{\textwidth}{0.4pt}}

For $n\geq 1$
\begin{align*}
E[R_{k,i}W_{k,i}(\mwidebar W_{k,m}&-A_{k,m})^2|m_{k,i}=m, \mwidebar N_{k,m}=n, R_{k,i}=1]\\
&\leq E[(\mwidebar W_{k,m}-A_{k,m})^2|m_{k,i}=m, \mwidebar N_{k,m}=n, R_{k,i}=1]\\
&\leq \frac{E[A_{k,m}(1-A_{k,m})]}{n}.
\end{align*}
Because $\Pr(R_{k,i}=1|\mwidebar N_{k,m}=n, m_{k,i}=m)=n/n_{k,m}$, we obtain
\[
E[R_{k,i}W_{k,i}(\mwidebar W_{k,m}-A_{k,m})^2|m_{k,i}=m,\mwidebar N_{k,m}=n]\leq 
\frac{E[A_{k,m}(1-A_{k,m})]}{n_{k,m}},
\]
which implies
\[
E[R_{k,i}W_{k,i}(\mwidebar W_{k,m}-A_{k,m})^2|m_{k,i}=m,\mwidebar N_{k,m}\geq 1]\leq 
\frac{E[A_{k,m}(1-A_{k,m})]}{n_{k,m}}.
\]
Therefore, 
\begin{align*}
E[R_{k,i}W_{k,i}&(\mwidebar W_{k,m}-A_{k,m})^2|m_{k,i}=m]\\&=
E[R_{k,i}W_{k,i}(\mwidebar W_{k,m}-A_{k,m})^2|m_{k,i}=m, \mwidebar N_{k,m}\geq 1]\Pr(\mwidebar N_{k,m}\geq 1|m_{k,i}=m)\\
&\leq 
q_k\frac{E[A_{k,m}(1-A_{k,m})]}{n_{k,m}}.
\end{align*}

\noindent\makebox[\linewidth]{\rule{\textwidth}{0.4pt}}

Conditional on $\mwidebar N_{k,m}=n$ and $A_{k,m}$, the variable $N_{k,m,1}$ has a binomial distribution with parameters $(n,A_{k,m})$. Then, using the formulas for the moments of a binomial distribution, we find that for any integer $n$, such that $1\leq n\leq n_{k,m}$, 
\begin{align*}
E\Bigg[\Bigg(\sum_{i=1}^{n_k}&1\{m_{k,i}=m\}R_{k,i}W_{k,i}(W_{k,i}-\mwidebar W_{k,m})\Bigg)^2\Big|A_{k,m}=a, \mwidebar N_{k,m}=n\Bigg]\\ &= E[(N_{k,m,1}-N_{k,m,1}^2/n)^2|A_{k,m}=a, \mwidebar N_{k,m}=n]\\
&=n^2a^2(1-a)^2+na(1-a)(1-6a+6a^2)+r_1(a) + r_2(a)/n,
\end{align*}
where $|r_1(a)|$ and $|r_2(a)|$ are uniformly bounded in $a\in[0,1]$. Therefore, 
\begin{align*}
E\Bigg[\Bigg(\sum_{i=1}^{n_k}&1\{m_{k,i}=m\}R_{k,i}W_{k,i}(W_{k,i}-\mwidebar W_{k,m})\Bigg)^2\Big|\mwidebar N_{k,m}=n\Bigg]\\
&=n^2E[A_{k,m}^2(1-A_{k,m})^2]+nE[A_{k,m}(1-A_{k,m})(1-6A_{k,m}+6A_{k,m}^2)]\\&+E[r_1(A_{k,m})] + E[r_2(A_{k,m})]/n.
\end{align*}
It follows that
\begin{align*}
E\Bigg[\sum_{m=1}^{m_k}&(\tau_{k,m}-\tau_k)^2\Bigg(\sum_{i=1}^{n_k}1\{m_{k,i}=m\}R_{k,i}W_{k,i}(W_{k,i}-\mwidebar W_{k,m})\Bigg)^2\Bigg]\\ 
&=\Bigg(\sum_{m=1}^{m_k}(\tau_{k,m}-\tau_k)^2(n_{k,m}(n_{k,m}-1)p_k^2q_k+n_{k,m}p_kq_k)\Bigg)E[A_{k,m}^2(1-A_{k,m})^2]\\&+\sum_{m=1}^{m_k}(\tau_{k,m}-\tau_k)^2n_{k,m}p_kq_kE[A_{k,m}(1-A_{k,m})(1-6A_{k,m}(1-A_{k,m}))]+\mathcal O(m_kq_k).
\end{align*}
Therefore,
\begin{align*}
\frac{1}{n_kp_kq_k}E\Bigg[\sum_{m=1}^{m_k}&(\tau_{k,m}-\tau_k)^2\Bigg(\sum_{i=1}^{n_k}1\{m_{k,i}=m\}R_{k,i}W_{k,i}(W_{k,i}-\mwidebar W_{k,m})\Bigg)^2\Bigg]\\
\longrightarrow& (E[A_{k,m}(1-A_{k,m})]-(5+p_k)E[A^2_{k,m}(1-A_{k,m})^2])\sum_{m=1}^{m_k}\frac{n_{k,m}}{n_k}(\tau_{k,m}-\tau_k)^2\\
&+p_kE[A^2_{k,m}(1-A_{k,m})^2]\sum_{m=1}^{m_k}\frac{n^2_{k,m}}{n_k}(\tau_{k,m}-\tau_k)^2.
\end{align*}
\noindent\makebox[\linewidth]{\rule{\textwidth}{0.4pt}}
Notice that,
\begin{align*}
E\Bigg[\Bigg(\sum_{i=1}^{n_k}1\{m_{k,i}=m\}R_{k,i}W_{k,i}&(W_{k,i}-\mwidebar W_{k,m})\Bigg)^4\Big|A_{k,m}=a, \mwidebar N_{k,m}=n\Bigg]\\ &= E[(N_{k,m,1}(1-N_{k,m,1}/n))^4|A_{k,m}=a, \mwidebar N_{k,m}=n]\\
&\leq  E[N^4_{k,m,1}|A_{k,m}=a, \mwidebar N_{k,m}=n]\\
&\leq n^4,
\end{align*}
Therefore, 
\begin{align*}
E\Bigg[\Bigg(\sum_{i=1}^{n_k}1\{m_{k,i}=m\}R_{k,i}W_{k,i}&(W_{k,i}-\mwidebar W_{k,m})\Bigg)^4\Bigg]=n_{k,m}^4p_k^4q_k\left(1+\mathcal O\left(\frac{1}{p_k\min_{m}n_{k,m}}\right)\right),
\end{align*}
uniformly in $m$.

\noindent\makebox[\linewidth]{\rule{\textwidth}{0.4pt}}

Suppose $X_{k,m} = (Z_{k,m,1}+Z_{k,m,2})^2$. Let $X_{k,m,1}=Z_{k,m,1}^2$ and $X_{k,m,2}=Z_{k,m,2}^2$. Now suppose,
\[
\sum_{m=1}^{m_k} E[X_{k,m,1}^2]\longrightarrow 0,
\]
and
\[
\sum_{m=1}^{m_k} E[X_{k,m,2}^2]\longrightarrow 0.
\]
Using the binomial theorem and H\"{o}lder's inequality, we obtain
\begin{align*}
\sum_{m=1}^{m_k} E[X_{k,m}^2] & = \sum_{m=1}^{m_k}\sum_{p=0}^{4} c_p E[Z_{k,m,1}^pZ_{k,m,2}^{(4-p)}]\\
&\leq c \sum_{m=1}^{m_k}\sum_{p=0}^{4} E[|Z_{k,m,1}|^p|Z_{k,m,2}|^{(4-p)}]\\
&\leq c \sum_{m=1}^{m_k}\sum_{p=0}^{4} (E[X_{k,m,1}^2])^{p/4}(E[X_{k,m,2}^2])^{(4-p)/4}\\
&\leq c\sum_{p=0}^{4} \left(\sum_{m=1}^{m_k} E[X_{k,m,1}^2]\right)^{p/4}\left(\sum_{m=1}^{m_k}E[X_{k,m,2}^2]\right)^{(4-p)/4}\longrightarrow 0.
\end{align*}

\end{document}